\documentclass[11pt, twoside,reqno]{amsart}

\usepackage[top=1in,bottom=1in,left=1.2in,right=1.2in]{geometry}
\usepackage{amsmath,amssymb,amsthm,amsfonts,enumerate, mathrsfs, mathtools, bm}
\usepackage{appendix}
\usepackage{relsize}
\usepackage{stmaryrd}
\usepackage{esint}
\usepackage{xcolor}
\usepackage{graphicx}
\usepackage{needspace}
\usepackage[numbers]{natbib}

\usepackage[
colorlinks=true,
linkcolor=blue,
filecolor=blue,
urlcolor=red,
citecolor=magenta]{hyperref}
\usepackage[nameinlink]{cleveref}

\newcounter{counterConstant}

\numberwithin{equation}{section}
\def\theequation{\arabic{section}.\arabic{equation}}

\newtheorem{theorem}{Theorem}[section]
\newtheorem{lemma}[theorem]{Lemma}
\newtheorem{remark}[theorem]{Remark}
\newtheorem{definition}[theorem]{Definition}
\newtheorem{proposition}[theorem]{Proposition}
\newtheorem{corollary}[theorem]{Corollary}
\newtheorem{assumption}{Assumption}

\def\cA{{\mathcal A}}

\def\cG{{\mathcal G}}
\def\cH{{\mathcal H}}

\def\cK{{\mathcal K}}
\def\cL{{\mathcal L}}

\def\cN{{\mathcal N}}

\def\cP{{\mathcal P}}

\def\cR{{\mathcal R}}

\def\cW{{\mathcal W}}

\def\mD{{\mathbb D}}

\def\mN{{\mathbb N}}

\def\mP{{\mathbb P}}

\def\mR{{\mathbb R}}
\def\mS{{\mathbb S}}

\def\mZ{{\mathbb Z}}

\def\bE{{\mathbf E}}
\def\bF{{\mathbf F}}

\def\sF{{\mathscr F}}

\def\sH{{\mathscr H}}

\def\sK{{\mathscr K}}

\def\geq{\geqslant}
\def\leq{\leqslant}

\def\1{{\mathbf{1}}}

\def\d{\mbox{\rm{d}}}
\def\dd{\,\d}

\def\eps{\varepsilon}

\def\div{\mathord{{\rm div}\,}}

\def\esssup{\mathop{\mathrm{ess\,sup}}}

\def\tr{\mathrm{Tr}}
\def\Id{\mathrm{Id}}
\def\Law{{\mathord{{\rm Law}}}}

\usepackage{fancyhdr}
{} 

\allowdisplaybreaks

\begin{document}

\title{Strong solvability and its failure for critical SDEs}
\author{Dejun Luo and Guohuan Zhao}

\address{State Key Laboratory of Mathematical Sciences, Academy of Mathematics and Systems Science, CAS, Beijing, 100190, China}
\email{luodj@amss.ac.cn}

\address{State Key Laboratory of Mathematical Sciences, Academy of Mathematics and Systems Science, CAS, Beijing, 100190, China}
\email{gzhao@amss.ac.cn}

\begin{abstract}
We revisit the construction of strong solutions to singular SDEs by compactness on Wiener space. For spatially varying diffusion coefficients, assuming sufficiently small form bounds on the drift and diffusion gradient, we use basic energy estimates and an extended compactness criterion for Wiener functionals to construct an approximation-independent stochastic flow of homeomorphisms. Both forward and inverse maps have local Sobolev moment bounds, and the resulting strong solutions are pathwise unique in a class of solutions satisfying a Krylov estimate.

We also show that the endpoint critical  condition $b\in L^\infty_tL^d_x$, together with $\mathrm{div } b=0$ does not guarantee the existence of a strong solution. For every $d\geq2$, we construct a SDE with a divergence-free drift in this class and an additive Brownian noise, started from the origin, has a weak solution unique in law, but no strong solution. The proof uses logarithmic rescaling to characterize strong solvability by whether the sample measures of an associated stationary diffusion are Dirac masses. Positivity of the top Lyapunov exponent, combined with random SRB theory, rules out this concentration and hence strong solvability.
\end{abstract}

\maketitle
\enlargethispage{2pt}


\setcounter{tocdepth}{1}
\tableofcontents 

\section{Introduction}

We consider stochastic differential equations in $\mR^d$ on $[s,T]$, where $T>0$ and $0\leq s<T$, written in It\^o form as
\begin{equation}\label{eq:intro}
 \d X_t=\widehat b(t,X_t)\dd t+\sigma(t,X_t)\dd W_t,
 \quad X_s=x\in \mR^d.
\end{equation}
Here $W$ is a standard $d$-dimensional Brownian motion, and the coefficients $\widehat b:[0,T]\times\mR^d\to\mR^d$ and $\sigma:[0,T]\times\mR^d\to\mR^{d\times d}$.

The parabolic scaling of \eqref{eq:intro} identifies the natural critical spaces for the drift and diffusion. For $\lambda>0$, the rescaled coefficients are
\[
 \widehat b_\lambda(t,x):=\lambda\widehat b(\lambda^2t,\lambda x),
 \qquad \sigma_\lambda(t,x):=\sigma(\lambda^2t,\lambda x).
\]
For $1\leq p,r\leq\infty$, we write the Lebesgue spaces 
\[
 L^r_tL^p_x:=L^r(0,T;L^p(\mR^d)),\qquad
 C_tL^p_x:=C([0,T];L^p(\mR^d)).
\]
Under this scaling, $\|\sigma\|_{L^\infty_{t,x}}$ is invariant, while the $L^r_tL^p_x$ norms of $\widehat b$ and $\nabla\sigma$ both scale by $\lambda^{1-\frac{d}{p}-\frac{2}{r}}$ on the corresponding time intervals. The critical conditions are therefore
\[
 \widehat b,\nabla\sigma\in L^r_tL^p_x,\quad
 \sigma\in L^\infty_{t,x},\quad {d}/{p}+{2}/{r}=1.
\]
The endpoint of interest is $(p,r)=(d,\infty)$, corresponding to $\widehat b,\nabla\sigma\in L^\infty_tL^d_x$. 

Our main results in this manuscript address both strong solvability and its failure at this endpoint. In Part~\ref{part:well}, we simplify the Wiener-space compactness method used in \cite{RZ25sdes} and extend it to multiplicative noise. Under small form bounds on the drift and diffusion gradient (including all cases with $(p,r)\neq (d,\infty)$), relatively simple energy estimates yield the existence of strong solutions and a stochastic flow. In Part~\ref{part:ill}, we construct an equation with additive noise and a divergence-free drift in $L^\infty_tL^d_x$ that has a weak solution unique in law but no strong solution. 

\medskip 

Recently, Krylov \cite{K25weak,K25strong} studies strong solvability for equations with rough drift and diffusion coefficients, with uniqueness in classes of solutions satisfying Krylov-type estimates. In particular, \cite{K25strong} establishes strong existence under small spatial Morrey bounds on both the drift and the diffusion gradient. For additive noise, R\"ockner--Zhao \cite[Theorem~1.1]{RZ25sdes} proved strong existence for drift in $C_tL^d_x$ or $L^r_tL^p_x$ with $r\in (2,\infty)$, $d\geq3$. Their construction combines smooth approximation, parabolic estimates, and Malliavin compactness on the original Wiener space. Kinzebulatov--Madou \cite{KM24strong} extended this approach to small form-bounded drifts with additive noise, measuring the size of the drift through its quadratic form bound. The latter two works provide the starting point for Part~\ref{part:well}.

To implement this approach for spatially varying diffusion, we draw on \cite{BFGM19stochastic} and study the (inverse) stochastic flows generated by spatially smooth approximations of \eqref{eq:intro}. The corresponding stochastic transport equations yield a closed system for moments of spatial derivatives, which we control by weighted energy inequalities. Assuming merely measurable time dependence and small form bounds on the drift and diffusion gradient, we combine these estimates with an extended compactness criterion for Wiener functionals to obtain a stochastic flow with Sobolev bounds for both the forward and inverse maps.

During the preparation of this manuscript, we became aware of Boutiah--Kinzebulatov \cite{BK26strong}. Their main theorem establishes strong existence for driftless SDEs with spatially varying diffusion under a small form bound on $\nabla\sigma$ and an additional H\"older regularity assumption in time. Their construction also employs the Wiener-space compactness approach used in \cite{RZ25sdes} and \cite{KM24strong}, with derivative estimates obtained by combining simplex estimates with a stochastic Gronwall argument. 

\medskip

However, the strong existence results discussed above rely on smallness in suitable norms or quadratic form bounds, and do not cover the unrestricted endpoint class $L^\infty_tL^d_x$. This endpoint difficulty was emphasized by Beck--Flandoli--Gubinelli--Maurelli \cite[Section~1.4]{BFGM19stochastic}. In their study of stochastic transport and continuity equations and the associated stochastic flows, they distinguish $C_tL^d_x$ from the unrestricted $L^\infty_tL^d_x$ class and leave the latter case open. They obtain Sobolev regularity through parabolic equations for moments of spatial derivatives, then use duality to prove uniqueness and construct stochastic flows.

The divergence-free condition raises the possibility of removing the smallness assumption at the endpoint. Indeed, for divergence-free drifts in $L^\infty_tL^d_x$, weak existence and uniqueness in classes of solutions satisfying Krylov estimates were proved without a smallness assumption in \cite[Theorem~1.1(b)]{RZ23sdes} for $d\geq3$ and \cite[Theorem~1.1]{RZZ25mckeanvlasovb} for $d=2$. The two-dimensional result also covers the larger Lorentz class $L^\infty_tL^{2,\infty}_x$. Both works use regularity estimates for the backward Kolmogorov equation.

Nonuniqueness in law is known in several supercritical regimes. In the Brownian case, Galeati--Gerencs\'er \cite[Remark~1.11]{GG25solution} construct drifts in $L^r_tL^p_x$ with $p>d$ and $d/p+2/r>1$ for which uniqueness in law fails. For $d\geq2$, L\"u--R\"ockner \cite[Theorem~1.1]{LR25sharp} obtain nonuniqueness in law for divergence-free drifts in $L^r_tL^p_x\cap C_tL^s_x$, $1<s<d$, under the condition $d/p+1/r>1$. Neither result includes the endpoint $(p,r)=(d,\infty)$ considered here.

Examples with no strong solution were constructed by Hess-Childs--Rowan \cite{HR26sharp}. In their Theorem~1.3, they construct divergence-free drifts in $L^\infty_tC^{-\rho}_x\cap L^1_tC^0_x$, $0<\rho<1/2$, for which the additive SDE is weakly unique but not strongly solvable. Their construction uses random alternating shears at several scales and also produces deterministic examples. However, it does not establish an $L^\infty_tL^d_x$ endpoint example.

Our example concerns divergence-free drifts in the endpoint class $L^\infty_tL^d_x$, without a smallness assumption. For $d\geq3$, strong well-posedness for divergence-free drifts in this class was conjectured in \cite[Remark~1.2(ii)]{RZ25sdes}. We show that a strong-solution theorem for the entire class is impossible. For every $d\geq2$, we construct a smooth compactly supported vector field $B_K$, defined in Section~\ref{sec:profile}. This vector field describes planar rotation with angular speed varying with the radius. Viewed at the parabolic scale $\sqrt t$, this profile defines a deterministic divergence-free drift with no strong solution. Its spatial $L^d$ norm stays constant as $t\downarrow0$, while its spatial mass concentrates near the initial point. Nevertheless, a weak solution exists and is unique in law, with Gaussian one-time marginals. 

\begin{figure}[htbp]
 \centering
 \includegraphics[width=0.9\textwidth]{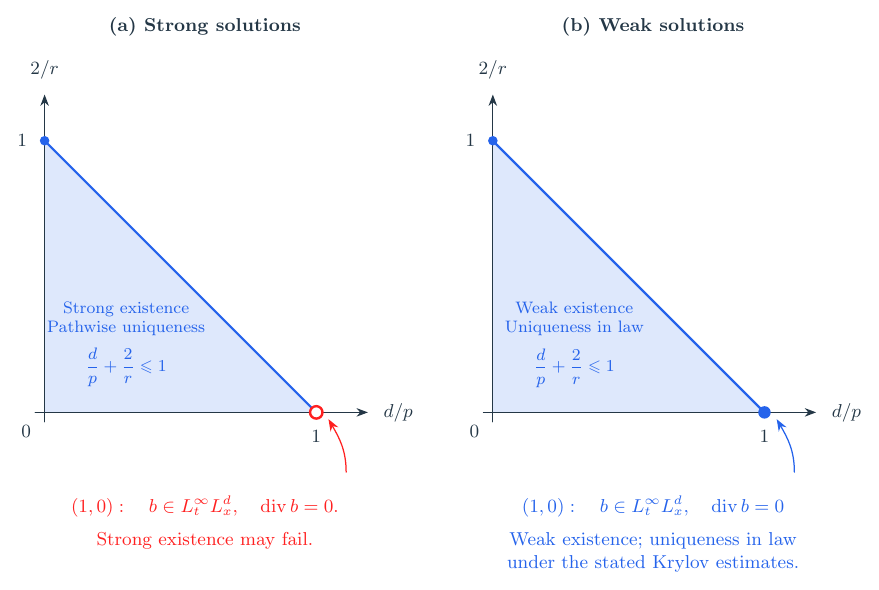}
\caption{\scriptsize Strong and weak solvability for $\sigma=I_d$, $d\geq2$, and $b\in L^r_tL^p_x$.}
 \label{fig:range}
\end{figure}

Figure~\ref{fig:range} compares strong and weak solvability in the $(d/p,2/r)$-plane when $\sigma=I_d$ and $\widehat{b}=b\in L^r_t L^p_x$.

\medskip

\subsection{Main results}\label{sec:results}

Before stating the main results, we specify the solution concepts for equation~\eqref{eq:intro} and the underlying probability spaces. Throughout, filtered probability spaces satisfy the usual conditions. A \emph{weak solution} with initial data $(s,x)$ is a tuple
\[
 (\Omega,\sF,(\sF_t)_{s\leq t\leq T},\mP;X,W) 
\]
satisfying the following properties. Under $\mP$, the process $W$ is continuous and adapted, and $W^s_t:=W_t-W_s$, $s\leq t\leq T$, is a standard $d$-dimensional Brownian motion started from zero at time $s$, with respect to $(\sF_t)$. The process $X$ is $(\sF_t)$-adapted with continuous paths on $[s,T]$. We also require
\[
 \int_s^T\bigl(|\widehat b(r,X_r)|+|\sigma(r,X_r)|^2\bigr)\dd r<\infty
 \quad\mP\mbox{-a.s.},
\]
and
\[
 X_t=x+\int_s^t\widehat b(r,X_r)\dd r
           +\int_s^t\sigma(r,X_r)\dd W_r
\]
for every $t\in[s,T]$ on an event of $\mP$-probability one.

A \emph{strong solution} is a weak solution for which $X$ is adapted to the filtration $\sigma\{W^s_r:s\leq r\leq t\}\vee\cN$, where $\cN$ denotes the family of $\mP$-null sets. 

\emph{Uniqueness in law} means that all weak solutions with the same initial data have the same distribution of $X$. On a space with probability measure $\mP$, we write $\Law(X):=\mP\circ X^{-1}$ for the path law and distinguish it from $\Law(X,W^s)$ which is  the joint law of $X$ and $W^s$. 

\emph{Pathwise uniqueness} means that two solutions $X,\widetilde X$ on the same filtered probability space, with the same Brownian motion and initial data, satisfy
\[
 \mP\bigl(X_t=\widetilde X_t\mbox{ for every }t\in[s,T]\bigr)=1.
\]
When uniqueness is restricted to a given class, both competing solutions are required to belong to that class. 

For an initial time $s<T$ and point $x\in \mR^d$, let $\sK_2(s,x)$ be the class of weak solutions of \eqref{eq:intro} with initial data $(s,x)$ such that, for every $a\in(s,T)$, there is a finite constant $C_{a,X}$ with
\begin{equation}\label{eq:k2}
 \bE_{\mP}\int_a^T |f(r,X_r)|\dd r
 \leq C_{a,X}\|f\|_{L^2((a,T)\times\mR^d)}, 
\end{equation}
for every Borel $f\in L^2((a,T)\times\mR^d)$. The constant $C_{a,X}$ is independent of $f$, may depend on the solution, and may blow up as $a\downarrow s$. We write $(X,W)\in\sK_2(s,x)$ for short.

\medskip

We work on the classical Wiener space
\[
 \Omega:=C([0,T];\mR^d),\quad \mP:=\cW_{0,T},
\]
where $\cW_{s,T}$ is the Wiener measure on $\{w\in C([s,T];\mR^d):w(s)=0\}$ and $W_t(\omega)=\omega(t)$ for $\omega\in\Omega$. $\sF:=\sigma(W_t:0\leq t\leq T)\vee \cN$, and $(\sF_t)$ is the natural filtration. In the rest of the paper, we use the Stratonovich form of \eqref{eq:intro}: for $b,\sigma_k:[0,T]\times\mR^d\to\mR^d$, $1\leq k\leq d$, consider
\begin{equation}\label{eq:sde}
 X_{s,t}(x)=x+\int_s^t b(r,X_{s,r}(x))\dd r
       +\sum_{k=1}^d\int_s^t\sigma_k(r,X_{s,r}(x))\circ\d W_r^k.
\end{equation}
The Stratonovich notation in \eqref{eq:sde} can be defined by the It\^o equation
\begin{equation}\label{eq:ito}
 X_{s,t}(x)=x+\int_s^t\widehat b(r,X_{s,r}(x))\dd r +\int_s^t\sigma(r,X_{s,r}(x))\dd W_r,
\end{equation}
where $\widehat b:=b+\frac12\sum_{k=1}^d(\nabla\sigma_k)\sigma_k$. In particular, when $\sigma$ is constant, $\widehat b=b$.

To state our first theorem, we need to introduce some notation and definitions. The form bound condition used to control $b$ and $\nabla\sigma$ is defined as follows.
\begin{definition}\label{def:form}
Let $E$ be a finite-dimensional Euclidean space. Given $\delta\ge 0$ and $g\in L^1([0,T]; \mR_+)$, we say that a Borel field $v\in L^2_{\mathrm{loc}}([0,T]\times\mR^d;E)$ is form-bounded with $\delta$ and $g$, denoted by $\bF_\delta(g)$, if for almost every $t$, 
\begin{equation*}\label{eq:form}
 \int|v(t,x)|^2\phi(x)^2\dd x
 \leq\delta\int|\nabla\phi|^2
          +g(t)\int\phi^2,
 \quad\phi\in H^1(\mR^d). \tag{FB}
\end{equation*}
\end{definition}

We emphasize that when $d\geq3$, $d<p<\infty$ and $r:=2p/(p-d)$, every $b\in L^{r}(0,T;L^p)$ belongs to $\bF_\delta(g_\delta)$ for each $\delta>0$, with $ g_\delta(t)=C\delta^{-\frac{d}{p-d}}\|b(t)\|_p^{r} \in L^1(0,T)$; and every $b$ satisfies
\begin{equation}\label{eq:tails}
 \lim_{R\to\infty}\mathop{\rm ess\,sup}_{0<t<T} \|b(t)\1_{\{|b(t)|>R\}}\|_d=0, 
\end{equation}
also belongs to $\bF_\delta(g_\delta)$ for each $\delta>0$, with a constant $g_\delta$. These inclusions are fully discussed in \cite[Section~1]{KM24strong} and reference cited therein. 

We make the following basic assumptions on the coefficients $\sigma$ and $b$.

\begin{assumption}\label{aspt:main}
There is $\lambda\in(0,1)$ such that
\begin{equation}\label{eq:ellip}
 \lambda|\xi|^2\leq\frac12\sum_{k=1}^d|\sigma_k(t,x)\cdot\xi|^2 \leq\lambda^{-1}|\xi|^2, \quad \xi\in\mR^d.
\end{equation}
Moreover, 
\[
 b,\nabla\sigma\in\bF_\delta(g)\mbox{ with }g\in L^1(0,T).
\]
Thus, 
\[
\widehat b=b+\frac{1}{2}\sum_{k=1}^d(\nabla\sigma_k)\sigma_k \in\bF_{C\delta}(Cg) \mbox{ with }C=C(d,\lambda)>0. 
\]
\end{assumption}
In the rest of this paper, we put
\begin{equation}\label{def-notation-Theta}
 q:=d+1,\quad P:=50q,
 \quad \Theta(t):=t+\int_0^t g(r)\dd r.
\end{equation}
For a field $f$ and $R>0$, define
\begin{equation*}
 \mathfrak{n}_f(u,v;R)
 :=\Big( \sup_{z\in\mR^d}\int_u^v\int_{B_{2R}(z)}|f(r,y)|^2\dd y\dd r \Big)^{\frac{1}{2}},
 \qquad \mathfrak{n}_f(u,v):=\mathfrak{n}_f(u,v;\sqrt{v-u}).
\end{equation*}

\begin{theorem}\label{thm:flow}
There exists $\overline\delta=\overline\delta(d,\lambda)>0$ such that, whenever Assumption~\ref{aspt:main} holds with $0\leq\delta<\overline\delta$, there are random fields $X_{s,t}(x)$ and $Y_{s,t}(x)$ on the above Wiener space, defined for $0\leq s\leq t\leq T$ and $x\in\mR^d$, with the following properties.
\begin{enumerate}[(i)]
\item For every $(s,x)$, $X_{s,\cdot}(x)$ is a strong solution of \eqref{eq:ito}. The equation holds $\mP$-a.s. for all $t\in[s,T]$; the $\mP$-null set may depend on $(s,x)$. On an event of $\mP$-probability one, both fields are jointly continuous in $(s,t,x)$ and, for every $0\leq s\leq u\leq t\leq T$,
\begin{equation}\label{eq:cocycle}
\begin{gathered}
 X_{u,t}\circ X_{s,u}=X_{s,t},\qquad
 Y_{s,u}\circ Y_{u,t}=Y_{s,t},\\
 X_{s,s}=Y_{s,s}=\Id,\qquad
 X_{s,t}\circ Y_{s,t}=Y_{s,t}\circ X_{s,t}=\Id.
\end{gathered}
\end{equation}
In particular, $X_{s,t}$ is a homeomorphism with inverse $Y_{s,t}$.
\item For every $0\leq s\leq t\leq T$ and $R>0$,
\[
 X_{s,t},Y_{s,t}\in L^P\bigl(\Omega,\mP;W^{1,P}(B_R)\bigr).
\]
Both fields are $\mP$-a.s. locally $1/8$-H\"older continuous in $(s,t,x)$. For every $0<\alpha<1/2$, the trajectories $t\mapsto X_{s,t}(x)$ for fixed $(s,x)$ and $s\mapsto Y_{s,t}(x)$ for fixed $(t,x)$ are $\mP$-a.s. $\alpha$-H\"older continuous on $[s,T]$ and $[0,t]$, respectively.
\item If $f\in\bF_\beta(g_f)$, then, for every integer $m\geq1$ and $s\leq u<v\leq T$,
\begin{equation}\label{eq:occ}
 \bE\left(\int_u^v|f(r,X_{s,r}(x))|\dd r\right)^m
 \leq m!\,C^m(v-u)^{\frac{m}{2}(1-\frac{d}{q})}\mathfrak{n}_f(u,v)^{\frac{2m}{q}},
\end{equation}
where $C$ depends only on $d,\lambda,\beta$, $\|g\|_{L^1(0,T)}$ and $\|g_f\|_{L^1(0,T)}$.
\item For every $(s,x)$ with $s<T$, the solution $X_{s,\cdot}(x)$ belongs to $\sK_2(s,x)$. Moreover, joint uniqueness in law and pathwise uniqueness hold within this class. 
\end{enumerate}
\end{theorem}

\begin{remark}
Krylov \cite{K25strong} decomposes the drift and $\nabla\sigma$ into a part with a small local spatial Morrey bound, uniformly in time, and a part in $L^2_tL^\infty_x$. The localized Fefferman--Phong inequality converts these Morrey conditions into form bounds \cite[Section~1]{KM24strong}. Both works obtain strong existence and pathwise uniqueness among solutions satisfying the corresponding Krylov-type estimates. Under additional assumptions, Krylov also obtains pathwise uniqueness without imposing these estimates on the solutions.
\end{remark}

\medskip

Our next result shows that membership in $L^\infty_tL^d_x$ alone does not ensure strong solvability, even when the diffusion is constant and the drift is divergence-free. Thus, the above theorem cannot be extended to the whole endpoint class without further assumptions.

\begin{theorem}\label{thm:counter}
For every $d\geq2$, there exists a drift
\[
 b\in L^\infty_tL^d_x,\quad
 \div b(t,\cdot)=0,
\]
such that the equation
\begin{equation}\label{eq:counter}
 X_t=\int_0^t b(r,X_r)\dd r+W_t,\quad X_0=0,
\end{equation}
has a continuous weak solution with a unique law on $C([0,1];\mR^d)$, but has no solution adapted to the natural filtration of $W$.
\end{theorem}

The drift can be chosen smooth on $(0,1]\times\mR^d$. The weak solution has the same Gaussian density as Brownian motion at every positive time; see \eqref{eq:marg}. It satisfies the subcritical Krylov estimates in \eqref{eq:krylov} and belongs to $\sK_2(0,0)$. Thus, the example has all these weak-solution properties even though no solution is adapted to the driving Brownian filtration.

\subsection{Main ideas and proof strategy}

In Part~\ref{part:well}, we adapt two ideas to multiplicative noise: parabolic equations for derivative moments, as in \cite{BFGM19stochastic}, and compactness on Wiener space, as in \cite{RZ25sdes,KM24strong}. Using the stochastic transport equation for the inverse flow, we derive a closed parabolic system for tensor moments of its first spatial derivatives. In the weighted energy estimate, integration by parts removes derivatives of $b$ and second derivatives of $\sigma$. The remaining terms are controlled by the form bounds on $b$ and $\nabla\sigma$, whose smallness allows the gradient terms to be absorbed by ellipticity. This leads to uniform moment bounds for the inverse Jacobian and its short-time deviation from the identity; by time reversal, we obtain the corresponding forward bounds.

The representation \eqref{eq:mall} below and a weighted change-of-variable estimate control the Malliavin derivatives using the Jacobian and displacement moment bounds. For coefficients merely measurable in time, Proposition~\ref{prop:mall} also provides uniform continuity under time translations in a weighted $L^2$ norm, without a prescribed power rate. Its proof separates short-time Jacobian increments from the temporal and spatial variation of $\sigma^n$. Together with the spatial derivative and tail bounds, the extended compactness criterion in Lemma~\ref{lem:compact} then implies strong relative compactness of $X^n_{s,t}-\Id$ and $Y^n_{s,t}-\Id$ in a weighted $L^2(\Omega\times\mR^d)$ space at each fixed $(s,t)$.

By the joint H\"older estimates, we extract a common subsequence converging in probability, locally uniformly in $(s,t,x)$, on the original Wiener space. The flow and inverse identities pass to the limit, leading to a stochastic flow of homeomorphisms. Krylov-type estimates identify the limiting trajectories as strong solutions for every fixed starting point. These solutions belong to $\sK_2(s,x)$, where joint uniqueness in law and strong existence yield pathwise uniqueness within this class by the argument of Cherny~\cite[Theorem~3.2]{C02uniqueness}. This implies convergence of the full sequence to a flow independent of the admissible approximation.

\smallskip

In Part~\ref{part:ill}, we show that strong solvability can fail in the endpoint class $L^\infty_tL^d_x$ without a smallness assumption, even for divergence-free drifts and additive noise. We construct the counterexample using a self-similar drift $b(t,x)=t^{-\frac{1}{2}}B(x/\sqrt t)$, whose spatial $L^d$-norm remains constant as $t\downarrow0$. Our use of self-similar variables is inspired by the work of Gallagher--Gallay \cite{GG05uniqueness} on uniqueness for the two-dimensional Navier--Stokes equation with a measure as initial vorticity, together with the related work \cite{GW05global, GGL05uniqueness}. In these works, logarithmic time and parabolic spatial rescaling turn Oseen vortices into stationary Gaussian vorticity profiles and allow one to study vorticity concentrated at the initial time. In particular, the treatment of atomic initial vorticity connects $t\downarrow0$ with dynamics as $\tau:=\log t\to-\infty$. This viewpoint motivates us to pass from an SDE with drift singular at $t=0$ to a stationary diffusion on $(-\infty,0]$.

For a smooth compactly supported profile $B$, the change of variables $Z_s:=e^{-\frac{s}{2}}X_{e^s}$ turns the equation with drift $t^{-\frac{1}{2}}B(x/\sqrt t)$ into
\[
 \d Z_s=\bigl(-Z_s/2+B(Z_s)\bigr)\dd s+\dd\beta_s,
 \quad -\infty<s\leq0.
\]
The initial time $t=0$ becomes the infinite past of a smooth stationary diffusion with invariant law $\mu_B$. For any weak solution with probability measure $\mP$, the conditional law of $Z_0$ under $\mP$, given the Brownian increment paths over the last $N$ unit time intervals, is the pushforward of $\mu_B$ by the composition of the maps driven by these paths. Proposition~\ref{prop:strong} shows that the original SDE is strongly solvable exactly when the limit of these conditional laws is a Dirac mass $\mP$-a.s.. The limiting conditional distribution is the sample measure $\mu_{\underline\omega^-}$ of the stationary diffusion, where $\underline\omega^-$ denotes the past Brownian noise realization.

For the profile $B_K$ constructed in Section~\ref{sec:profile}, a moving orthonormal frame expresses the growth of tangent vectors in $\mR^2\times\{0\}$ in terms of a uniformly elliptic auxiliary diffusion. Its invariant density relative to Gaussian measure satisfies $\lambda_K=-\frac12+\frac14I_K$, where $\lambda_K$ is the largest Lyapunov exponent on this subspace and $I_K$ is the relative Fisher information of the auxiliary invariant density. We derive this identity directly from the stationary equation of the auxiliary diffusion, using an integration-by-parts argument similar to that in \cite[Section~3.1]{BBP22regularity}. A local compactness argument in the stationary equation shows that $I_K\to\infty$ as $K\to\infty$, so $\lambda_K>0$ for all sufficiently large $K$. The random SRB theorem of Ledrappier--Young \cite{LY88entropy}, with the direct geometric proof of Blumenthal--Young \cite{BY19equivalence}, then ensures that the sample measures have absolutely continuous conditional measures along positive-dimensional unstable manifolds. Such sample measures cannot be Dirac, so Proposition~\ref{prop:strong} rules out a strong solution. To the best of our knowledge, this seems to be the first use of positive Lyapunov exponents and random SRB theory to prove nonexistence of strong solutions for SDEs.

The construction of sample measures used here mainly follows \cite{BY19equivalence}. For further background, we refer to Le Jan \cite{L87equilibre} on statistical equilibrium and its relation to Lyapunov exponents, and to Flandoli--Gess--Scheutzow \cite{FGS17synchronization} on synchronization by noise.

\smallskip

\paragraph{\em Organization of the paper.} Part~\ref{part:well} contains the proof of Theorem~\ref{thm:flow}. Section~\ref{sec:pde} introduces the approximations and establishes the preliminary estimates. Section~\ref{sec:deriv} derives the spatial and Malliavin derivative estimates, which are used in Section~\ref{sec:compact} to construct the limiting stochastic flow and prove uniqueness. Part~\ref{part:ill} contains the proof of Theorem~\ref{thm:counter}. Section~\ref{sec:entry} characterizes strong solvability by whether the sample measures are Dirac masses. Section~\ref{sec:lyap} constructs a profile with a positive top Lyapunov exponent, and Section~\ref{sec:srb} uses random SRB theory to rule out strong solutions. Appendix~\ref{app:pde} collects some proofs of PDE estimates used in Part~\ref{part:well}.

\part{Well-posedness}\label{part:well}

We now turn to the proof of Theorem~\ref{thm:flow}. Throughout Part~\ref{part:well}, we always assume that the coefficients $b$ and $\sigma$ satisfy Assumption \ref{aspt:main}. And the smallness thresholds denoted by $\overline\delta$ may differ between statements, with their dependence given in each case, and are always chosen at most one. The threshold in Theorem~\ref{thm:flow} is chosen to satisfy all these requirements, including those for the approximations with ellipticity parameter $\lambda/2$.

\section{Approximation}\label{sec:pde}

The proof of Theorem~\ref{thm:flow} uses three preliminary results. The Krylov-type estimate controls the drift along the approximating solutions $X^n$ and permits passage to the limit in the SDE. We use the backward equation to prove uniqueness in law for solutions in $\sK_2(s,x)$. Finally, an $L^2$-to-$L^\infty$ estimate places the constructed solutions in this class. We first construct approximations and prove the Krylov-type estimates and moment bounds for time increments. The proofs of the PDE estimates and of the results concerning transition semigroups and the class $\sK_2(s,x)$ defined by \eqref{eq:k2} are collected in Appendix~\ref{app:pde}.

All approximating flows and their limits are constructed on the fixed classical Wiener space $(\Omega,\sF,\mP)$ from the introduction. Throughout Sections~\ref{sec:pde}--\ref{sec:compact}, $\bE$ denotes expectation under $\mP$, and convergence in probability refers to $\mP$, unless a different measure is explicitly given. We write $\mP\times\rho$ for the product measure $\mP(\dd\omega)\rho(x)\dd x$.

\subsection{Weights and admissible approximations}

For $\kappa>d$ and $z\in\mR^d$, set
\begin{equation*}
 \rho_z(x):=(1+|x-z|^2)^{-\frac{\kappa}{2}},\quad \rho:=\rho_0.
\end{equation*}
For $1\leq p<\infty$, write $L^p_{\rho_z}:=L^p(\mR^d,\rho_z(x)\dd x)$, with $\|u\|_{p,\rho_z}^p:=\int|u|^p\rho_z$. The weighted Sobolev space $W^{1,p}_{\rho_z}$ consists of functions with $u,\nabla u\in L^p_{\rho_z}$, with norm
\[
 \|u\|_{W^{1,p}_{\rho_z}}^p
       :=\|u\|_{p,\rho_z}^p+\|\nabla u\|_{p,\rho_z}^p.
\]
We write $\langle U,V\rangle_{\rho_z}:=\int U\cdot V\,\rho_z$, where the dot denotes the componentwise Euclidean inner product.

We first record a weighted form inequality and local $L^2$ bounds for form-bounded fields.

\begin{lemma}\label{lem:wform}
Let $\alpha\geq0$, let $h\in L^1(0,T)$ be nonnegative, and suppose that $v\in\bF_\alpha(h)$ in the sense of Definition~\ref{def:form}. Fix $\kappa>d$ and use the weights $\rho_z$ defined above. Then the following estimates hold.
\begin{enumerate}[(i)]
\item For almost every $t\in(0,T)$, for all $z\in\mR^d$, $\eps>0$, and $U\in W^{1,2}_{\rho_z}$,
\begin{equation}\label{eq:wform}
\int |v(t)|^2|U|^2\rho_z
\leq (1+\eps)\alpha\int|\nabla U|^2\rho_z +\bigg(h(t)+\frac{\kappa^2}{4}(1+\eps^{-1}) \alpha \bigg)\int|U|^2\rho_z.
\end{equation}
\item There is a constant $C=C(d,\kappa)$ such that
\begin{equation}\label{eq:loc-src}
\begin{aligned}
\sup_{z\in\mR^d}\int |v(t)|^2\rho_z
& \leq C\bigl(\alpha+h(t)\bigr),\\
\sup_{z\in\mR^d}\int_s^t\int_{B_2(z)}|v(r)|^2\dd r
& \leq C\left(\alpha(t-s)+\int_s^t h(r)\dd r\right).
\end{aligned}
\end{equation}
\end{enumerate}
\end{lemma}
\begin{proof}
For $U\in W^{1,2}_{\rho_z}$, the function $\phi:=|U|\sqrt{\rho_z}$ belongs to $H^1$, and
\[
 \nabla\phi=\sqrt{\rho_z}\,\nabla|U|
       +\frac{|U|}{2\sqrt{\rho_z}}\nabla\rho_z.
\]
Using $|\nabla|U||\leq|\nabla U|$, $|\nabla\rho_z|\leq\kappa\rho_z$, and Young's inequality, we obtain
\begin{align*}
 & \int|v|^2|U|^2\rho_z
  \leq \alpha\int|\nabla\phi|^2+h\int\phi^2\\
 \leq{} & (1+\eps)\alpha\int|\nabla U|^2\rho_z
       +\left(h+\frac{\kappa^2}{4}(1+\eps^{-1})\alpha\right)
          \int|U|^2\rho_z.
\end{align*}
This proves \eqref{eq:wform}.

For the local estimate, choose $\chi\in C_c^\infty(B_3)$ with $0\leq\chi\leq1$ and $\chi=1$ on $B_2$, and put $\chi_z(y):=\chi(y-z)$. Applying the form bound to $\chi_z$, we get
\[
 \int_{B_2(z)}|v(r,y)|^2\dd y
 \leq\alpha\|\nabla\chi\|_2^2 +h(r)\|\chi\|_2^2
 \leq C\bigl(\alpha+h(r)\bigr).
\]
Integrating over the time interval $[s,t]$ and taking the supremum over $z$ proves the second estimate in \eqref{eq:loc-src}.

Finally, since $\kappa>d$, we have $\sqrt{\rho_z}\in H^1$. Applying the form bound directly to $\sqrt{\rho_z}$ yields
\[
 \int|v(t)|^2\rho_z
 \leq\alpha\int|\nabla\sqrt{\rho_z}|^2+h(t)\int\rho_z
 \leq\left(\frac{\kappa^2}{4}\alpha+h(t)\right)\int\rho
 \leq C\bigl(\alpha+h(t)\bigr).  \qedhere
\]
\end{proof}

We next describe the smooth approximations used throughout the proof. All subsequent estimates use only the properties in the following definition.

\begin{definition}\label{def:approx}
We say a sequence $(b^n,\sigma^n)$ is an admissible approximation, if the coefficients and their spatial derivatives are bounded on $[0,T]\times\mR^d$, and they satisfy the following conditions:
\begin{enumerate}[(a)]
\item $b^n,\nabla\sigma^n\in\bF_\delta(g)$;
\item \eqref{eq:ellip} holds with $\lambda$ replaced by $\lambda/2$;
\item $b^n\to b$, $\sigma^n\to\sigma$, and $\nabla\sigma^n\to\nabla\sigma$ locally in space-time $L^2$.
\end{enumerate}

\end{definition}

\begin{proposition}\label{prop:approx}
There exists $\overline\delta=\overline\delta(d,\lambda)>0$ such that, for every $0\leq\delta<\overline\delta$, an admissible approximation $(b^n,\sigma^n)$ exists, and
\begin{equation}\label{eq:ito-lim}
 \widehat b^n:=b^n+\frac12\sum_{k=1}^d(\nabla\sigma_k^n)\sigma_k^n
       \to\widehat b
       \quad\mbox{in }L^2_{\mathrm{loc}}([0,T]\times\mR^d).
\end{equation}
\end{proposition}
\begin{proof}

Fix a nonnegative $\eta\in C_c^\infty(B_1)$ with $\int\eta=1$, chosen once for each dimension, and put $\eta_n(x):=n^d\eta(nx)$. We combine space convolution with a cutoff on times when $g$ is large. Since convolution alone does not necessarily produce a bounded drift, and it may also destroy the lower ellipticity bound, the cutoff is essential. For a constant $a_0>0$ to be chosen below, set
\[
\begin{aligned}
 E_n:=\{t:g(t)\leq a_0n^2\}, \quad b^n:=\1_{E_n}(b*_x\eta_n),\quad \sigma^n:=\1_{E_n}(\sigma*_x\eta_n)+\1_{E_n^c}\sqrt{2\lambda}\,I.
\end{aligned}
\]
We first verify (a). Let $v$ denote either $b$ or $\nabla\sigma$, and put $v^n:=\1_{E_n}(v*_x\eta_n)$. Jensen's inequality implies 
\[
 |v^n(t,x)|^2
 \leq\1_{E_n}(t)\int\eta_n(z)|v(t,x-z)|^2\dd z.
\]
For $\phi\in C_c^\infty(\mR^d)$, applying \eqref{eq:form} to the translated test function $y\mapsto\phi(y+z)$, we obtain that 
\begin{align*}
& \int|v^n(t,x)|^2\phi(x)^2\dd x
 \leq \1_{E_n}(t)\int\eta_n(z)\int|v(t,y)|^2\phi(y+z)^2\dd y\dd z\\
\leq{} & \1_{E_n}(t)\bigl(\delta\|\nabla\phi\|_2^2+g(t)\|\phi\|_2^2\bigr)\\
\leq{} & \delta\|\nabla\phi\|_2^2+g(t)\|\phi\|_2^2.
\end{align*}
By density, the estimate holds for every $\phi\in H^1(\mR^d)$. Since $E_n$ depends only on time, and $\sigma^n$ is spatially constant outside $E_n$, we have
\[
 \nabla\sigma^n=\1_{E_n}\bigl((\nabla\sigma)*_x\eta_n\bigr),
 \qquad b^n,\nabla\sigma^n\in\bF_\delta(g).
\]
This proves (a).
For every spatial multi-index $\alpha=(\alpha_1,\ldots,\alpha_d)\in\mN_0^d$, write $\partial^\alpha:=\partial_1^{\alpha_1}\cdots\partial_d^{\alpha_d}$. For $t\in E_n$, we obtain from the form bound
\[
 \sup_x|\partial^\alpha b^n(t,x)|
 \leq C_\alpha n^{|\alpha|}\bigl(\delta n^2+g(t)\bigr)^{\frac{1}{2}}
 \leq C_\alpha n^{|\alpha|+1}(\delta+a_0)^{\frac{1}{2}}.
\]
As $\sigma$ is bounded, all spatial derivatives of $\sigma^n$ are also bounded. Thus the required smoothness and boundedness hold.

To check ellipticity (b), first apply \eqref{eq:form} to a cutoff equal to one on $B_r(x)$ and supported in $B_{2r}(x)$. For $v=b,\nabla\sigma$, and $0<r\leq1$, we have
\begin{equation*}
 \int_{B_r(x)}|v(t,y)|^2\dd y
       \leq C\bigl(\delta r^{d-2}+g(t)r^d\bigr).
\end{equation*}
Let
\[
 \overline\sigma_{n,x}(t):=\frac{1}{|B_{1/n}(x)|}
       \int_{B_{1/n}(x)}\sigma(t,y)\dd y.
\]
Since $\sigma^n(t,x)$ is the mean with respect to $\eta_n(x-y)\dd y$, replacing it by $\overline\sigma_{n,x}(t)$ can only increase the weighted variance. Using Poincar\'e's inequality, we obtain
\begin{align*}
 &\int|\sigma(t,y)-\sigma^n(t,x)|^2\eta_n(x-y)\dd y\\
 \leq{} & \int|\sigma(t,y)-\overline\sigma_{n,x}(t)|^2\eta_n(x-y)\dd y \leq Cn^d\int_{B_{1/n}(x)}|\sigma(t,y)-\overline\sigma_{n,x}(t)|^2\dd y\\
 \leq{} & Cn^{d-2}\int_{B_{1/n}(x)}|\nabla\sigma(t,y)|^2\dd y 
 \leq C\bigl(\delta+n^{-2}g(t)\bigr)
 \leq C(\delta+a_0).
\end{align*}
By the variance identity, for every $t\in E_n, x\in\mR^d$ and $\xi\in\mR^d$, we have
\begin{align*}
 & |(\sigma^n(t,x))^\top\xi|^2
  = \int|\sigma(t,y)^\top\xi|^2\eta_n(x-y)\dd y-\int|(\sigma(t,y)-\sigma^n(t,x))^\top\xi|^2 \eta_n(x-y)\dd y\\
 \geq{} & \bigl[2\lambda-C(\delta+a_0)\bigr]|\xi|^2.
\end{align*}
Here $C$ depends only on $d$ and the fixed mollifier $\eta$. Set $a_0:=\lambda/(2C)$ and $\overline\delta:=\min\{1,a_0\}$. Then $C(\delta+a_0)\leq\lambda$ whenever $\delta<\overline\delta$, which proves the lower bound in (b). The upper bound follows from Jensen's inequality; outside $E_n$, both bounds follow from the choice of the constant matrix.

For (c), we only need to note that $\1_{E_n}\to1$ almost everywhere in time, so for any locally square-integrable field $v$ we have
\[
 \1_{E_n}(v*_x\eta_n)-v=\1_{E_n}(v*_x\eta_n-v)-\1_{E_n^c}v \to0\quad\mbox{in local }L^2.
\]

Finally, for any admissible approximation, write
\[
 (\nabla\sigma_k^n)\sigma_k^n-(\nabla\sigma_k)\sigma_k
   =(\nabla\sigma_k^n-\nabla\sigma_k)\sigma_k^n
       +(\nabla\sigma_k)(\sigma_k^n-\sigma_k).
\]
The first term tends to zero locally in $L^2$ by (b) and (c). The second term does so because $\nabla\sigma_k$ is locally square integrable, while $\sigma_k^n-\sigma_k$ is uniformly bounded and tends to zero locally in measure. Together with $b^n\to b$, we get \eqref{eq:ito-lim}.
\end{proof}

Put
\begin{equation*}
a:=\frac12\sigma\sigma^\top,\quad
c:=b-\frac12\sum_k(\div\sigma_k)\sigma_k, \quad L:=\nabla\cdot(a\nabla)+c\cdot\nabla.
\end{equation*}
Equivalently, $L\varphi=a:\nabla^2\varphi+\widehat b\cdot\nabla\varphi$, where $a:\nabla^2\varphi:=\sum_{i,j=1}^d a_{ij}\partial_{ij}\varphi$. The same notation with superscript $n$ is used for the approximations. The fields $c,\nabla a$ belong to $\bF_{C\delta}(Cg)$. Their approximations converge locally in $L^2$ and belong to $\bF_{C\delta}(Cg)$, where $C$ depends only on $d,\lambda$. These statements follow directly from boundedness of $\sigma$ and the product rule.

\subsection{Uniform estimates}

We first state the parabolic estimates used below, with coefficients $a,c$ as above. The estimates hold uniformly for admissible approximations and after reversing time. The proof is given in Appendix~\ref{app:energy}.

\begin{lemma}\label{lem:pde}
Fix $\kappa>d$ and let $q=d+1$. There exists $\overline\delta=\overline\delta(d,\lambda,\kappa)>0$ such that, for every $0\leq\delta<\overline\delta$, the following estimates hold. Let $F\in\bF_\beta(g_F)$ and let $u$ solve
\[
 \partial_tu=\nabla\cdot(a\nabla u)+c\cdot\nabla u+F,
 \quad u(0)=\varphi,
\]
with coefficients smooth in space. Then, for every $S\leq T$,
\begin{align}
 &\sup_{t\leq S}\int(|u(t)|^q+|\nabla u(t)|^q)\rho_z
 +\int_0^S\int\rho_z\left|\nabla|\nabla u|^{\frac{q}{2}}\right|^2\dd t
 \nonumber\\
 \leq{} & C\left[\int(|\varphi|^q+|\nabla\varphi|^q)\rho_z
             +\int_0^S\int|F|^2\rho_z\dd t\right].
 \label{eq:energy}
\end{align}
If $\varphi=0$, then
\begin{equation}\label{eq:pot}
 \|u\|_{L^\infty([0,S]\times\mR^d)}
 \leq C \mathfrak{n}_F(0,S;1)^{\frac{2}{q}}.
\end{equation}
Here $C$ depends only on $d,\lambda,\kappa,T,\beta$, $\|g\|_{L^1(0,T)}$ and $\|g_F\|_{L^1(0,T)}$. No smallness of $\beta$ is required.
\end{lemma}

Let $X^n$ denote the smooth flow for \eqref{eq:ito} with the $n$-th coefficients $(\sigma^n, \widehat b^n)$, driven by the canonical Brownian motion $W$ under $\mP$, and let $Y^n:=(X^n)^{-1}$. The existence of these smooth spatial flows for bounded coefficients with bounded spatial derivatives follows from Kunita's classical theory on smooth SDE, see \cite{K90stochastic}. The smooth flow and inverse identities are understood on a common event of $\mP$-measure one.

We now derive the Krylov-type estimates and moment bounds from the $L^\infty$ estimate for the parabolic equation in Lemma~\ref{lem:pde} through It\^o's formula and the Markov property.

\begin{lemma}\label{lem:occ}
There exists $\overline\delta=\overline\delta(d,\lambda)>0$ such that, for every $0\leq\delta<\overline\delta$, the following estimates hold uniformly in $n$.
\begin{enumerate}[(i)]
\item If $f\in\bF_\beta(g_f)$, then, for every integer $m\geq1$ and $s\leq u<v\leq T$,
\[
\sup_{x}\bE\left(\int_u^v|f(r,X^n_{s,r}(x))|\dd r\right)^m
\leq m!\,C^m(v-u)^{\frac{m}{2}(1-\frac{d}{q})}\mathfrak{n}_f(u,v)^{\frac{2m}{q}}.
\]
Here $C$ depends only on $d,\lambda,\beta$, $\|g\|_{L^1(0,T)}$ and $\|g_f\|_{L^1(0,T)}$.
\item For every real $p\geq1$,
\begin{equation}\label{eq:disp}
 \sup_x\bE\sup_{u\leq r\leq v}
       |X^n_{s,r}(x)-X^n_{s,u}(x)|^p
       \leq C_p(v-u)^{\frac{p}{2}},
 \quad s\leq u\leq v\leq T.
\end{equation}
Here $C_p$ depends only on $d,\lambda,p$ and $\|g\|_{L^1(0,T)}$. The same time-increment estimate holds for $r\mapsto Y^n_{t-r,t}(x)$, $0\leq r\leq t$, with $t$ fixed and with the same constant.
\end{enumerate}
Therefore, for every $p\geq1$,
\begin{equation}\label{eq:disp-u}
 \sup_{n,s,t,x}\bE\bigl(
 |X^n_{s,t}(x)-x|^p+|Y^n_{s,t}(x)-x|^p\bigr)<\infty.
\end{equation}
\end{lemma}
\begin{proof}
\emph{(i) Krylov-type estimates.}
By a standard approximation argument, it suffices to consider bounded $f$ that are smooth in space. It is enough to prove the estimates for $s=u$, the case $s\leq u$ then follows from the Markov property by conditioning on $\sF_u$.

Put $\ell:=v-u$ and, for $0\leq\tau\leq1$, set
\[
 B_\tau:={(W_{u+\ell\tau}-W_u)}/{\sqrt\ell},\quad
 Z_\tau:={(X^n_{u,u+\ell\tau}(z)-z)}/{\sqrt\ell}.
\]
The process $B$ is a Brownian motion with respect to $(\sF_{u+\ell\tau})_{0\leq\tau\leq1}$. Define
\[
b_\ell(\tau,y):=\sqrt\ell\,b^n(u+\ell\tau,z+\sqrt\ell\,y), 
\quad \sigma_\ell(\tau,y):=\sigma^n(u+\ell\tau,z+\sqrt\ell\,y)
\] 
and $f_\ell(\tau,y):=\sqrt\ell\,f(u+\ell\tau,z+\sqrt\ell\,y)$. Then $Z$ solves the It\^o equation with diffusion $\sigma_\ell$, initial point zero, and drift $\widehat b_\ell(\tau,y):=\sqrt\ell\,\widehat b^n(u+\ell\tau,z+\sqrt\ell\,y)$.

Obviously, ellipticity is preserved. Applying \eqref{eq:form} to $x\mapsto\phi((x-z)/\sqrt\ell)$ and changing variables, we obtain
\[
 b_\ell,\nabla_y\sigma_\ell\in\bF_\delta(k_\ell),
 \quad f_\ell\in\bF_\beta(h_\ell),
\]
where $k_\ell(\tau):=\ell g(u+\ell\tau)$ and $h_\ell(\tau):=\ell g_f(u+\ell\tau)$. Their time integrals satisfy
\[
 \int_0^1 k_\ell(\tau)\dd\tau\leq\|g\|_{L^1(0,T)},\qquad
 \int_0^1 h_\ell(\tau)\dd\tau\leq\|g_f\|_{L^1(0,T)}.
\]
Under the same change of variables, we also have
\[
 \mathfrak{n}_{f_\ell}(0,1;1)^2
 =\ell^{-\frac{d}{2}}\mathfrak{n}_f(u,v)^2.
\]
Let $V$ solve the backward equation for $Z$ with source $-|f_\ell|$ and terminal value $V(1)=0$. It\^o's formula and Lemma~\ref{lem:pde}, applied on $[0,1]$ with $\kappa=2q$ yield 
\begin{equation}\label{eq:occ-aux}
\begin{aligned}
 & \bE\int_u^v|f(r,X^n_{u,r}(z))|\dd r
  \leq \sqrt\ell\,\|V\|_\infty
 \leq C\sqrt\ell\,\mathfrak{n}_{f_\ell}(0,1;1)^{\frac{2}{q}}\\
 ={} & C\ell^{\frac{1}{2}(1-\frac{d}{q})}\mathfrak{n}_f(u,v)^{\frac{2}{q}}.
\end{aligned}
\end{equation}

To obtain higher moments, write $X_r:=X^n_{s,r}(x)$ and $I_r:=\int_r^v|f(a,X_a)|\dd a$. Applying \eqref{eq:occ-aux} on $[r,v]$ and using the Markov property, we obtain
\[
 \bE[I_r\mid\sF_r]\leq A
 :=C\ell^{\frac{1}{2}(1-\frac{d}{q})}\mathfrak{n}_f(u,v)^{\frac{2}{q}},
 \quad u\leq r\leq v,
\]
since $q>d$, $v-r\leq\ell$, and $\mathfrak{n}_f(r,v)\leq\mathfrak{n}_f(u,v)$ for $r<v$; for $r=v$ the left-hand side vanishes.
Using the identity $I_r^m=m\int_r^v|f(a,X_a)|I_a^{m-1}\dd a$, conditioning at $a$, and induction on $m$, we obtain
\[
 \bE[I_r^m\mid\sF_r]
 \leq m(m-1)!A^{m-1}\bE[I_r\mid\sF_r]
 \leq m!A^m.
\]
Taking $r=u$ and expectations proves (i) for bounded spatially smooth $f$.

We record two consequences of (i). A covering by balls of radius two implies
\[
 \mathfrak{n}_f(u,v)^2
 \leq C_d\max\{1,\ell^{\frac{d}{2}}\}\mathfrak{n}_f(u,v;1)^2.
\]
Consequently, the right-hand side of (i) is bounded by $m!\,C_T^m\mathfrak{n}_f(u,v;1)^{\frac{2m}{q}}$, where $C_T$ depends on the parameters in (i) and on $T$. Next, applying \eqref{eq:form} to cutoffs that equal one on $B_{2\sqrt\ell}(z)$ and vanish outside $B_{3\sqrt\ell}(z)$ yields
\[
 \mathfrak{n}_f(u,v)^2
 \leq C_d\ell^{\frac{d}{2}}\left(\beta+\int_u^v g_f(r)\dd r\right).
\]
Substitution into (i) bounds its left-hand side by $m!\,C^m\ell^{\frac{m}{2}}$, with the same constant dependence as in (i).

\emph{(ii) Time increments.} Since $\widehat b^n\in\bF_{C\delta}(Cg)$ with $C$ depending only on $d,\lambda$, the preceding time bound applies to the It\^o drift. For real $p\geq1$, choose $m=\lceil p\rceil$. Hence, we have
\[
 \bE\left(\int_u^v|\widehat b^n(r,X^n_{s,r}(x))|\dd r\right)^p
 \leq (m!C^m)^{\frac{p}{m}}(v-u)^{\frac{p}{2}}
 \leq C_p(v-u)^{\frac{p}{2}},
\]
where $C_p$ depends only on $d,\lambda,p$ and $\|g\|_{L^1(0,T)}$. The boundedness of $\sigma^n$, and the Burkholder--Davis--Gundy inequality now yield
\begin{align*}
& \bE\sup_{u\leq r\leq v}|X^n_{s,r}(x)-X^n_{s,u}(x)|^p
 \leq C_p\bE\left(\int_u^v|\widehat b^n(r,X^n_{s,r}(x))|\dd r\right)^p +C_p(v-u)^{\frac{p}{2}}\\
\leq{} & C_p(v-u)^{\frac{p}{2}}.
\end{align*}
Taking the supremum over $x$ proves \eqref{eq:disp}.

For the inverse flow, fix $t$, reverse time by $r\mapsto t-r$, and set
\[
 \overline W_r:=W_t-W_{t-r},\quad
 \overline b(r,x):=-b^n(t-r,x),\quad
 \overline\sigma(r,x):=-\sigma^n(t-r,x).
\]
If $\overline X$ is the flow with these coefficients driven by $\overline W$, then $Y^n_{s,t}=\overline X_{0,t-s}$. Its It\^o drift is $-b^n+\frac12\sum_k(\nabla\sigma_k^n)\sigma_k^n$ evaluated at time $t-r$, so the same form bounds apply. The forward estimate just proved therefore applies to the inverse flow and proves \eqref{eq:disp-u}.
\end{proof}

We next collect some results about the backward equation. The proof is given in Appendix~\ref{app:semi}.

\begin{proposition}\label{prop:semi}
There exists $\overline\delta=\overline\delta(d,\lambda)>0$ such that, for every $0\leq\delta<\overline\delta$, the following assertions hold.
For $\varphi\in C_c^\infty(\mR^d)$ and fixed $t\leq T$, put
\[
 u_n(s,x):=P^n_{s,t}\varphi(x):=\bE\varphi(X^n_{s,t}(x)).
\]
Here $(P^n_{s,t})_{0\leq s\leq t\leq T}$ is the time-inhomogeneous transition semigroup associated with $X^n$.
\begin{enumerate}[(i)]
\item The functions $u_n$ converge locally uniformly on $[0,t]\times\mR^d$ to a function $P_{s,t}\varphi(x)$ independent of the admissible approximation. For each $(s,x)$, all path-space subsequential limits of $X^n_{s,\cdot}(x)$ have the same finite-dimensional distributions, given by these limiting operators.
\item For $0\leq r<t\leq T$,
\begin{equation*}
 \|P^n_{r,t}\varphi\|_\infty+
 \|P_{r,t}\varphi\|_\infty
 \leq C(t-r)^{-\frac{d}{4}}\|\varphi\|_2.
\end{equation*}
Every weak limit $\mathsf Q$ of the path laws $\Law(X^n_{s,\cdot}(x))$ has time-$t$ marginal densities $p_{s,t}(x,\cdot)$ for $t>s$, with
\begin{equation}\label{eq:dens}
 \|p_{s,t}(x,\cdot)\|_2\leq C(t-s)^{-\frac{d}{4}}.
\end{equation}
Therefore, every weak solution with path law $\mathsf Q$ belongs to $\sK_2(s,x)$.
\item The function $u(r):=P_{r,t}\varphi$ is bounded and continuous, and
\[
\partial_ru+a:\nabla^2u+\widehat b\cdot\nabla u=0 \quad\mbox{in }L^2((0,t)\times\mR^d),\quad
u(t)=\varphi.
\]
Moreover,
\begin{equation}\label{eq:l2back}
\sup_{0\leq r\leq t}\|u(r)\|_{H^1}^2
+\int_0^t\bigl(\|\nabla^2u(r)\|_2^2+\|\partial_ru(r)\|_2^2 +\|\widehat b(r)\cdot\nabla u(r)\|_2^2\bigr)\dd r 
\leq C\|\varphi\|_{H^1}^2.
\end{equation}
The estimate also holds uniformly for $u_n(r)=P^n_{r,t}\varphi$.
\end{enumerate}
The constants $C$ in (ii) and (iii) depend only on $d,\lambda,T$ and $\|g\|_{L^1(0,T)}$.
\end{proposition}

\section{Derivative moments of the stochastic flows}\label{sec:deriv}

In this section, we prove the estimates needed for compactness of $X^n$ and $Y^n$ on the original Wiener space. The argument has three parts: bounds on spatial derivatives, bounds and time translations for Malliavin derivatives, and joint H\"older bounds in $(s,t,x)$. Section~\ref{sec:compact} will combine the first two to obtain compactness at fixed $(s,t)$, and use the third to pass to locally uniform convergence.

\medskip

Recall $q=d+1$ and $P=50q$. We define the maximal moment order
\[
M:=10^3q^2
\]
to cover the applications of H\"older's inequality below. Set 
\[
J^n_{s,t}:=\nabla X^n_{s,t}~\mbox{ and }~ K^n_{s,t}:=\nabla Y^n_{s,t}. 
\]
We use
\begin{equation}\label{def-weight-funct}
\rho(x):=(1+|x|^2)^{-q}~ \textit{ and }~ \rho_*(x):=(1+|x|^2)^{-\frac{q}{2}}.
\end{equation}

\subsection{The moment system and its energy estimate}
In this subsection, we derive a closed parabolic system for the inverse flows, which will be later used to obtain the desired moment estimates for the derivatives of the flows. 

Recall that $(b^n,\sigma^n)$ is an admissible approximation in the
sense of Definition~\ref{def:approx}. In particular, $b^n,\nabla\sigma^n\in\bF_\delta(g)$, with the same $\delta>0$ and $g\in L^1(0,T)$
for all $n$. All estimates below are uniform in $n$. For notational simplicity, we suppress the superscript $n$ on the approximating coefficients and flows. 

\subsubsection{Transport equation and derivative moments}

Let $\varphi\in C^\infty(\mR^d)$ be deterministic, with all derivatives of positive order bounded. In particular, $\varphi$ is globally Lipschitz and has at most linear growth. Put
\begin{equation}\label{def-u}
 u(t,y):=\varphi\bigl(Y_{s,t}(y)\bigr).    
\end{equation}
For each fixed $n$, the smooth approximating coefficients ensure that the spatial derivatives of $u$ have the integrability needed to justify the differentiations and expectations below.

For a smooth vector field $v$, write
\[
A_v u:=v\cdot\nabla u
\]
and let $A_v$ act componentwise. 

For coefficients that are also smooth in time, the method of
stochastic characteristics \cite[Section~6.1]{K90stochastic}
shows that $u$ defined in \eqref{def-u} satisfies
\begin{equation}\label{transp-eq-Ito}
 \d u=\Big[-A_bu+\frac12\sum_k A_{\sigma_k}^{2}u\Big]\dd t
       -\sum_k A_{\sigma_k}u\,\d W_t^k,
\end{equation}
with $u(s,\cdot)=\varphi$. For the coefficients used here, which are merely measurable in time, one can directly use the equation for the inverse flow established in \cite[Theorem~2.4]{LM15space}, as the coefficients $(\widehat b^n,\sigma^n)$ satisfy the regularity assumptions therein.

To estimate moments of the Jacobian matrices $\nabla X$ and $\nabla Y$, we need to derive the equation for $\nabla u$. For this purpose, we introduce some notation. Given vector fields $v, z:\mR^d\to\mR^d$, define
\[
 \cL_v z:=A_vz+(\nabla v)^\top z=v\cdot \nabla z+(\nabla v)^\top z.
\]
For a smooth function $f:\mR^d\to\mR$, we have 
\begin{equation}\label{commut-nabla}
 \nabla(A_vf)=A_v\nabla f+(\nabla v)^\top\nabla f=\cL_v(\nabla f), 
\end{equation}
and 
\[
 \cL_v^2z:=\cL_v(\cL_vz)=A_v^2z+2(\nabla v)^\top A_vz
       +\bigl[A_v((\nabla v)^\top)+((\nabla v)^\top)^2\bigr]z.
\]

For an index tuple 
\[
I=(i_1,\ldots,i_m)\in\{1,\ldots,d\}^m,
\]
write
\[
 I^{\alpha\to j}:=(i_1,\ldots,i_{\alpha-1},j,i_{\alpha+1},\ldots,i_m),
 \quad 1\leq\alpha\leq m,\quad 1\leq j\leq d.
\]
Thus, $I^{\alpha\to j}$ replaces the $\alpha$-th entry of $I$ by $j$. 

For a smooth vector field $v$ and a tensor field $U=(U_I)_I$ of order $m$, define
\begin{equation}\label{eq:def_abl}
 (A_vU)_I:=v\cdot\nabla U_I,
 \quad (B_vU)_I:=\sum_{\alpha=1}^m\sum_{j=1}^d \partial_{i_\alpha}v_j\,U_{I^{\alpha\to j}},
 \quad L_v:=A_v+B_v.
\end{equation}
For $m=1$, we have $B_vU=(\nabla v)^\top U$ and $L_v=\cL_v$. The operator norm of $B_v$ is at most $C_dm|\nabla v|$.

We now derive the equation for the product of derivatives. Put $w:=\nabla u$, so that $w(s,y)=\nabla\varphi(y)$. Noting that, by \eqref{commut-nabla},
\[
 \nabla(A_v^2u)=\cL_v\bigl(\nabla(A_vu)\bigr)
                =\cL_v(\cL_vw)=\cL_v^2w, 
\]
the transport equation \eqref{transp-eq-Ito} yields
\begin{equation}\label{eq:w}
\d w=\Big[-\cL_bw+\frac12\sum_k\cL_{\sigma_k}^2w \Big]\dd t -\sum_k\cL_{\sigma_k}w\, \d W_t^k.
\end{equation}
Set $U:=w^{\otimes m}$, so that $U_I=\prod_{\alpha=1}^m w_{i_\alpha}$ and $U(s)=(\nabla\varphi)^{\otimes m}$. For each $I$, It\^o's formula implies
\begin{equation}\label{eq:dUI}
\begin{aligned}
 & \d U_I = \sum_{\alpha=1}^m
       \Big(\prod_{\beta\ne\alpha}w_{i_\beta}\Big)\d w_{i_\alpha}+\sum_{\alpha<\beta}
       \Big(\prod_{\gamma\notin\{\alpha,\beta\}}w_{i_\gamma}\Big)
       \d\langle w_{i_\alpha},w_{i_\beta}\rangle_t\\
       ={} & \sum_{\alpha=1}^m
       \Big(\prod_{\beta\ne\alpha}w_{i_\beta}\Big)\d w_{i_\alpha}+\sum_{\alpha<\beta}
       \Big(\prod_{\gamma\notin\{\alpha,\beta\}}w_{i_\gamma}\Big) \sum_k(\cL_{\sigma_k}w)_{i_\alpha} (\cL_{\sigma_k}w)_{i_\beta}\dd t.
\end{aligned}    
\end{equation}

To collect these terms, let $z^1,\ldots,z^m$ be smooth vector fields and set $Z:=z^1\otimes\cdots\otimes z^m$, so that $Z_I=\prod_{\alpha=1}^m(z^\alpha)_{i_\alpha}$. Here superscripts label the vector fields, while subscripts denote their components. By definition, 
\begin{align*}
 (A_vZ)_I =\sum_{\alpha=1}^m\bigl[v\cdot\nabla(z^\alpha)_{i_\alpha}\bigr]
\prod_{\beta\ne\alpha}(z^\beta)_{i_\beta}, \quad 
(B_vZ)_I =\sum_{\alpha=1}^m\Big[\sum_{j=1}^d(\partial_{i_\alpha}v_j)(z^\alpha)_j\Big] \prod_{\beta\ne\alpha}(z^\beta)_{i_\beta}.
\end{align*}
For each $\alpha$, the sum of the two expressions in square  brackets is $(\cL_vz^\alpha)_{i_\alpha}$. Thus, 
\[
 L_v(z^1\otimes\cdots\otimes z^m)
 =\sum_{\alpha=1}^m z^1\otimes\cdots\otimes (\cL_vz^\alpha)\otimes\cdots\otimes z^m.
\]
Taking all factors equal to $w$, one sees that
\begin{equation*}
(L_vU)_I=\sum_{\alpha=1}^m(\cL_vw)_{i_\alpha}\prod_{\beta\ne\alpha}w_{i_\beta}.
\end{equation*}
Applying this rule once more, including to the factor $\cL_vw$, we obtain
\begin{align}\label{eq:L2vU}
 (L_v^2U)_I=\sum_{\alpha=1}^m(\cL_v^2w)_{i_\alpha} \prod_{\beta\ne\alpha}w_{i_\beta}+2\sum_{\alpha<\beta}(\cL_vw)_{i_\alpha}(\cL_vw)_{i_\beta} \prod_{\gamma\notin\{\alpha,\beta\}}w_{i_\gamma}.
\end{align}
Thus, substituting $v$ by $\sigma_k$ in \eqref{eq:L2vU}, we see that the two sums in $\frac12(L_{\sigma_k}^2U)_I$ are exactly the individual It\^o drift corrections and the quadratic covariation terms in \eqref{eq:dUI}. Combining \eqref{eq:w}--\eqref{eq:L2vU}, we get
\[
 \d U=\Big[-L_bU+\frac12\sum_kL_{\sigma_k}^{2}U \Big]\dd t
                       -\sum_kL_{\sigma_k}U\, \d W_t^k.
\]
Since the coefficients are deterministic, taking expectation closes the system: $V:=\bE U$ solves, for almost every $t$,
\begin{equation}\label{eq:t-pde}
 \partial_tV=-L_bV+\frac12\sum_{k=1}^d L_{\sigma_k}^{\,2}V,
 \quad V(s)=(\nabla\varphi)^{\otimes m}.
\end{equation}
The initial tensor has finite weighted $L^2$ norm for both weights defined in \eqref{def-weight-funct}, since
\[
 \|V_0\|_{2,\eta}^2
 =\int|\nabla\varphi|^{2m}\eta
 \leq\|\nabla\varphi\|_\infty^{2m}\int\eta<\infty,
 \quad \eta\in\{\rho,\rho_*\}.
\]

\subsubsection{Energy estimates}

Recall that $M=10^3q^2$, $q=d+1$, and $\Theta$ is defined in \eqref{def-notation-Theta}.

\begin{lemma}\label{lem:tensor}
There exists $\overline\delta=\overline\delta(d,\lambda,M)>0$ such that, for every $0\leq\delta<\overline\delta$, the following statements hold. For $1\leq m\leq M$, let $V\in C([s,T];L^2_\rho)\cap L^2(s,T;W^{1,2}_\rho)$ solve \eqref{eq:t-pde} in the sense of distributions, with arbitrary deterministic initial tensor $V_0\in L^2_\rho$. Then, for almost every $t\geq s$,
\begin{equation}\label{eq:t-en}
 \frac{\d}{\d t}\|V(t)\|_{2,\rho}^2
       +c\|\nabla V(t)\|_{2,\rho}^2
 \leq C(1+g(t))\|V(t)\|_{2,\rho}^2.
\end{equation}
Additionaly, if $V_0\in W^{1,2}_\rho$, then
\begin{equation}\label{eq:t-inc}
\begin{aligned}
 \|V(t)-V_0\|_{2,\rho}^2
 &\leq C(t-s)\|\nabla V_0\|_{2,\rho}^2+C\bigl(\Theta(t)-\Theta(s)\bigr)\|V_0\|_{2,\rho}^2.
\end{aligned}
\end{equation}
The same statements hold with $\rho$ replaced throughout by $\rho_*$.
\end{lemma}

\begin{proof}
Put $h(t):=1+g(t)$ and assume $\delta\leq1$. Since $|\nabla\log\rho|\leq C$, we obtain from Lemma~\ref{lem:wform}, for every tensor field $Z\in W^{1,2}_\rho$,
\begin{align}\label{eq:weight_coe}
\||b|Z\|_{2,\rho}^2+\||\nabla\sigma|Z\|_{2,\rho}^2\leq C\delta\|\nabla Z\|_{2,\rho}^2 +Ch(t)\|Z\|_{2,\rho}^2.
\end{align}

We first take smooth compactly supported initial tensors $V_0$. Testing \eqref{eq:t-pde} with $2V\rho$, we get
\begin{equation}\label{eq:ten-en}
 \frac{\d}{\d t}\|V\|_{2,\rho}^2
 =-2\langle V,L_bV\rangle_\rho
       +\sum_k\langle V,L_{\sigma_k}^2V\rangle_\rho.
\end{equation}

\emph{The drift term.} Recall that $L_v$ is given by \eqref{eq:def_abl}.The transport term satisfies 
\[
|\langle V,A_bV\rangle_\rho| \leq\int\rho|b|\, |V|\, |\nabla V|. 
\]
For each term in $\langle V,B_bV\rangle_\rho$, integration by parts yields that
\begin{equation}\label{eq-drift-IBP}
 \int\rho V_I(\partial_{i_\alpha}b_j)V_{I^{\alpha\to j}}=-\int\rho b_j\Bigl[
       (\partial_{i_\alpha}V_I)V_{I^{\alpha\to j}}
       +V_I\partial_{i_\alpha}V_{I^{\alpha\to j}} +(\partial_{i_\alpha}\log\rho)V_IV_{I^{\alpha\to j}}\Bigr].
\end{equation}
Summing over $I,\alpha,j$, we obtain
\[
 |\langle V,L_bV\rangle_\rho|
 \leq C\int\rho|b|\, |V|(|\nabla V|+|V|).
\]
Cauchy--Schwarz and \eqref{eq:weight_coe} yield that 
\begin{equation*}
\begin{aligned}
& |\langle V,L_bV\rangle_\rho|
  \leq C\||b|V\|_{2,\rho}(\|\nabla V\|_{2,\rho}+\|V\|_{2,\rho})\\
 \leq{} & C(\sqrt\delta\|\nabla V\|_{2,\rho}+\sqrt h\|V\|_{2,\rho})(\|\nabla V\|_{2,\rho}+\|V\|_{2,\rho})\\
 \leq{} & \left[\eps+C(\sqrt\delta+\delta)\right]\|\nabla V\|_{2,\rho}^2+C_{\eps}h\|V\|_{2,\rho}^2.
\end{aligned}
\end{equation*}

\emph{The diffusion term.} Put $A_k:=A_{\sigma_k}$, $B_k:=B_{\sigma_k}$, and $c_k:=\div\sigma_k+\sigma_k\cdot\nabla\log\rho$. For tensor fields $U,Z$, we have
\[
 \langle U,A_kZ\rangle_\rho
 =-\langle A_kU,Z\rangle_\rho-\int\rho c_k\,U\cdot Z.
\]
So the adjoint operator of $A_k$ is $A_k^{*,\rho}=-A_k-c_k$. For a smooth vector field $v$, let $B_v^\top$ denote the pointwise transpose on the tensor product space $(\mR^d)^{\otimes m}$, i.e.
\[
 (B_v^\top Z)_I=\sum_{\alpha=1}^m\sum_{j=1}^d(\partial_jv_{i_\alpha})Z_{I^{\alpha\to j}}.
\]
For fixed $\alpha$ and change the summation variables by setting $J:=I^{\alpha\to j}$ and $\ell:=i_\alpha$. If $j_\alpha$ denotes the $\alpha$-th entry of $J$, then $j_\alpha=j$ and $I=J^{\alpha\to\ell}$. Thus, $(I,j)\mapsto(J,\ell)$ is a bijection of $\{1,\ldots,d\}^m\times\{1,\ldots,d\}$, and the factors $U_I$, $\partial_{i_\alpha}v_j$, and $Z_{I^{\alpha\to j}}$ become $U_{J^{\alpha\to\ell}}$, $\partial_\ell v_{j_\alpha}$, and $Z_J$, respectively. Thus, for each fixed $\alpha$, we have 
\[
\sum_{I,j}\int\rho U_I(\partial_{i_\alpha}v_j)Z_{I^{\alpha\to j}}=\sum_{J,\ell}\int\rho U_{J^{\alpha\to\ell}}(\partial_\ell v_{j_\alpha})Z_J =\sum_{I,j}\int\rho(\partial_jv_{i_\alpha})U_{I^{\alpha\to j}}Z_I
\]
Therefore, $\langle U,B_k Z\rangle_\rho=\langle B_k^\top U,Z\rangle_\rho$, i.e. $B_k^{*,\rho}=B_k^\top$, and the adjoint of $L_{\sigma_k}=A_k+B_k$ in the weighted inner product is
\[
 L_{\sigma_k}^{*,\rho}
 =-A_{\sigma_k}-c_k+B_{\sigma_k}^{\top}.
\]
Therefore, we have
\begin{equation*}
\begin{aligned}
 \langle V,L_{\sigma_k}^2V\rangle_\rho
 & =\langle -A_kV-c_kV+B_k^\top V,\,A_kV+B_kV\rangle_\rho =-\|A_kV\|_{2,\rho}^2\\
 &\quad +\underbrace{\langle(B_k^\top-B_k)V,A_kV\rangle_\rho-\langle c_kV,A_kV\rangle_\rho +\langle B_k^\top V,B_kV\rangle_\rho-\langle c_kV,B_kV\rangle_\rho}_{=:\cR_k}.
\end{aligned}   
\end{equation*}

Boundedness of $\sigma$ and the definition of $B_k$ imply
\[
 |A_kV|\leq C|\nabla V|,\quad
 |B_kV|+|B_k^\top V|\leq C|\nabla\sigma_k|\, |V|,\quad
 |c_k|\leq C(1+|\nabla\sigma_k|).
\]
Therefore,
\[
 \sum_k|\mathcal R_k|
 \leq C\int\rho\left[
 (1+|\nabla\sigma|)|V|\, |\nabla V|
              +(1+|\nabla\sigma|+|\nabla\sigma|^2)|V|^2\right],
\]
and no second derivative of $\sigma$ occurs. Using the Cauchy inequality and weighted form bound \eqref{eq:wform}, we now estimate
\begin{align*}
 & \int\rho(1+|\nabla\sigma|)|V|\, |\nabla V|
  \leq \|\nabla V\|_{2,\rho}
        (\|V\|_{2,\rho}+\||\nabla\sigma|V\|_{2,\rho})\\
 \leq{} & (\eps+C\sqrt\delta)\|\nabla V\|_{2,\rho}^2
            +C_{\eps} h\|V\|_{2,\rho}^2,\\
 & \int\rho(1+|\nabla\sigma|+|\nabla\sigma|^2)|V|^2
  \leq C(\|V\|_{2,\rho}^2+\||\nabla\sigma|V\|_{2,\rho}^2) \leq C\delta\|\nabla V\|_{2,\rho}^2+Ch\|V\|_{2,\rho}^2.
\end{align*}
So we get 
\begin{equation*}
\sum_k|\cR_k|\leq\left[\eps+C(\sqrt\delta+\delta)\right]\|\nabla V\|_{2,\rho}^2+C_{\eps}h\|V\|_{2,\rho}^2.
\end{equation*}
On the other hand, ellipticity yields
\begin{equation}\label{eq:ten-ell}
\sum_k\|A_kV\|_{2,\rho}^2
 =\sum_I\int\rho\sum_k|\sigma_k\cdot\nabla V_I|^2
 \geq2\lambda_0\|\nabla V\|_{2,\rho}^2,\quad\lambda_0\geq\lambda/2.   
\end{equation}
Combining \eqref{eq:ten-en} and \eqref{eq:ten-ell}, we obtain
\[
 \frac{\d}{\d t}\|V\|_{2,\rho}^2
 +\left[2\lambda_0-3\eps-C(\sqrt\delta+\delta)\right]
       \|\nabla V\|_{2,\rho}^2
 \leq C_{\eps}h\|V\|_{2,\rho}^2.
\]
Choose $\eps=\eps(d,\lambda)>0$ and then $\overline\delta=\overline\delta(d,\lambda,M)\in(0,1]$ small enough that the bracket is bounded below by a positive constant $c>0$ for all $m\leq M$ and $\delta<\overline\delta$. We thus obtain \eqref{eq:t-en}, with constants $c,C$ depending only on $d,\lambda,M$. 

\emph{The increment.} Put $F:=V-V_0$, so $F(s)=0$ and
\[
 \partial_tF=-L_bF+\frac12\sum_kL_{\sigma_k}^2F+H_0,
 \quad H_0:=-L_bV_0+\frac12\sum_kL_{\sigma_k}^2V_0.
\]
Testing the equation for $F$ with $2F\rho$ and using the energy estimate for the homogeneous terms, we obtain
\begin{equation}\label{eq:inc-energy}
 \frac{\d}{\d t}\|F\|_{2,\rho}^2+c\|\nabla F\|_{2,\rho}^2
 \leq Ch\|F\|_{2,\rho}^2+2\langle F,H_0\rangle_\rho.
\end{equation}
It remains to bound $|\langle F,H_0\rangle_\rho|$, keeping at most one spatial derivative on $V_0$. For the drift term, integrating by parts in $B_bV_0$ (cf. \eqref{eq-drift-IBP}), we obtain
\begin{align*}
 \langle F,-L_bV_0\rangle_\rho ={} & -\langle F,A_bV_0\rangle_\rho +\sum_{I,\alpha,j}\int\rho b_j\Bigl[(\partial_{i_\alpha}F_I)(V_0)_{I^{\alpha\to j}}\\
 &+F_I\partial_{i_\alpha}(V_0)_{I^{\alpha\to j}}+(\partial_{i_\alpha}\log\rho)F_I(V_0)_{I^{\alpha\to j}}\Bigr].
\end{align*}
By Cauchy--Schwarz and \eqref{eq:weight_coe}, applied to $F$ and $V_0$, we obtain
\begin{align*}
 & |\langle F,-L_bV_0\rangle_\rho|
  \leq C\||b|F\|_{2,\rho}\|\nabla V_0\|_{2,\rho}+C\||b|V_0\|_{2,\rho}\bigl(\|\nabla F\|_{2,\rho}+\|F\|_{2,\rho}\bigr)\\
 \leq{} & C\bigl(\sqrt\delta\|\nabla F\|_{2,\rho}+\sqrt h\|F\|_{2,\rho}\bigr)\|\nabla V_0\|_{2,\rho}\\
 &+C\bigl(\sqrt\delta\|\nabla V_0\|_{2,\rho}+\sqrt h\|V_0\|_{2,\rho}\bigr) \bigl(\|\nabla F\|_{2,\rho}+\|F\|_{2,\rho}\bigr)\\
 \leq{} & C\bigl(\|\nabla F\|_{2,\rho}+\sqrt h\|F\|_{2,\rho}\bigr) \bigl(\|\nabla V_0\|_{2,\rho}+\sqrt h\|V_0\|_{2,\rho}\bigr),
\end{align*}
where the last step uses $\delta\leq1$ and $h\geq1$.
For the diffusion forcing, one integration by parts yields 
\begin{align*}
 & |\langle F,L_{\sigma_k}^2V_0\rangle_\rho|
  = |\langle L_{\sigma_k}^{*,\rho}F,L_{\sigma_k}V_0\rangle_\rho|\\
 ={} & |\langle-A_kF-c_kF+B_k^\top F,\,A_kV_0+B_kV_0\rangle_\rho|\\
 \leq{} & \bigl(\|A_kF\|_{2,\rho}+\|c_kF\|_{2,\rho}+\|B_k^\top F\|_{2,\rho}\bigr)\bigl(\|A_kV_0\|_{2,\rho}+\|B_kV_0\|_{2,\rho}\bigr)\\
 \leq{} & C\bigl(\|\nabla F\|_{2,\rho}+\sqrt h\|F\|_{2,\rho}\bigr) \bigl(\|\nabla V_0\|_{2,\rho}+\sqrt h\|V_0\|_{2,\rho}\bigr).
\end{align*}
Combining the drift and diffusion estimates, we obtain
\begin{align*}
 |\langle F,H_0\rangle_\rho|\leq\eps\|\nabla F\|_{2,\rho}^2
       +C_{\eps}h\|F\|_{2,\rho}^2+C_{\eps}
       (\|\nabla V_0\|_{2,\rho}^2+h\|V_0\|_{2,\rho}^2). 
\end{align*}
Choosing $0<\eps<c/4$ and substituting this bound into \eqref{eq:inc-energy}, we obtain 
\begin{align*}
 \frac{\d}{\d t}\|F\|_{2,\rho}^2+ \frac{c}{2} \|\nabla F\|_{2,\rho}^2
 \leq Ch\|F\|_{2,\rho}^2+C\|\nabla V_0\|_{2,\rho}^2+Ch\|V_0\|_{2,\rho}^2.
\end{align*}
Since $F(s)=0$, we obtain by Gronwall's inequality
\begin{align*}
 & \|F(t)\|_{2,\rho}^2
  \leq C\int_s^t\bigl(
       \|\nabla V_0\|_{2,\rho}^2+h(r)\|V_0\|_{2,\rho}^2\bigr)\dd r\\
 ={} & C(t-s)\|\nabla V_0\|_{2,\rho}^2
       +C\bigl(\Theta(t)-\Theta(s)\bigr)\|V_0\|_{2,\rho}^2, 
\end{align*}
where $C$ only depends on $d,\lambda,M,T$ and $\|g\|_{L^1(0,T)}$. This proves \eqref{eq:t-inc}. 

By density of smooth compactly supported initial tensors, the estimates extend to $L^2_\rho$ and $W^{1,2}_\rho$ initial data, respectively. Both weights $\rho$ and $\rho_*$ are integrable with bounded logarithmic gradient, so all the calculations apply to $\rho_*$. The change of variables $r\mapsto t-r$ preserves the form and ellipticity bounds and the relevant time integral of $1+g$, so the same conclusions hold for the coefficients $\overline b,\overline\sigma$ at the end of the proof of Lemma~\ref{lem:occ}. 
\end{proof}

\subsection{Spatial and Malliavin derivative estimates}
We now use the transport equations and energy estimates to control the Jacobian matrices of the flows and their Malliavin derivatives. The bounds are uniform in the approximation index $n$.

\subsubsection{Moment bounds for Jacobians}
\begin{proposition}\label{prop:jac}
There exists $\overline\delta=\overline\delta(d,\lambda,M)>0$ such that, for every $0\leq\delta<\overline\delta$ and every even integer $2\leq m\leq M$,
\begin{align}\label{eq:jac}
 \sup_{n,s,t}\int\bE\bigl(|J^n_{s,t}(x)|^m+|K^n_{s,t}(x)|^m\bigr)\rho(x)\dd x\leq C_m,
\end{align}
\begin{align}\label{eq:jac-inc}
\int\bE\bigl(|J^n_{s,t}(x)-I|^m+|K^n_{s,t}(x)-I|^m\bigr)\rho(x)\dd x \leq C_m\bigl(\Theta(t)-\Theta(s)\bigr)^{\frac{1}{2}}.
\end{align}
Here $C_m$ depends only on $d,\lambda,m,T$ and $\|g\|_{L^1(0,T)}$. These bounds also hold with $\rho_*$.
\end{proposition}
\begin{proof}
We first estimate $K=\nabla Y$. Fix $j\in\{1,\ldots,d\}$ and now take $\varphi(x):=x_j$. Then $u(t,y)=Y^j_{s,t}(y)$ by \eqref{def-u}, and the initial tensor in \eqref{eq:t-pde} is $V_0=e_j^{\otimes m}$. For the matrix entry $a:=\partial_iY^j=K_{ji}$, the tensor component with $I=(i,\ldots,i)$ equals $\bE a^m$. Since $m$ is even, we obtain by Cauchy--Schwarz and \eqref{eq:t-en} that
\[
 \int\bE|a|^m\rho
 \leq\Big(\int\rho\Big)^{\frac{1}{2}}
             \|V_I(t)\|_{2,\rho}\leq C_m.
\]
Summing over the matrix entries, we obtain the bound for $K$.

For \eqref{eq:jac-inc}, put $a_0:=\delta_{ij}$. The binomial identity 
\[
 \bE|a-a_0|^m
 =\sum_{\ell=1}^m\binom{m}{\ell}(-a_0)^{m-\ell} (\bE a^\ell-a_0^\ell)
\]
reduces the problem to the moment equations of orders $1,\ldots,m$. For each such order $\ell$, let $V^{(\ell)}$ denote its moment tensor. Then $\bE a^\ell=V^{(\ell)}_{(i,\ldots,i)}$ and the corresponding initial component is $a_0^\ell$. Since $e_j^{\otimes\ell}$ is constant, we deduce from \eqref{eq:t-inc} and Cauchy--Schwarz that
\[
 \int|\bE a^\ell-a_0^\ell|\rho
 \leq\Big(\int\rho\Big)^{\frac{1}{2}}
           \|V^{(\ell)}-e_j^{\otimes\ell}\|_{2,\rho}
 \leq C_\ell\bigl(\Theta(t)-\Theta(s)\bigr)^{\frac{1}{2}}.
\]
Plugging these bounds into the binomial identity and summing over the matrix entries, we obtain \eqref{eq:jac-inc} for $K$.

Finally, $X_{s,t}$ is the inverse of the flow with coefficients $\overline b,\overline\sigma$ from Lemma~\ref{lem:occ}. The driving Brownian motion $\overline W_r=W_t-W_{t-r}$ is considered with respect to its natural filtration. Applying the bounds to the spatial derivatives of the inverse flow of this flow proves the statements for $J=\nabla X$. The corresponding time integral is $\int_0^{t-s}(1+g(t-r))\dd r=\Theta(t)-\Theta(s)$.
\end{proof}

\subsubsection{Estimates for Malliavin derivatives}
We will express the Malliavin derivative through the Jacobian and $\sigma$, and then use the spatial derivative bounds to control Malliavin derivatives and joint regularity in $(s,t,x)$. Since the Malliavin derivative contains Jacobians evaluated at random positions, we first record a change-of-variables estimate. It allows us to use the spatially integrated moment bounds even when the field and the flow depend on the same Brownian path.

\begin{lemma}\label{lem:comp}
There exists $\overline\delta=\overline\delta(d,\lambda)>0$ such that, for every $0\leq\delta<\overline\delta$, $Z=X^n_{s,t}$ or $Y^n_{s,t}$, and any random measurable field $F$,
\begin{equation}\label{eq:comp}
 \left(\bE\int |F(Z(x))|^p\rho(x)\dd x\right)^{\frac{1}{p}}
 \leq C_p\left(\bE\int|F(y)|^{3p}\rho(y)\dd y\right)^{\frac{1}{3p}},
 \quad p\geq1.
\end{equation}
Here $C_p$ depends only on $d,\lambda,p,T$ and $\|g\|_{L^1(0,T)}$, and the field $F$ need not be independent of $Z$.
\end{lemma}
\begin{proof}
Put $U:=Z^{-1}$ and change variables $y=Z(x)$. Before taking the $p$-th root, the left-hand side is
\[
 \bE\int |F(y)|^p\rho(y)
        \frac{\rho(U(y))}{\rho(y)}|\det \nabla U(y)|\dd y.
\]
The weight ratio satisfies
\[
 \frac{\rho(U(y))}{\rho(y)}
       \leq C(1+|U(y)-y|)^{2q}.
\]
Since $U$ is either $X^n_{s,t}$ or $Y^n_{s,t}$, we apply the displacement estimate \eqref{eq:disp-u} with exponent $6q$ to obtain
\[
 \bE\int\left(\frac{\rho(U(y))}{\rho(y)}\right)^3\rho(y)\dd y
 \leq C\bE\int(1+|U(y)-y|)^{6q}\rho(y)\dd y
 \leq C.
\]
To control the determinant, let $m_*:=2\lceil 3d/2\rceil$. Since $m_*\leq M$, \eqref{eq:jac} yields
\[
 \bE\int|\det\nabla U(y)|^3\rho(y)\dd y
 \leq C\bE\int\bigl(1+|\nabla U(y)|^{m_*}\bigr)\rho(y)\dd y
 \leq C.
\]
Applying H\"older's inequality, we obtain \eqref{eq:comp}. 
\end{proof}

We next establish a moment bound for the Malliavin derivative and continuity under translations of its time variable.

On the fixed $d$-dimensional Wiener space, let $E$ be a finite-dimensional Euclidean space. For $F:\Omega\to E$, denote its Malliavin derivatives in the directions of $W^k$ by $\mathrm D_r^kF$, $1\leq k\leq d$.

It is well known that the Malliavin derivative of $X^n_{s,t}(x)$ satisfies
\begin{equation}\label{eq:mall}
 \mathrm D_rX^n_{s,t}(x)
=J^n_{r,t}(X^n_{s,r}(x))\sigma^n(r,X^n_{s,r}(x))\1_{r\in[s,t]},
\end{equation}
for $(\dd r\times\mP\times\dd x)$-almost every $(r,\omega,x)$ (see for instance \cite{N06malliavin}).

\begin{proposition}\label{prop:mall}
There exists $\overline\delta=\overline\delta(d,\lambda)>0$ such that, for every $0\leq\delta<\overline\delta$ and every fixed $0\leq s\leq t\leq T$, it holds that
\begin{align}\label{eq:mall-bnd}
 &\sup_n\esssup_{r\in[0,T]}
      \bE\int |\mathrm D_rX^n_{s,t}(x)|^4\rho(x)\dd x<\infty,\\
 \label{eq:mall-sh}
 &\lim_{h\downarrow0}\sup_n\int_0^{T-h}
    \bE\int\left|\mathrm D_{r+h}X^n_{s,t}(x)
                         -\mathrm D_rX^n_{s,t}(x)\right|^2
                                    \rho(x)\dd x\dd r=0.
\end{align}
The same statements hold for $Y^n_{s,t}$.
\end{proposition}
\begin{proof}
By \eqref{eq:mall} and the uniform boundedness of $\sigma^n$, applying Lemma~\ref{lem:comp} with $p=4$, $F=J^n_{r,t}$, and $Z=X^n_{s,r}$, one sees that 
\[
 \bE\int |\mathrm D_rX^n_{s,t}(x)|^4\rho(x)\dd x
 \leq C\left(\bE\int |J^n_{r,t}(y)|^{12}\rho(y)\dd y\right)^{\frac{1}{3}}
 \leq C.
\]
The last inequality follows from \eqref{eq:jac}, with $C$ depending only on $d,\lambda,T$ and $\|g\|_{L^1(0,T)}$, and proves \eqref{eq:mall-bnd}.

For \eqref{eq:mall-sh}, the case $s=t$ is immediate. Fix $0<h<t-s$ and put $X_r:=X^n_{s,r}$. For $s\leq r<r+h\leq t$, applying the chain rule to the identity $X^n_{r,t}=X^n_{r+h,t}\circ X^n_{r,r+h}$, we obtain $J^n_{r,t}= J^n_{r+h,t}(X^n_{r,r+h}) J^n_{r,r+h}$, hence 
  $$J^n_{r,t}(X_r) = J^n_{r+h,t}(X_{r+h}) J^n_{r,r+h}(X_r). $$
Recalling \eqref{eq:mall}, we have
\begin{equation}\label{eq:mall-dif}
\begin{aligned}
&\mathrm D_rX^n_{s,t}-\mathrm D_{r+h}X^n_{s,t}\\
 ={} & J^n_{r+h,t}(X_{r+h}) \bigl\{\underbrace{(J^n_{r,r+h}(X_r)-I)\sigma^n(r,X_r)}_{=:I_1(r)}+\underbrace{ [\sigma^n(r,X_r)-\sigma^n(r+h,X_r)]}_{=:I_2(r)}\\
 &\quad +\underbrace{[\sigma^n(r+h,X_r)-\sigma^n(r+h,X_{r+h})]}_{=:I_3(r)}\bigr\}.
\end{aligned}
\end{equation}
The factor $J^n_{r+h,t}(X_{r+h})$ in \eqref{eq:mall-dif} is uniformly bounded in $L^4(\mP\times\rho)$ by Lemma \ref{lem:comp} and \eqref{eq:jac}. By H\"older's inequality and \eqref{eq:mall-bnd}, we therefore have
\begin{equation}\label{eq:mall-sum}
\begin{aligned}
\int_0^{T-h}
\|\mathrm D_{r+h}X^n_{s,t}-\mathrm D_rX^n_{s,t}\|_{L^2(\mP\times\rho)}^2
\dd r \leq Ch+C \sum_{k=1}^3 \int_s^{t-h}
\|I_k(r)\|_{L^4(\mP\times\rho)}^2 \dd r.
\end{aligned}
\end{equation}
Here $Ch$ bounds the contribution from $([s-h,s]\cup[t-h,t])\cap[0,T-h]$, where one of the two Malliavin derivatives vanishes.

For $I_1$, Lemma \ref{lem:comp} and \eqref{eq:jac-inc} yield that
\begin{align*}
 \|I_1(r)\|_{L^4(\mP\times\rho)}^2 \leq C\left(\bE\int|J^n_{r,r+h}(y)-I|^{12}\rho(y)\dd y\right)^{\frac{1}{6}} \leq C\bigl(\Theta(r+h)-\Theta(r)\bigr)^{\frac{1}{12}}.
\end{align*}

For $I_2$, by \eqref{eq:comp} with $p=4$ and the uniform boundedness of $\sigma^n$,
\begin{equation*}
\begin{aligned}
 & \int_{s}^{t-h}\|I_2(r)\|_{L^4(\mP\times\rho)}^2 \dd r 
  \leq C \int_{0}^{T-h}\left(\int
 |\sigma^n(r+h,y)-\sigma^n(r,y)|^{12}\rho(y)\dd y\right)^{\frac{1}{6}} \dd r \\
 \leq{} & C a(h)^{\frac{1}{6}}, 
\end{aligned}
\end{equation*}
where 
\[
 a(h):=\sup_n\int_0^{T-h}\int
       |\sigma^n(r+h,y)-\sigma^n(r,y)|^2\rho(y)\dd y\dd r.
\]
Since $\rho$ is integrable and $\sigma^n$ is uniformly bounded, $\sigma^n\to\sigma$ in $\cH:=L^2(0,T;L^2_\rho)$. For $0<h<T$, put $\cH_h:=L^2(0,T-h;L^2_\rho)$. By the triangle inequality, we have
\[
 \|\sigma^n(h+\cdot)-\sigma^n\|_{\cH_h}
 \leq2\|\sigma^n-\sigma\|_{\cH}
       +\|\sigma(h+\cdot)-\sigma\|_{\cH_h}.
\]
The first term is uniformly small for all sufficiently large $n$, and the second tends to zero as $h\downarrow0$ by continuity of time translations in $L^2$. The same continuity applies to the finitely many remaining $\sigma^n$, so $a(h)\to0$ as $h\downarrow0$.

For $I_3$, let $\sigma^n_\eps:=\sigma^n*_x\eta_\eps$, where $\eta_\eps(x):=\eps^{-d}\eta(x/\eps)$, $0<\eps\leq1$, and $\eta$ is the mollifier from Proposition~\ref{prop:approx}. Using $\rho(y+z)\leq C\rho(y)$ for $|z|\leq1$ and \eqref{eq:loc-src}, we have
\begin{align*}
 \sup_n\int_0^T\int|\sigma^n-\sigma^n_\eps|^2\rho\dd r
 &\leq C\eps^2\sup_n\int_0^T\int|\nabla\sigma^n|^2\rho\dd r
 \leq C\eps^2,\\
 \sup_n\|\nabla\sigma^n_\eps\|_\infty
 &\leq\sup_n\|\sigma^n\|_\infty\|\nabla\eta_\eps\|_1
 \leq C\eps^{-1}.
\end{align*}

By \eqref{eq:comp} and \eqref{eq:disp} with $p=4$, and H\"older's inequality in time,
\begin{equation}\label{eq:mall_3}
\begin{aligned}
 &\int_s^{t-h}\|I_3(r)\|_{L^4(\mP\times\rho)}^2\dd r\\
 \leq{} & C\sum_{\theta\in\{0,h\}}\int_s^{t-h}\| (\sigma^n-\sigma^n_\eps)(r+h,X_{r+\theta})\|_{L^4(\mP\times\rho)}^2\dd r\\
 &\quad+C\int_s^{t-h}\|\sigma^n_\eps(r+h,X_r)-\sigma^n_\eps(r+h,X_{r+h})\|_{L^4(\mP\times\rho)}^2\dd r\\
 \leq{} & C\left(\int_0^T\int|\sigma^n-\sigma^n_\eps|^2\rho\dd r\right)^{\frac{1}{6}}
 +C\|\nabla\sigma^n_\eps\|_\infty^2h \leq C\bigl(\eps^{\frac{1}{3}}+\eps^{-2}h\bigr).
\end{aligned}
\end{equation}

Combining \eqref{eq:mall-sum}--\eqref{eq:mall_3}, we get
\begin{align*}
&\sup_n\int_0^{T-h}
\|\mathrm D_{r+h}X^n_{s,t}-\mathrm D_rX^n_{s,t}\|_{L^2(\mP\times\rho)}^2
\dd r\\
\leq{} & Ch+C\int_s^{t-h}
\bigl(\Theta(r+h)-\Theta(r)\bigr)^{\frac{1}{12}}\dd r
+C\bigl(a(h)^{\frac{1}{6}}+\eps^{\frac{1}{3}}+\eps^{-2}h\bigr), 
\end{align*}
where $C$ depends only on $d,\lambda,T$, $\|g\|_{L^1(0,T)}$ and the fixed mollifier $\eta$. Since $\Theta$ is absolutely continuous, first let $h\downarrow0$ with $\eps$ fixed, and then let $\eps\downarrow0$, proving \eqref{eq:mall-sh}.

For the inverse flow $Y^n_{s,t}$. By Lemma~\ref{lem:occ}, $Y^n_{s,t}(W)=\overline X^n_{0,t-s}(\cR_tW)$, where $(\cR_tW)(a):=W_t-W_{t-a}$ for $0\leq a\leq t$. Differentiating this identity in Cameron--Martin directions, we find
\[
 \mathrm D_rY^n_{s,t}(W)
   =(\mathrm D_{t-r}\overline X^n_{0,t-s})(\cR_tW),
 \quad s\leq r\leq t;
\]
the identity holds almost everywhere, and the derivative vanishes on $[0,T]\setminus[s,t]$. The map $\cR_t$ preserves Wiener measure on $[0,t]$, and the coefficients $\overline b,\overline\sigma$ satisfy the same bounds and convergence assumptions, so the forward estimate implies \eqref{eq:mall-bnd}. Changing variables by $a=t-r-h$ on $[s,t-h]$, we obtain \eqref{eq:mall-sh}.
\end{proof}

\subsection{Joint regularity}
We combine spatial regularity with the moment bounds for time increments. The flow identities
\[
 X^n_{s,t}=X^n_{u,t}\circ X^n_{s,u},\quad
 Y^n_{s,t}=Y^n_{s,u}\circ Y^n_{u,t},\quad s\leq u\leq t,
\]
allow us to estimate changes in the initial time of $X^n$ and the terminal time of $Y^n$ using spatial H\"older bounds and \eqref{eq:disp}.

\begin{proposition}\label{prop:holder}
Put $\Delta_T:=\{(s,t):0\leq s\leq t\leq T\}$ and $r:=M/8$. There exists $\overline\delta=\overline\delta(d,\lambda)>0$ such that, for every $0\leq\delta<\overline\delta$ and every $R<\infty$, the $r$-th moments of the $C^{\frac{1}{8}}(\Delta_T\times B_R)$ norms of $X^n,Y^n$ are bounded uniformly in $n$. 
\end{proposition}
\begin{proof}
Fix $R\geq1$. Since $r=M/8=125q^2$ and $1-d/(4r)>1/2$, Morrey's inequality, together with comparison of the weight $\rho$ on balls implies that for every smooth map $Z$,
\begin{equation}\label{eq:morrey}
 |Z(x)-Z(y)|
 \leq C\|\nabla Z\|_{L^{4r}_\rho}
       (1+|x|+|y|)^L|x-y|^{\frac{1}{2}},
\end{equation}
where $C,L$ depend only on $d,r$. By \eqref{eq:jac}, for $Z=X^n$ or $Y^n$ and $x,y\in B_R$,
\[
 \bE|Z_{s,t}(x)-Z_{s,t}(y)|^r
 \leq C_R\left(\bE\int|\nabla Z_{s,t}|^{4r}\rho\right)^{\frac{1}{4}}
             |x-y|^{\frac{r}{2}}
 \leq C_R|x-y|^{\frac{r}{2}}.
\]

For $s<s'\leq t$ and $x\in B_R$, put $Y:=X^n_{s,s'}(x)$ and $Z:=X^n_{s',t}$. Using the flow identity and \eqref{eq:morrey}, we estimate
\[
\begin{aligned}
 & |X^n_{s,t}(x)-X^n_{s',t}(x)|^r
  = |Z(Y)-Z(x)|^r\\
 \leq{} & C_R\|\nabla Z\|_{L^{4r}_\rho}^r
       (1+|Y-x|)^{Lr}|Y-x|^{\frac{r}{2}}.
\end{aligned}
\]
By \eqref{eq:jac} with $4r=M/2$ and \eqref{eq:disp-u}, the first two factors on the right are uniformly bounded in $L^2(\mP)$ and $L^4(\mP)$, respectively. H\"older's inequality and \eqref{eq:disp}, therefore yield
\begin{align*}
 & \bE|X^n_{s,t}(x)-X^n_{s',t}(x)|^r
  \leq C_R\left(\bE\|\nabla Z\|_{L^{4r}_\rho}^{2r}\right)^{\frac{1}{2}}
       \left(\bE(1+|Y-x|)^{4Lr}\right)^{\frac{1}{4}}
       \left(\bE|Y-x|^{2r}\right)^{\frac{1}{4}}\\
 \leq{} & C_R\left(\bE|Y-x|^{2r}\right)^{\frac{1}{4}}
 \leq C_R|s'-s|^{\frac{r}{4}}.
\end{align*}

Using $Y^n_{s,t'}=Y^n_{s,t}\circ Y^n_{t,t'}$ and repeating the argument, we obtain, for $s\leq t<t'$,
\begin{align*}
 \bE|Y^n_{s,t'}(x)-Y^n_{s,t}(x)|^r \leq  C_R|t'-t|^{\frac{r}{4}}.
\end{align*}
The other time increments follow directly from \eqref{eq:disp}:
\begin{align*}
 \bE|X^n_{s,t'}(x)-X^n_{s,t}(x)|^r
 &\leq C|t'-t|^{\frac{r}{2}},\quad s\leq t\leq t',\\
 \bE|Y^n_{s,t}(x)-Y^n_{s',t}(x)|^r
 &\leq C|s'-s|^{\frac{r}{2}},\quad s\leq s'\leq t.
\end{align*}

Combining all the estimates above, we obtain, for $x,y\in B_R$,
\begin{align*}
 \bE\bigl(|X^n_{s,t}(x)-X^n_{s',t'}(y)|^r+|Y^n_{s,t}(x)-Y^n_{s',t'}(y)|^r\bigr) \leq C_R\bigl(|s-s'|+|t-t'|+|x-y|\bigr)^{\frac{r}{4}}.
\end{align*}
Here $C_R$ depends only on $d,\lambda,T,R$ and $\|g\|_{L^1(0,T)}$, since $r=M/8$ is fixed by $d$. Noticing that $\frac14-\frac{d+2}{r} =\frac14-\frac{q+1}{125q^2}>\frac18$, we conclude by Kolmogorov's criterion that
\[
 \sup_n\bE\Bigl(\|X^n\|_{C^{\frac{1}{8}}(\Delta_T\times B_R)}^r
              +\|Y^n\|_{C^{\frac{1}{8}}(\Delta_T\times B_R)}^r\Bigr)<\infty.
\]
The joint tightness of $X^n,Y^n$ follows.
\end{proof}
\section{Compactness and uniqueness of the limiting flow}\label{sec:compact}

We now complete the proof of Theorem~\ref{thm:flow}. First, an extended compactness criterion for Wiener functionals turns the estimates of Section~\ref{sec:deriv} into strong $L^2$ compactness for each fixed pair $(s,t)$. By the joint H\"older bounds, we then obtain locally uniform subsequential convergence in probability, jointly for the flow and its inverse, on the original Wiener space. The Krylov-type estimate then permits passage to the equation at every starting point. Finally, the backward equation implies uniqueness in law within $\sK_2(s,x)$; reconstructing the Brownian motion from the solution path turns this into pathwise uniqueness in the same class and ensures convergence of the whole sequence.

\subsection{A compactness criterion on Wiener space}

The next lemma extends the compactness criterion for Wiener functionals of \cite{DMN92compact} to spatially dependent random fields in a weighted $L^2$ space. We state it with separate spatial and noise dimensions. In the application to $X^n_{s,t}-\Id$ and $Y^n_{s,t}-\Id$ at each fixed pair $(s,t)$, the noise dimension is $m=d$ and the fields take values in $E=\mR^d$. Recall that the weight functions $\rho$ and $\rho_\ast$ are defined in \eqref{def-weight-funct}.

\begin{lemma}\label{lem:compact}
Let $d,m\geq1$, $T>0$, and let $E$ be a finite-dimensional Euclidean space. Let $\Omega=C_0([0,T];\mR^m)$, equipped with Wiener measure $\mP$ and the completed Borel sigma field $\sF$, and let $W$ be the canonical Brownian motion. Let $F_n:\Omega\times\mR^d\to E$ be jointly measurable random fields that belong to the Malliavin Sobolev space $L^2_\rho(\mR^d;\mD^{1,2})$ and have weak spatial gradients $\nabla F_n$. Assume that
\[
 \sup_n\left(
 \bE\int|F_n|^2\rho_*
 +\bE\int|\nabla F_n|^2\rho
 +\bE\int_0^T\int|\mathrm D_rF_n|^2\rho\dd r
 \right)<\infty,
\]
and
\[
 \lim_{h\downarrow0}\sup_n\int_0^{T-h}\bE\int
   |\mathrm D_{r+h}F_n-\mathrm D_rF_n|^2\rho\dd r=0.
\]
Then $(F_n)_n$ is relatively compact in the strong topology of $L^2(\mP\times\rho;E)$.
\end{lemma}
\begin{proof}
By It\^o--Wiener chaos expansion, we have 
\[
 F_n=\sum_{\ell=0}^\infty I_\ell(f_{n,\ell})
 \quad\mbox{in }L^2(\mP\times\rho;E),
\]
where $I_\ell$ is the $\ell$-fold It\^o--Wiener integral with respect to $W$, acting on the noise variables, and $I_0(f)=f$. 
The kernels are symmetric in the pairs $(t_j,k_j)$: for every permutation $\pi\in\mathfrak S_\ell$,
\[
 f_{n,\ell}(\boldsymbol t,\boldsymbol k,x)
 =f_{n,\ell}(\boldsymbol t_\pi,\boldsymbol k_\pi,x)
 \quad\mbox{a.e.},
\]
where $\boldsymbol t_\pi=(t_{\pi(1)},\ldots,t_{\pi(\ell)})$ and $\boldsymbol k_\pi=(k_{\pi(1)},\ldots,k_{\pi(\ell)})$.
For $\boldsymbol t=(t_1,\ldots,t_\ell)$ and $\boldsymbol k=(k_1,\ldots,k_\ell)$, define
\[
 \|f\|_\ell^2:=
 \sum_{k_1,\ldots,k_\ell=1}^m
 \int_{[0,T]^\ell}\!\int
 |f(\boldsymbol t,\boldsymbol k,x)|^2\rho(x)\dd x\,\dd\boldsymbol t.
\]
For $\ell=0$, this is $\|f\|_0^2=\int|f|^2\rho$. The chaos isometry and the formula for the Malliavin derivative yield
\begin{align*}
 \|F_n\|_{L^2(\mP\times\rho;E)}^2
 =\sum_{\ell\geq0}\ell!\|f_{n,\ell}\|_\ell^2,\quad 
 \bE\int_0^T\int|\mathrm D_rF_n|^2\rho\dd r
 =\sum_{\ell\geq1}\ell\,\ell!\|f_{n,\ell}\|_\ell^2.
\end{align*}
Consequently, the contribution of the higher chaos orders satisfies
\[
 \sup_n \Big\| \sum_{\ell>N}^\infty I_\ell(f_{n,\ell}) \Big\|_{L^2(\mP\times\rho)}^2= \sup_n\sum_{\ell>N}\ell!\|f_{n,\ell}\|_\ell^2
 \leq\frac{C}{N+1}\to 0, \quad N\to\infty.
\]

It therefore suffices to prove relative compactness of $(f_{n,\ell})_n$ in $\|\cdot\|_\ell$ for each fixed $\ell\geq0$. We first control the spatial variables. Since chaos projections commute with weak derivatives, we obtain from the isometry that $\sup_n\|\nabla_x f_{n,\ell}\|_\ell^2\leq C/\ell!$. For $\tau_z^xf(\cdot,x):=f(\cdot,x+z)$, the comparison $\rho(x-z)\leq C\rho(x)$ for $|z|\leq1$ and the assumed $\rho_*$ bound yield
\begin{align*}
 \sup_n\|\tau_z^xf_{n,\ell}-f_{n,\ell}\|_\ell^2\leq\frac{C}{\ell!}|z|^2, 
 \quad
 \sup_n\|\1_{B_R^c}f_{n,\ell}\|_\ell^2
 \leq\frac{C}{\ell!}\sup_{|x|>R}\frac{\rho(x)}{\rho_*(x)}=\frac{C}{\ell!}(1+R^2)^{-\frac{q}{2}}.
\end{align*}

For time translations, suppose $\ell\geq1$. Let $\tau_h^{(j)}$ translate $t_j$ by $h$, and let $\|\cdot\|_{\ell,h,j}$ be the kernel norm restricted to $t_j\in[0,T-h]$. Using the identity $\mathrm D_r^kI_\ell(f)=\ell I_{\ell-1}(f((r,k),\cdot))$, the chaos isometry and symmetry in the pairs $(t_j,k_j)$, we obtain, for every $1\leq j\leq\ell$,
\[
 \sup_n\|\tau_h^{(j)}f_{n,\ell}-f_{n,\ell}\|_{\ell,h,j}^2
 \leq (\ell\,\ell!)^{-1}\sup_n\int_0^{T-h}\bE\int
 |\mathrm D_{r+h}F_n-\mathrm D_rF_n|^2\rho\dd r
 \to0\quad(h\downarrow0).
\]

Since $\rho$ is comparable to $1$ on each $B_R$, we can apply the Fr\'echet--Kolmogorov criterion to these translation estimates to obtain relative compactness on $[0,T]^\ell\times B_R$. Combining this with the spatial tail bound, we conclude that $(f_{n,\ell})_n$ is relatively compact in $\|\cdot\|_\ell$, as required.
\end{proof}

\begin{proposition}\label{prop:limit}
There exists $\overline\delta=\overline\delta(d,\lambda)>0$ such that, for every $0\leq\delta<\overline\delta$, every admissible sequence (cf. Definition \ref{def:approx}) has a subsequence for which $X^n,Y^n$ converge in $\mP$-probability, locally uniformly in $(s,t,x)$, on the original Wiener space. On an event of $\mP$-measure one, any such pair of limits satisfies, for all $0\leq s\leq u\leq t\leq T$,
\[
 X_{u,t}\circ X_{s,u}=X_{s,t},\quad
 Y_{s,u}\circ Y_{u,t}=Y_{s,t},\quad
 X_{s,t}\circ Y_{s,t}=Y_{s,t}\circ X_{s,t}=\Id,
\]
and has the regularity and Sobolev bounds stated in Theorem~\ref{thm:flow}(ii).
\end{proposition}
\begin{proof}
By \eqref{eq:disp-u}, Proposition \ref{prop:jac} and Proposition~\ref{prop:mall}, Lemma~\ref{lem:compact} applies to $X^n_{s,t}-\Id$ and $Y^n_{s,t}-\Id$ for each fixed $(s,t)$. By a diagonal argument, we select a common subsequence converging in $L^2(\mP\times\rho)$ for every $(s,t)$ in a countable dense set. By the Arzel\`a--Ascoli theorem, on sets where the local $C^{\frac{1}{8}}$ norms are uniformly bounded, we obtain a unique locally uniform limit from this convergence. The uniform moment bounds for these norms in Proposition~\ref{prop:holder} therefore imply that the subsequence is Cauchy in probability in the local uniform topology. Passing to a further subsequence, we obtain
\[
 (X^n,Y^n)\to(X,Y)
 \quad\mP\mbox{-a.s., locally uniformly in }(s,t,x).
\]
The flow identities follow by passing to the limit in
\[
 X^n_{u,t}\circ X^n_{s,u}=X^n_{s,t},\quad
 Y^n_{s,u}\circ Y^n_{u,t}=Y^n_{s,t},\quad
 X^n_{s,t}\circ Y^n_{s,t}=Y^n_{s,t}\circ X^n_{s,t}=\Id.
\]
Indeed, the inner maps send each compact set into a common random compact set. Thus, the identities hold simultaneously on the event of local uniform convergence.

Fatou's lemma and Proposition~\ref{prop:holder} imply the joint $1/8$-H\"older regularity. For each $(s,x)$ and every $p\geq1$, we also obtain from \eqref{eq:disp}
\[
 \bE\sup_{u\leq r\leq v}|X_{s,r}(x)-X_{s,u}(x)|^p
 \leq C_p(v-u)^{\frac{p}{2}},\quad s\leq u\leq v\leq T.
\]
Here $C_p$ depends only on $d,\lambda,p$ and $\|g\|_{L^1(0,T)}$. Kolmogorov's criterion yields H\"older continuity for every exponent $\alpha<1/2$. The same argument applies to $s\mapsto Y_{s,t}(x)$ for fixed $(t,x)$. Finally, \eqref{eq:disp-u} and \eqref{eq:jac} yield uniform bounds in $L^P(\Omega,\mP;W^{1,P}(B_R))$ at each fixed $(s,t)$. Using weak compactness and the aforementioned convergence, we obtain the stated Sobolev bounds for $X_{s,t}$ and $Y_{s,t}$.
\end{proof}

\subsection{The limiting equation and uniqueness}\label{sec:weak}

We first show that each subsequential limit $X$ in Proposition~\ref{prop:limit} is a strong solution of \eqref{eq:ito} and satisfies the Krylov-type estimate \eqref{eq:occ}. We then establish joint uniqueness in law within $\sK_2(s,x)$. These results yield pathwise uniqueness in this class and convergence of the full sequence to a flow independent of the admissible approximation, completing the proof of Theorem~\ref{thm:flow}.

\begin{proposition}\label{prop:ident}
There exists $\overline\delta=\overline\delta(d,\lambda)>0$ such that, for every $0\leq\delta<\overline\delta$, the following holds. Let $X$ be a subsequential limit of $X^n$ in Proposition~\ref{prop:limit}. For every $(s,x)$, $X_{s,\cdot}(x)$ is a strong solution of \eqref{eq:ito} and satisfies \eqref{eq:occ}.
\end{proposition}
\begin{proof}
Fix $(s,x)$. Since $X^n_{s,\cdot}(x)\to X_{s,\cdot}(x)$ uniformly in probability, $X_{s,\cdot}(x)$ is continuous and adapted.

For $f\in\bF_\beta(g_f)$, put $f_N:=|f|\wedge N$ and $f_{N,\eps}:=f_N*_x\eta_\eps$. Truncation and convolution preserve the form bound; by Jensen's inequality, we also have
\begin{equation}\label{eq:src-moll}
 \mathfrak{n}_{f_{N,\eps}}(u,v;R)\leq\mathfrak{n}_f(u,v;R),\qquad R>0.
\end{equation}
For $s\leq u<v\leq T$, since $f_{N,\eps}$ is bounded and uniformly Lipschitz in space, we may pass to the limit along the paths in Lemma~\ref{lem:occ}(i) and use \eqref{eq:src-moll} to get
\[
 \bE\left(\int_u^v f_{N,\eps}(r,X_{s,r}(x))\dd r\right)^m
 \leq m!\,C^m(v-u)^{\frac{m}{2}(1-\frac{d}{q})}\mathfrak{n}_f(u,v)^{\frac{2m}{q}},
\]
where $C$ depends only on $d,\lambda,\beta$, $\|g\|_{L^1(0,T)}$ and $\|g_f\|_{L^1(0,T)}$. Along a sequence $\eps\downarrow0$, $f_{N,\eps}\to f_N$ Lebesgue-a.e. By Proposition~\ref{prop:semi}(ii), $X_{s,r}(x)$ has a density for every $r>s$. Thus, $f_{N,\eps}(r,X_{s,r}(x))\to f_N(r,X_{s,r}(x))$ for $(\dd r\times\mP)$-almost every $(r,\omega)$. Dominated convergence as $\eps\downarrow0$, followed by monotone convergence as $N\to\infty$, implies \eqref{eq:occ}.

Now choose $\chi_R\in C_c^\infty(\mR^d)$ with $0\leq\chi_R\leq1$, $\chi_R=1$ on $B_R$ and $\chi_R=0$ outside $B_{2R}$, and set
\[
 f_n:=\chi_R\widehat b^n,\quad f:=\chi_R\widehat b.
\]
By \eqref{eq:ito-lim}, $\mathfrak{n}_{f_n-f}(s,T;1)\to0$. Truncating $f$ and then mollifying in space, choose bounded fields $h_k$, uniformly Lipschitz in space for each $k$, such that
\[
 \mathfrak{n}_{h_k-f}(s,T;1)\to0.
\]
The fields $f_n,f,h_k$, and hence their differences, satisfy a common form-bound inequality, uniformly in $n,k$. Applying Lemma~\ref{lem:occ}(i) and \eqref{eq:occ} with $m=1$, together with the covering estimate in the proof of Lemma~\ref{lem:occ}, we obtain
\begin{align*}
 &\bE\int_s^T|f_n(r,X^n_{s,r}(x))-f(r,X_{s,r}(x))|\dd r\\
 \leq{} & C\mathfrak{n}_{f_n-f}(s,T;1)^{\frac{2}{q}}
       +C\mathfrak{n}_{h_k-f}(s,T;1)^{\frac{2}{q}}+\bE\int_s^T|h_k(r,X^n_{s,r}(x))-h_k(r,X_{s,r}(x))|\dd r.
\end{align*}
Here $C$ depends only on $d,\lambda,T$ and $\|g\|_{L^1(0,T)}$. For fixed $k$, the last term tends to zero by convergence of the paths and boundedness of $h_k$. Letting first $n\to\infty$ and then $k\to\infty$, we obtain
\begin{equation*}
 \bE\sup_{s\leq t\leq T}\left|
 \int_s^t\bigl[f_n(r,X^n_{s,r}(x))-f(r,X_{s,r}(x))\bigr]\dd r\right|
 \to0.
\end{equation*}

Put $G_n(r):=(\chi_R\, \sigma^n)(r,X^n_{s,r}(x))$ and $G(r):=(\chi_R\, \sigma)(r,X_{s,r}(x))$. Since $\sigma^n,\sigma$ are bounded, the functions $|\chi_R(\sigma^n-\sigma)|^2$ belong to $\bF_0(C)$, where $C$ depends only on $d,\lambda$, and
\[
 \mathfrak{n}_{|\chi_R(\sigma^n-\sigma)|^2}(s,T;1)^2
 \leq C\|\sigma^n-\sigma\|_{L^2((s,T)\times B_{2R})}^2
 \to0.
\]
Applying the same approximation argument with Lemma~\ref{lem:occ}(i) and \eqref{eq:occ}, again with $m=1$, we get
\[
 \bE\int_s^T|G_n(r)-G(r)|^2\dd r\to0.
\]
The Burkholder--Davis--Gundy inequality then yields
\[
 \bE\sup_{s\leq t\leq T}\left|
 \int_s^t[G_n(r)-G(r)]\dd W_r\right|^2\to0.
\]
By \eqref{eq:disp} and Fatou's lemma,
\[
 \lim_{R\to\infty}\sup_n\mP\left(
 \sup_{s\leq r\leq T}(|X^n_{s,r}(x)|\vee|X_{s,r}(x)|)>R\right)=0.
\]
Letting $R\to\infty$ removes the cutoff. By \eqref{eq:occ}, we also have $\int_s^T|\widehat b(r,X_{s,r}(x))|\dd r<\infty$ almost surely. Passing to the limit in the equation for $X^n_{s,\cdot}(x)$ proves \eqref{eq:ito} as an identity of continuous processes.
\end{proof}

We next prove joint uniqueness in law within $\sK_2(s,x)$.

\begin{proposition}\label{prop:weak}
There exists $\overline\delta=\overline\delta(d,\lambda)>0$ such that, for every $0\leq\delta<\overline\delta$, the following assertions hold. Let $(Z,B)\in\sK_2(s,x)$ be defined on $(\Omega',\sF',(\sF'_r),\mP')$. Then, for every $s\leq r\leq t\leq T$ and $\varphi\in C_c^\infty(\mR^d)$,
\begin{equation}\label{eq:trans}
 \bE_{\mP'}[\varphi(Z_t)\mid\sF'_r]
 =P_{r,t}\varphi(Z_r),\quad\mP'\mbox{-a.s.}
\end{equation}
Moreover, the joint law $\Law_{\mP'}(Z,B^s)$ is the same for all $(Z,B)\in\sK_2(s,x)$, where $B^s_t:=B_t-B_s$.
\end{proposition}
\begin{proof}
By \eqref{eq:k2}, $Z$ spends zero time in every fixed space-time Lebesgue-null set. In particular, \eqref{eq:ellip} holds along $Z$ for $(\dd r\times\mP')$-almost every $(r,\omega)$. Fix $t$ and $\varphi$, and put
\[
 u(r):=P_{r,t}\varphi,\quad u_\eps:=u*_x\eta_\eps,
 \quad R_\eps:=\partial_ru_\eps
                  +a:\nabla^2u_\eps+\widehat b\cdot\nabla u_\eps.
\]
By Proposition~\ref{prop:semi}(iii) and the form bound of $\widehat b$,
\begin{align*}
\int_0^t\|\widehat b\cdot(\nabla u_\eps-\nabla u)\|_2^2\dd r \leq C\delta\int_0^t\|\nabla^2u_\eps-\nabla^2u\|_2^2\dd r +C\int_0^t g(r)\|\nabla u_\eps-\nabla u\|_2^2\dd r \to 0.
\end{align*}
The last term tends to zero, since $u\in L^\infty(0,t;H^1)$ and $g\in L^1(0,t)$. Combining this with boundedness of $a$ and $\partial_ru_\eps\to\partial_ru$ in $L^2$, we obtain
\[
 R_\eps\to0
 \quad\mbox{in }L^2((0,t)\times\mR^d).
\]

Let $s<r<t$ and $H$ be a bounded $\sF'_r$-measurable random variable. For fixed $\eps$, $u_\eps$ has bounded spatial derivatives and $\partial_ru_\eps\in L^1(0,t;C_b)$, so It\^o's formula applies by time regularization. Taking expectations and using \eqref{eq:k2}, we get
\begin{align*}
 & \left|\bE_{\mP'}\bigl[H\{u_\eps(t,Z_t)-u_\eps(r,Z_r)\}\bigr]\right|\\
 ={} & \left|\bE_{\mP'}\left[H\int_r^tR_\eps(\tau, Z_\tau)\dd \tau\right] \right|\leq\|H\|_\infty C_{r,Z}
       \|R_\eps\|_{L^2((r,t)\times\mR^d)}
 \to0.
\end{align*}
Here $C_{r,Z}$ is the constant in \eqref{eq:k2}, depending on $r$ and the chosen solution $Z$ on $[s,T]$. Since $u$ is bounded and continuous and $u(t)=\varphi$, letting $\eps\downarrow0$ proves \eqref{eq:trans} for $r>s$. For $r=s$, take $H\in L^\infty(\sF'_s,\mP')$ in the same identity and let $r\downarrow s$. By continuity, we have $u(r,Z_r)\to u(s,x)$; we then conclude by bounded convergence.

A standard monotone-class argument and repeated conditioning in \eqref{eq:trans} show that all solutions in $\sK_2(s,x)$ have the same path law $\mathsf Q$.

On canonical path space with law $\mathsf Q$, let $Z$ be the canonical process and put
\[
 N_t:=Z_t-x-\int_s^t\widehat b(r,Z_r)\dd r.
\]
Then $N$ is a square-integrable martingale in the natural filtration, with
\[
 \dd\langle N\rangle_r=\sigma\sigma^\top(r,Z_r)\dd r\leq CI\dd r.
\]
Set $\sigma^{-1}=0$ where \eqref{eq:ellip} fails, so that $\sigma^{-1}$ is bounded. Approximating by elementary predictable integrands and applying Doob's inequality, we construct a Borel path functional $\Psi$ such that
\[
 \Psi(Z)_t=\int_s^t\sigma(r,Z_r)^{-1}\dd N_r,
 \quad\mathsf Q\mbox{-a.s.}
\]
More precisely, take an almost surely uniformly convergent subsequence of the elementary integrals and set its limit to zero off the convergence set. For any $(Z,B)\in\sK_2(s,x)$, we obtain from the same approximations on its probability space that
\begin{equation*}
 \begin{aligned}
 \Psi(Z)_t
 =\int_s^t\sigma(r,Z_r)^{-1}\sigma(r,Z_r)\dd B_r
 =B_t-B_s,\quad\mP'\mbox{-a.s.}
 \end{aligned}
\end{equation*}
Thus, $\Law_{\mP'}(Z,B^s)=(\Id,\Psi)_\#\mathsf Q$, proving joint uniqueness in law.
\end{proof}

\begin{proof}[Proof of Theorem~\ref{thm:flow}]
Choose $\overline\delta$ to be the minimum of the thresholds in the preceding results and in Lemma~\ref{lem:pde} with $\kappa=2q$, also allowing ellipticity parameter $\lambda/2$ for the approximations. Since $q=d+1$, $P=50q$ and $M=10^3q^2$ are fixed by $d$, this minimum depends only on $d,\lambda$.

Choose an admissible approximation from Proposition~\ref{prop:approx}, and let $(X,Y)$ be a subsequential limit given by Proposition~\ref{prop:limit}. Proposition~\ref{prop:ident} shows that $X_{s,\cdot}(x)$ is a strong solution of \eqref{eq:ito} for every $(s,x)$. The composition and inverse identities in \eqref{eq:cocycle} hold on the common event of local uniform convergence; the diagonal identities follow from $X^n_{s,s}=Y^n_{s,s}=\Id$. Since $X_{s,t}$ and $Y_{s,t}$ are continuous and mutually inverse, they are homeomorphisms. This proves (i). The regularity assertions in (ii) follow from Proposition~\ref{prop:limit}, and the estimate in (iii) follows from Proposition~\ref{prop:ident}.

For (iv), fix $(s,x)$ with $s<T$. By Proposition~\ref{prop:semi}(ii), we have $X_{s,\cdot}(x)\in\sK_2(s,x)$. Since the approximations $X^n_{s,\cdot}(x)$ depend only on the Brownian increments $W^s_t:=W_t-W_s$, their limit admits a measurable solution map
\[
 X_{s,\cdot}(x)=F_{s,x}(W^s)\quad\mP\mbox{-a.s.}.
\]

Let $(Z,B)\in\sK_2(s,x)$, with probability measure $\mP'$. By the joint uniqueness in law in Proposition~\ref{prop:weak},
\[
 \Law_{\mP'}(Z,B^s)
 =\Law_{\mP}(X_{s,\cdot}(x),W^s)
 =\Law_{\mP}(F_{s,x}(W^s),W^s).
\]
The last measure is supported on the graph of $F_{s,x}$. Hence
\[
 Z=F_{s,x}(B^s)\quad\mP'\mbox{-a.s. in }C([s,T];\mR^d).
\]
This representation establishes pathwise uniqueness in $\sK_2(s,x)$ and completes the proof of (iv).

The argument also shows that the construction is independent of the admissible approximation. Indeed, any two subsequential limits of $X^n$ agree at each $(s,x)$ by pathwise uniqueness. Using a countable dense set and joint continuity, we conclude that the random fields are equal on one event of probability one, and their inverses agree as well. We then deduce from Proposition~\ref{prop:limit} that the whole sequence converges in probability in the local uniform topology. Interleaving any two admissible approximations shows that they converge to the same pair $(X,Y)$.
\end{proof}

\part{Ill-posedness}\label{part:ill}

Part~\ref{part:well} establishes strong solvability under small form bounds. We now ask whether strong solvability can hold throughout the critical endpoint class $L^\infty_tL^d_x$ without a size or time-continuity condition. Theorem~\ref{thm:counter} answers this question negatively, already for additive noise and divergence-free drifts: this class contains equations with no strong solution. To prove it, we use logarithmic rescaling and the criterion for strong solvability in Proposition~\ref{prop:strong}. For the profile, the Lyapunov calculation and the random SRB theorem show that conditioning on the entire past noise leaves a non-Dirac distribution of the state, which is incompatible with strong solvability.

\section{Self-similar equations and strong solvability}\label{sec:entry}

Throughout this section, $d\geq1$ and $B\in C_c^\infty(\mR^d;\mR^d)$ is arbitrary. Incompressibility is not assumed. Set
\begin{equation*}
 A(z):=-\tfrac12z+B(z),\quad
 b(t,x):=t^{-\frac{1}{2}}B(x/\sqrt t),\quad 0<t\leq1,
\end{equation*}
and consider
\begin{equation}\label{eq:scale}
 X_t=\int_0^t b(r,X_r)\dd r+W_t.
\end{equation}
The value of $b$ at $t=0$ does not affect the equation; we set it to zero. For a weak solution in this section, we write $(\Omega,\sF,(\sF_t)_{0\leq t\leq1},\mP)$ for its filtered probability space. We first prove weak existence and uniqueness in law, and then characterize when the solution is adapted to the natural filtration of its driving Brownian motion.

\subsection{The equation in logarithmic time}

For any weak solution of \eqref{eq:scale}, put
\begin{equation}\label{eq:def-Z}
 Z_s:=e^{-\frac{s}{2}}X_{e^s},\quad
 \beta((r,s]):=\int_{e^r}^{e^s}u^{-\frac{1}{2}}\dd W_u,\quad r<s\leq0.
\end{equation}
For every $r<0$, the process $(\beta((r,s]))_{r\leq s\leq0}$ is a standard Brownian motion started at zero at time $r$, with respect to $(\sF_{e^s})_{r\leq s\leq0}$. Then 
\begin{equation}\label{eq:log}
 Z_s-Z_r=\int_r^s A(Z_a)\dd a+\beta((r,s]).
\end{equation}
The inverse change of variables is
\begin{equation}\label{eq:inv-lamp}
 X_t=\sqrt t\,Z_{\log t},\quad
 W_t=\int_{(-\infty,\log t]}e^{\frac{a}{2}}\dd\beta_a.
\end{equation}
The last integral exists in $L^2(\mP)$ because $\int_{-\infty}^{\log t}e^a\dd a=t$.

\medskip

Let $P_t$ be the semigroup associated with \eqref{eq:log}. Since $A$ is globally Lipschitz, this equation has a pathwise unique strong solution on every bounded interval of logarithmic time.

Put
\begin{equation*}
 \gamma_d(\dd z):=(2\pi)^{-\frac{d}{2}}e^{-\frac{|z|^2}{2}}\dd z.
\end{equation*}

\begin{lemma}\label{lem:harris}
The semigroup $P_t$ has a unique invariant probability $\mu_B$, with smooth strictly positive density and finite second moment. There are $C,c>0$, depending only on $d$ and $B$, such that
\begin{equation}\label{eq:harris}
 \|\nu P_t-\mu_B\|_{\mathrm{TV}}
 \leq Ce^{-ct}\int(1+|z|^2)\nu(\dd z)
\end{equation}
for every probability measure $\nu$ with finite second moment. If $\div B(z)=B(z)\cdot z$, then $\mu_B=\gamma_d$.
\end{lemma}
\begin{proof}
For $V(z):=1+|z|^2$, the generator $L:=\frac12\Delta+A\cdot\nabla$ satisfies
\begin{equation*}
 LV=d-|z|^2+2B(z)\cdot z\leq C_B-\tfrac12V.
\end{equation*}
Here $C_B$ depends only on $d$ and $\|B\|_\infty$. Uniform ellipticity and smooth coefficients ensure the existence of a continuous, strictly positive transition density. Hence every compact set satisfies a minorization condition at each fixed positive time. Applying the continuous-time Foster--Lyapunov criteria of Meyn--Tweedie \cite[Theorems~4.2 and~6.1]{MT93stability} to the bound above, we obtain positive Harris recurrence, a unique invariant probability $\mu_B$ with $\int V\dd\mu_B<\infty$, and \eqref{eq:harris} for all $t\geq0$.

It remains to describe the invariant law. Since $\mu_B=\mu_BP_1$, it has a density $\rho_B$. The stationary equation $L^*\rho_B=0$ and elliptic regularity make $\rho_B$ smooth, while positivity of the transition density makes it strictly positive.

Finally, suppose $\div B(z)=B(z)\cdot z$. Then
\[
 \div(B\gamma_d)=\gamma_d\bigl(\div B-B\cdot z\bigr)=0.
\]
Thus, the additional drift $B$ preserves $\gamma_d$, as does the Ornstein--Uhlenbeck part of the generator. Therefore $\gamma_d$ is invariant for $P_t$, and uniqueness implies $\mu_B=\gamma_d$.
\end{proof}

\begin{proposition}\label{prop:entry}
Equation~\eqref{eq:scale} admits a weak solution, and its law on $C([0,1];\mR^d)$ is unique. For every weak solution, the process $Z$ given by \eqref{eq:def-Z} is stationary with invariant probability $\mu_B$.
\end{proposition}
\begin{proof}
For an arbitrary weak solution, boundedness of $B$ implies that
\begin{equation}\label{eq:entry-m}
 |X_t-W_t|\leq2\|B\|_\infty\sqrt t,\quad
 \sup_{s\leq0}\bE|Z_s|^2\leq2d+8\|B\|_\infty^2.
\end{equation}
The pathwise bound holds $\mP$-a.s. for all $t\in[0,1]$. Let $\nu_s:=\Law(Z_s)$. By pathwise uniqueness for the Lipschitz SDE \eqref{eq:log} on each bounded time interval, we have $\nu_s=\nu_rP_{s-r}$ for $r<s$, as well as the Markov property relative to the given filtration. Letting $r\to-\infty$, and using \eqref{eq:harris} and \eqref{eq:entry-m}, we obtain that $\nu_s=\mu_B$. Thus, all finite-dimensional distributions are
\[
 \mu_B(\dd z_1)P_{s_2-s_1}(z_1,\dd z_2)\cdots
                         P_{s_k-s_{k-1}}(z_{k-1},\dd z_k).
\]
Thus, there is at most one law for $X$ on continuous path space.

For existence, construct the stationary diffusion $Z$ on the whole real line on a probability space with measure $\mP$ from these consistent laws on bounded time intervals. Let $(\mathcal G_s)_{s\in\mR}$ be the usual augmentation of $(\sigma(Z_a:a\leq s))_{s\in\mR}$. Then
\[
 \beta((r,s]):=Z_s-Z_r-\int_r^s A(Z_a)\dd a
\]
are Brownian increments independent of $\mathcal G_r$. Define $X,W$ by \eqref{eq:inv-lamp}, with $W_0=0$. By Doob's inequality, we can choose a continuous square-integrable martingale version of $W$ with respect to $(\mathcal G_{\log t})_{t>0}$, extended at $t=0$ by $\bigcap_{s\leq0}\mathcal G_s$. Its quadratic variation is
\[
 \langle W^i,W^j\rangle_t
 =\delta_{ij}\int_{-\infty}^{\log t}e^a\dd a
 =\delta_{ij}t,\quad t>0.
\]
Thus, L\'evy's characterization shows that $W$ is a standard $d$-dimensional Brownian motion relative to this filtration.

For $0<u<t$, applying It\^o's formula, we get
\[
 X_t-X_u=W_t-W_u+\int_u^t b(r,X_r)\dd r.
\]
The drift is absolutely integrable at zero, since $\|b(r)\|_\infty\leq\|B\|_\infty r^{-\frac{1}{2}}$. Therefore $X_t$ has a limit $\mP$-a.s. as $t\downarrow0$. By stationarity, we have $\bE|X_t|^2=t\int|z|^2\mu_B(\dd z)\to0$, so that limit is zero $\mP$-a.s.. This implies that $X$ is continuous on $[0,1]$ and solves \eqref{eq:scale}.
\end{proof}

The restriction $s\leq0$ in \eqref{eq:def-Z} corresponds to the original time interval $[0,1]$. Any such $Z$ can be extended to a stationary process on $\mR$: on an enlarged probability space, add Brownian increments for $s>0$, independent of $\sF_1$, and solve \eqref{eq:log} forward from $Z_0$. Since $\Law(Z_0)=\mu_B$, the extended process is stationary. This extension is consistent with the forward dynamics and two-sided noise used below. The conditional laws below use only the past $s\leq0$, so the original equation may remain on $[0,1]$.

\subsection{Sample measures and strong solvability}\label{sec:sample}

For every $0<\tau<1$, the drift $b$ is globally Lipschitz in space on $[\tau,1]$. Hence $X_1$ is a measurable function of $X_\tau$ and the Brownian increment path $(W_t-W_\tau)_{\tau\leq t\leq1}$. The question is whether, as $\tau\downarrow0$, $X_1$ still depends on information from $X_\tau$ beyond that contained in the Brownian path. Although $X_\tau\to0$ almost surely, the solution map from time $\tau$ to time $1$ also depends on $\tau$, so this convergence alone does not show that $X_1$ is a measurable function of $W$. Precisely, we ask whether
\[
 \Law\bigl(X_1\mid\sigma(W_t:0\leq t\leq1)\bigr)
       \quad\mbox{is a Dirac mass almost surely.}
\]
This is equivalent to $X_1$ being measurable with respect to the completed Brownian sigma field. For our self-similar equation, Proposition~\ref{prop:strong} will show that this single-time condition is equivalent to strong solvability on $[0,1]$.

In logarithmic variables, the earlier state is $Z_{\log\tau}=\tau^{-\frac{1}{2}}X_\tau$, whose law is $\mu_B$ for every $0<\tau\leq1$ by Proposition~\ref{prop:entry}. We therefore study how this fixed distribution evolves when the noise between $\log\tau$ and $0$ is given, and then let $\log\tau\to-\infty$.

This is where random dynamical systems enter the argument. Once the noise is fixed, the equation for $Z$ defines smooth maps from earlier states to the present state. Evolving $\mu_B$ by these maps retains the dependence on the observed noise, whereas the Markov semigroup averages over it. The limiting conditional distributions are the sample measures constructed below; the invariant law $\mu_B$ alone does not tell us whether they are Dirac masses.

Lyapunov exponents describe the long-time growth of infinitesimal perturbations of the state under the same noise. For the profiles constructed below, we combine a positive top exponent with the random SRB theorem of Blumenthal--Young, whose hypotheses we check in Section~\ref{sec:srb}. The theorem ensures absolute continuity of conditional measures along unstable manifolds of positive dimension, which rules out Dirac sample measures and hence strong solvability.

We represent the dynamics by compositions of independent random diffeomorphisms, following \cite[Section~1]{BY19equivalence}. Each diffeomorphism is determined by a Brownian path on a unit interval:
\[
 \Omega_0:=C_0([0,1];\mR^d),\quad Q:=\mbox{Wiener measure}.
\]
For $\eta\in\Omega_0$, consider the integral equation
\[
 Z_t^{z,\eta}=z+\int_0^t A(Z_r^{z,\eta})\dd r+\eta(t),
       \quad 0\leq t\leq1.
\]
Define the map $f_\eta:\mR^d\to\mR^d$ by
\[
f_\eta(z):=Z_1^{z,\eta}.
\]
Since $A$ is a globally Lipschitz, smooth vector field, $f_\eta$ is a smooth diffeomorphism, and the map $\eta\mapsto f_\eta$ is continuous.

Write
\[
 \Omega_-:=\Omega_0^{\mZ_{\leq0}},\quad Q_-:=Q^{\mZ_{\leq0}},
 \quad
 \Omega_+:=\Omega_0^{\mN},\quad Q_+:=Q^{\mN}.
\]
The coordinate $\omega_j$ is the Brownian increment path over logarithmic time $[j-1,j]$. A past history is $\underline\omega^-:=(\omega_j)_{j\leq0}\in\Omega_-$. For $N\geq1$, the evolution from time $-N$ to time $0$ uses the noise paths $\omega_{-N+1},\ldots,\omega_0$. The corresponding map is
\begin{equation}\label{eq:pullback}
 F_N(\underline\omega^-)
    :=f_{\omega_0}\circ f_{\omega_{-1}}\circ\cdots
          \circ f_{\omega_{-N+1}}.
\end{equation}
If the state at time $-N$ has law $\mu_B$ and is independent of these noise inputs, then its conditional law at time $0$, given these inputs, is
\[
 \pi_N(\underline\omega^-)
    :=(F_N(\underline\omega^-))_\#\mu_B.
\]
As $N$ increases, we condition on a longer segment of the past noise. The limiting distribution is called the \emph{sample measure}; it describes the uncertainty about the present state that remains after the entire past noise is known. We denote it by $\mu_{\underline\omega^-}$, following the notation for sample measures in \cite[Lemma~2]{BY19equivalence}. The law of this random probability measure is the \emph{statistical equilibrium} in Le Jan's work \cite{L87equilibre}.

\begin{proposition}\label{prop:sample}
There is a Borel map $\underline\omega^-\mapsto\mu_{\underline\omega^-}$ from $\Omega_-$ to $\cP(\mR^d)$, equipped with the topology of weak convergence, such that
\begin{equation}\label{eq:sample}
 \pi_N(\underline\omega^-)
       \to \mu_{\underline\omega^-},
       \quad Q_-\mbox{-a.s.}
\end{equation}
\end{proposition}
\begin{proof}
Let $\sH_N:=\sigma(\omega_{-N+1},\ldots,\omega_0)$ and write $\pi_N(h):=\int h\dd\pi_N$. Stationarity of $\mu_B$ says
\begin{equation}\label{eq:one-step}
 \int_{\Omega_0}(f_\eta)_\#\mu_B\,Q(\dd\eta)=\mu_B.
\end{equation}
Since $F_{N+1}=F_N\circ f_{\omega_{-N}}$, we have $\pi_{N+1}(h)= \int h(F_N)\,\dd \big((f_{\omega_{-N}})_\# \mu_B \big)$ and $\omega_{-N}$ is independent of $\sH_N$; \eqref{eq:one-step} implies 
\[
 \bE_{Q_-}[\pi_{N+1}(h)\mid\sH_N]=\pi_N(h).
\]
For $h\in C_b(\mR^d)$ this is a bounded martingale. 
Martingale convergence for a countable family in $C_b(\mR^d)$ that separates probability measures makes all subsequential limits coincide, proving \eqref{eq:sample}. The limit $\mu_{\underline{\omega}^-}$ is probability measure-valued and Borel measurable.
\end{proof}

We add independent future inputs and take $\underline\omega:=(\omega_j)_{j\in\mZ}$ with law $Q^{\mZ}$. The shift $\theta$ is given by $(\theta\underline\omega)_j:=\omega_{j+1}$, the new past is $(\theta\underline\omega)^-=(\omega_{j+1})_{j\leq0}$. Shifting the finite compositions in \eqref{eq:pullback}, one sees that $F_N(\theta \underline{\omega}^-) = f_{\omega_1}\circ F_{N-1}(\underline{\omega}^-) $ and thus $\pi_N (\theta \underline{\omega}^-)= (f_{\omega_1})_\#\pi_{N-1} (\underline{\omega}^-) $; taking limits yields
\begin{equation*}
 \begin{aligned}
 & (f_{\omega_1})_\#\mu_{\underline\omega^-}=\mu_{(\theta\underline\omega)^-}, \quad Q^{\mZ}\mbox{-a.s.},\\
 & \int_{\Omega_-}\mu_{\underline\omega^-}\,Q_-(\dd\underline\omega^-)=\mu_B,
 \end{aligned}
\end{equation*}
where the second identity follows from \eqref{eq:one-step} and bounded convergence. Thus, the joint law
\[
 \mu_B^*(\dd z,\dd\underline\omega)
     :=\mu_{\underline\omega^-}(\dd z)Q^{\mZ}(\dd\underline\omega)
\]
is invariant under the skew product $(z,\underline\omega)\mapsto(f_{\omega_1}(z),\theta\underline\omega)$. The current state has law $\mu_B$ and is independent of the future inputs, because $\mu_{\underline\omega^-}$ depends only on the past. This is the invariant measure for the two-sided skew product in \cite[Lemmas~1 and~2]{BY19equivalence}; its conditional distributions given the noise are precisely the sample measures just constructed.

We now show that the sample measure is a conditional law for the original equation. The only change is the time variable: one unit of logarithmic time corresponds to an interval $[e^{j-1},e^j]$ in the original Brownian path. For $j\leq0$ and $v\in[0,1]$, put $t_j(v):=e^{j-1+v}$ and define, for $w\in\Omega_0$,
\begin{equation*}
 \Xi_j(w)(v)
   :=t_j(v)^{-\frac{1}{2}}w(t_j(v))
       -t_j(0)^{-\frac{1}{2}}w(t_j(0))
       +\frac12\int_{t_j(0)}^{t_j(v)}r^{-\frac{3}{2}}w(r)\dd r.
\end{equation*}
Each $\Xi_j:\Omega_0\to\Omega_0$ is continuous. If $W$ is Brownian, integration by parts yields
\[
 \Xi_j(W)(v)=\int_{e^{j-1}}^{e^{j-1+v}}r^{-\frac{1}{2}}\dd W_r.
\]
The processes $\Xi_j(W)$, $j\leq0$, are independent standard Brownian motions on $[0,1]$: the integration intervals are disjoint and $\int_{t_j(u)}^{t_j(v)}r^{-1}\dd r=v-u$. Therefore,
\begin{equation}\label{eq:blocks}
 \Xi:\Omega_0\to\Omega_-,\quad
 \Xi(w):=(\Xi_j(w))_{j\leq0},\quad \Xi_\#Q=Q_-.
\end{equation}
The map $\Xi$ converts $W$ into the sequence of past Brownian increment paths used above.

\begin{proposition}[Criterion for strong solvability]\label{prop:strong}
For every weak solution $(X,W)$ of \eqref{eq:scale}, with probability measure $\mP$ and natural filtration $(\sF_t^W)$ of $W$,
\begin{equation}\label{eq:cond}
 \Law(X_1\mid\sF_1^W)=\mu_{\Xi(W)}
       \quad\mP\mbox{-a.s.}
\end{equation}
Equation \eqref{eq:scale} has a strong solution on $[0,1]$ if and only if $\mu_{\underline\omega^-}$ is a Dirac mass for $Q_-$-a.s. $\underline\omega^-$. In that case pathwise uniqueness also holds.
\end{proposition}
\begin{proof}
\emph{The conditional law.} Let $Z$ be given by \eqref{eq:def-Z} and put
\[
 \cG_N:=\sigma(\Xi_{-N+1}(W),\ldots,\Xi_0(W)).
\]
By Proposition~\ref{prop:entry}, $Z_{-N}$ has law $\mu_B$. It is $\sF_{e^{-N}}$-measurable, while the processes generating $\cG_N$ are independent of that sigma field. Here $Z_{-N}=e^{\frac{N}{2}}X_{e^{-N}}$, so the evolution from $\tau=e^{-N}$ to $1$ reads $X_1=F_N(\Xi(W))(Z_{-N})$. Hence
\begin{equation}\label{eq:cond-n}
 \Law(X_1\mid\cG_N)=\pi_N(\Xi(W))
       \quad\mP\mbox{-a.s.}
\end{equation}
The paths $\Xi_j(W)$, $j\leq0$, determine the complete Brownian path. Indeed,
\[
 W_{t_j(v)}-W_{t_j(0)}
       =\int_0^v t_j(r)^{\frac{1}{2}}\,\dd\bigl(\Xi_j(W)(r)\bigr),
\]
where the integral is defined by integration by parts. Summing the increments and using continuity at $W_0=0$, we get $\bigvee_N\cG_N=\sF_1^W$, up to completion. By \eqref{eq:blocks}, $\Xi(W)$ has law $Q_-$ under $\mP$. Applying Proposition~\ref{prop:sample} and L\'evy's upward theorem in \eqref{eq:cond-n} against a countable family of bounded continuous functions separating probability measures, we obtain \eqref{eq:cond}.

\emph{Strong existence and pathwise uniqueness.} For a strong solution, $X_1$ is $\sF_1^W$-measurable, so its conditional law is $\delta_{X_1}$. Formula \eqref{eq:cond} and $\Xi(W)\sim Q_-$ then force $\mu_{\underline\omega^-}$ to be Dirac $Q_-$-almost surely.

Conversely, suppose that $\mu_{\underline\omega^-}$ is a Dirac mass $Q_-$-almost surely. The map $x\mapsto\delta_x$ has a continuous inverse on the closed set of Dirac measures, so there is a Borel $a:\Omega_-\to\mR^d$ such that $\mu_{\underline\omega^-}=\delta_{a(\underline\omega^-)}$ $Q_-$-a.s.. For $0<\tau\leq1$, put $(R_\tau w)(t):=\tau^{-\frac{1}{2}}w(\tau t)$. The identity $\sqrt\tau\,b(\tau r,\sqrt\tau x)=b(r,x)$ shows that $(R_\tau X,R_\tau W)$ is again a weak solution. Applying \eqref{eq:cond} to it, we obtain
\begin{equation*}
 X_\tau=\sqrt\tau\,a\bigl(\Xi(R_\tau W)\bigr)
       \quad\mP\mbox{-a.s.}
\end{equation*}
The right-hand side depends only on $W$ up to $\tau$. Taking this identity at all rational times and using continuity shows that every weak solution is adapted to the Brownian filtration. The same representation makes two solutions driven by the same Brownian motion coincide, proving pathwise uniqueness.

By the Yamada--Watanabe theorem \cite[Chapter~IX, Theorem~1.7]{RY99continuous}, weak existence from Proposition~\ref{prop:entry} and the pathwise uniqueness just proved imply strong existence for \eqref{eq:scale} on any prescribed filtered probability space carrying a Brownian motion.
\end{proof}

Formula \eqref{eq:cond} is the link used below: the original SDE is strongly solvable precisely when the entire past noise determines the present state uniquely.

\section{The Lyapunov spectrum}\label{sec:lyap}

We now construct a family of profiles $B_K$, $d\geq2$. By Proposition~\ref{prop:entry}, weak existence and uniqueness in law hold for the corresponding self-similar equations. To rule out strong solvability for large $K$, we study the Lyapunov spectrum and show that its top exponent is positive.

\subsection{The profile}\label{sec:profile}

Let $d\geq2$ and write $z=(v,u)\in\mR^2\times\mR^{d-2}$, where $v=(v_1,v_2)$ and, for $d>2$, $u=(u_1,\ldots,u_{d-2})$. Put
\[
 J:=\begin{pmatrix}0&-1\\1&0\end{pmatrix}.
\]
Choose $\chi\in C_c^\infty([0,\infty))$ with $0\leq\chi\leq1$, $\chi=1$ on $[0,1]$, and $\chi=0$ on $[2,\infty)$, and set $H(r):=r\chi(r)$. For $d>2$, choose $\zeta\in C_c^\infty(\mR^{d-2})$ with $0\leq\zeta\leq1$ and $\zeta=1$ on $B_1(0)$. For $d=2$, omit $u$ and set $\zeta:=1$. For $K>0$, define
\begin{equation}\label{eq:shear}
B_K(z):=\bigl(K\zeta(u)H(|v|^2)Jv,0_{d-2}\bigr)
\end{equation} 
and
\begin{equation}\label{eq:shear_b}
b_K(t,x):=t^{-\frac{1}{2}}B_K(x/\sqrt t), ~~ 0<t\leq1; \quad b_K(0,x):=0. 
\end{equation}
The cutoffs make $B_K$ smooth and compactly supported, and $b_K$ is smooth on $(0,1]\times\mR^d$. Since $Jv\cdot v=0$ and $\tr J=0$,
\[
 \div_v\bigl(H(|v|^2)Jv\bigr)
       =2H'(|v|^2)v\cdot Jv+H(|v|^2)\tr J=0.
\]
The factor $\zeta(u)$ contributes no divergence, because the last $(d-2)$-components of $B_K$ vanish. Thus
\begin{equation}\label{eq:prof}
 B_K(z)\cdot z=0,\quad \div B_K(z)=0.
\end{equation}
For $0<t\leq1$, we obtain the first identity below by the change of variables $x=\sqrt t\,z$ and the other two estimates from the pointwise bound on $B_K$:
\begin{equation*}
 \|b_K(t)\|_{L^d}=\|B_K\|_{L^d},\quad
 \|b_K(t)\|_\infty\leq t^{-\frac{1}{2}}\|B_K\|_\infty,
 \quad
 \int_0^t\|b_K(r)\|_\infty\dd r\leq2\|B_K\|_\infty\sqrt t.
\end{equation*}
For $K>0$ and fixed $R>0$,
\begin{equation*}
 \lim_{t\downarrow0}
   \|b_K(t)\1_{\{|b_K(t)|>R\}}\|_d=\|B_K\|_d>0.
\end{equation*}
Thus, condition \eqref{eq:tails} fails. Moreover, restricting to a shorter initial interval does not reduce the $L^\infty_tL^d_x$ norm. In particular, $b_K$ does not belong to $C([0,1];L^d(\mR^d))$, the class discussed in Section~\ref{sec:results}.

For $d\geq3$, the Sobolev inequality nevertheles implies $b_K\in\bF_{C_d\|B_K\|_d^2}(0)$. Hence the large-$K$ examples also show that finite form bounds alone cannot ensure strong existence. In dimension two, these nonzero self-similar drifts belong to no $\bF_\delta(g)$ with finite $\delta$ and $g\in L^1$: scaling at the initial time forces a homogeneous form inequality, which logarithmic cutoffs exclude. The two main theorems do not determine an optimal form-bound threshold.

By Lemma~\ref{lem:harris}, \eqref{eq:prof} implies $\mu_{B_K}=\gamma_d$. For any $K>0$, take a weak solution of \eqref{eq:counter} with $b=b_K$ and denote its probability measure by $\mP$. By Proposition~\ref{prop:entry}, we have $\Law(Z_{\log t})=\gamma_d$; using the relation $X_t=\sqrt t Z_{\log t}$, we therefore obtain
\begin{equation}\label{eq:marg}
 \Law(X_t)(\dd x)=(2\pi t)^{-\frac{d}{2}}e^{-\frac{|x|^2}{2t}}\dd x,
 \quad 0<t\leq1.
\end{equation}
For $p\in(1,\infty)$ and $p':=p/(p-1)$, this density has $L^{p'}$ norm $C_{d,p}t^{-\frac{d}{2p}}$. Applying H\"older's inequality first in space and then in time, we obtain, for $p,q\in(1,\infty)$ with $d/p+2/q<2$,
\begin{equation}\label{eq:krylov}
 \bE\int_0^\tau |f(t,X_t)|\dd t
 \leq C_{d,p,q}\tau^{1-\frac{1}{q}-\frac{d}{2p}}
       \|f\|_{L^q((0,\tau);L^p(\mR^d))},\quad 0<\tau\leq1.
\end{equation}
Here $C_{d,p,q}$ depends only on $d,p,q$, and the condition $d/p+2/q<2$ ensures the required time integrability. Likewise, for $0<a<1$, we have by Cauchy--Schwarz
\[
 \bE\int_a^1 |f(t,X_t)|\dd t
 \leq C_d\left(\int_a^1 t^{-\frac{d}{2}}\dd t\right)^{\frac{1}{2}}
                  \|f\|_{L^2((a,1)\times\mR^d)}.
\]
Thus, every weak solution satisfies the local $L^2$ Krylov estimate in \eqref{eq:k2}, and hence belongs to $\sK_2(0,0)$.

The weak solution therefore has Gaussian one-time marginals and the Krylov estimates above. It remains to check whether it is a measurable function of its driving Brownian motion. By Proposition~\ref{prop:strong}, this is equivalent to asking whether the conditional law of $X_1$ given $W$ is almost surely a Dirac mass. For our profiles, Proposition \ref{prop:BY} below will allow us to show that this conditional law is non-Dirac almost surely once the top Lyapunov exponent is positive; we present this argument in Section~\ref{sec:srb}.

We therefore study the derivative of the flow with respect to its initial point, keeping the Brownian path fixed. Our choice of $B_K$ makes the angular speed depend on the radius: it equals $K|v|^2$ where $0<|v|<1$ and $\zeta(u)=1$. Points at different radii therefore rotate at different speeds. Without noise, $\dot z_t=-z_t/2+B_K(z_t)$ satisfies $|z_t|=e^{-\frac{t}{2}}|z_0|$, so the radial derivative of the angular velocity tends to zero as the trajectory approaches the origin. With Brownian noise, the invariant law is $\gamma_d$, and ergodicity ensures that trajectories spend a positive fraction of time in this region.

In the next subsection, we differentiate the solution with respect to the initial planar coordinate $v$ in a fixed direction. The resulting tangent vector changes both direction and length. We express the planar position in coordinates parallel and perpendicular to this vector. The two resulting coordinates, together with the remaining coordinates, satisfy an autonomous SDE. The logarithmic rate of growth of the vector's length is a function of these variables. The Lyapunov exponent is the average of this function against the invariant law of this SDE.

\begin{figure}[htbp]
 \centering
 \includegraphics[width=\textwidth]{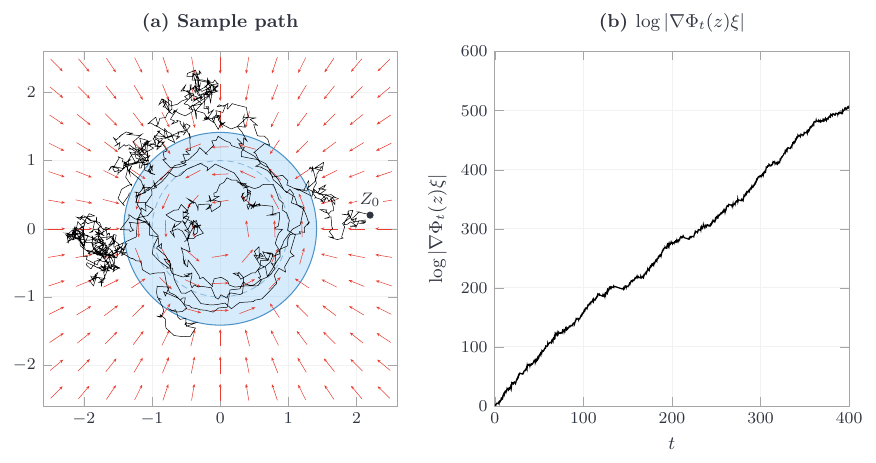}
\caption{\scriptsize Simulation for $d=2$, $K=20$, with $z=(2.2,0.2)$ and $\xi=(1,0)$. (a) A sample path $Z_t=\Phi_t(z)$ for $0\leq t\leq6$. (b) The logarithm of $|\nabla\Phi_t(z)\xi|$ for $0\leq t\leq400$.}
 \label{fig:profile}
\end{figure}

Figure~\ref{fig:profile} illustrates the planar trajectory and the growth of this tangent vector.

\Needspace{6\baselineskip}
\subsection{The auxiliary equation}\label{sec:frame}

For calculations on $[0,\infty)$, we work on $(\Omega_+,Q_+)$ from Section~\ref{sec:sample}. By concatenating the Brownian paths $\omega_j$, $j=1,2,\ldots$, we obtain a Brownian motion $\beta$ on $[0,\infty)$ under $Q_+$. Let $\Phi_t$ be the smooth flow of \eqref{eq:log} with $B=B_K$, driven by this Bxrownian motion and started at time zero. The coordinate subspace $E:=\mR^2\times\{0_{d-2} \}$ is invariant under $\nabla\Phi_t(z)$ for every $z$. We will calculate the growth of vectors in $E$.

Fix an initial point $z=(v,u)\in \mR^2\times \mR^{d-2}$ and a nonzero vector $\xi=(\xi_1,\xi_2)\in\mR^2$. Write the solution process
\[
Z_t:=\Phi_t(z)=(V_t,U_t)\in\mR^2\times\mR^{d-2}
\]
and the driving noise as $\beta=(\beta^v,\beta^u)$, where $\beta^v:=(\beta^1,\beta^2)$ and, for $d>2$, $\beta^u:=(\beta^3,\ldots,\beta^d)$. Let $\xi_t\in\mR^2$ be the first two components of $\nabla\Phi_t(z)(\xi,0_{d-2})$. To study the limit of $t^{-1}\log|\xi_t|$, we introduce the auxiliary process
\[
Y_t:=(a_t,c_t,U_t)\in \mR^d, 
\]
where
\[
 e_t:=\xi_t/|\xi_t|,\quad a_t:=V_t\cdot e_t,
 \quad c_t:=V_t\cdot Je_t. 
\]

\begin{lemma}\label{lem:frame}
There is a $d$-dimensional Brownian motion $\widehat\beta$ under $Q_+$, adapted to the natural filtration of $\beta$, such that
\begin{equation}\label{eq:aux}
 \d Y_t=\bigl(-\tfrac12Y_t+KF(Y_t)\bigr)\dd t+\dd\widehat\beta_t ~~\mbox{ and }~~
 \frac{\d}{\d t}\log|\xi_t|=G_K(Y_t), 
\end{equation}

where for $y=(a,c,u)\in \mR\times\mR\times\mR^{d-2}$,
\begin{equation*}
 \begin{aligned}
 F(y)&:=2\zeta(u)H'(a^2+c^2)a^2(c,-a,0_{d-2}),\\
 G_K(y)&:=-\frac12-2K\zeta(u)H'(a^2+c^2)ac=-\tfrac12-\tfrac K2\div F.
 \end{aligned}
\end{equation*}
Moreover, it holds that
\begin{equation}\label{eq:aux-id}
 y\cdot F(y)=0,\quad
 \div F(y)=4\zeta(u)H'(a^2+c^2)ac.
\end{equation}
\end{lemma}
\begin{proof}
For the profile $B_K$, the drift is
\[
 A(v,u)=\bigl(-\tfrac12v+K\zeta(u)H(|v|^2)Jv,
                   -\tfrac12u\bigr).
\]
Differentiating with respect to the initial point, we obtain
\begin{equation}\label{eq:dflow}
 \frac{\d}{\d t}\nabla\Phi_t(z)
       =\nabla A(\Phi_t(z))\nabla\Phi_t(z),\quad \nabla\Phi_0(z)=I_d.
\end{equation}
Writing $A=(A^v,A^u)$, the Jacobian matrix of the drift is
\begin{equation}\label{eq:DA}
 \begin{split}
 \nabla A(v,u)&=
 \begin{pmatrix}
 \nabla_v A^v(v,u)&KH(|v|^2)(Jv)(\nabla\zeta(u))^\top\\
 0&-\tfrac12I_{d-2}
 \end{pmatrix},\\
 \nabla_v A^v(v,u)&=-\tfrac12I_2
       +K\zeta(u)\bigl[H(|v|^2)J+2H'(|v|^2)(Jv)v^\top\bigr].
 \end{split}
\end{equation}
For $d=2$, the blocks involving $u$ are omitted. The zero lower-left block shows that $\nabla\Phi_t(z)(\xi,0_{d-2})$ remains in $E$. Equation~\eqref{eq:dflow} implies that $\xi_t$ solves the linear equation
\begin{equation}\label{eq:xi}
      \dot\xi_t=\nabla_v A^v(V_t,U_t)\xi_t,\quad \xi_0=\xi\neq 0.  
\end{equation}
We have $\xi_t\neq0$ for every $t$.

Let $v\in \mR^2$ and $e\in \mR^2$ with $|e|=1$. Put $a=v\cdot e$ and $c=v\cdot Je$. Then $v=ae+cJe$ and $Jv=-ce+aJe$. Applying \eqref{eq:DA} to $e$, we find
\begin{equation}\label{eq:DAv}
      \nabla_v A^v(v,u)e=\underbrace{\bigl(-\tfrac12-2K\zeta(u)H'(|v|^2)ac\bigr)}_{=\, G_K(y)} e+K\zeta(u)\bigl(H(|v|^2)+2H'(|v|^2)a^2\bigr)Je.
\end{equation}

Recall that $|V_t|^2=a_t^2+c_t^2$, and put
\begin{equation}\label{eq:vartheta}
      \vartheta_t:=K\zeta(U_t)[H(|V_t|^2)+2H'(|V_t|^2)a_t^2]. 
\end{equation}
Since $\xi_t=|\xi_t|e_t$, by \eqref{eq:xi} and \eqref{eq:DAv}, we get
\[
 \dot\xi_t=G_K(Y_t)\xi_t+\vartheta_tJ\xi_t=|\xi_t|\bigl(G_K(Y_t)e_t+\vartheta_tJe_t\bigr).
\]
Using $|e_t|=1$ and $e_t\cdot Je_t=0$, we obtain
\[
 \frac{\d}{\d t}|\xi_t|
 =\frac{\xi_t\cdot\dot\xi_t}{|\xi_t|}
 =|\xi_t|G_K(Y_t) ~~ \mbox{ and }~~
 \frac{\d}{\d t}\log|\xi_t|
 =\frac{1}{|\xi_t|}\frac{\d}{\d t}|\xi_t|=G_K(Y_t).
\]
Therefore, we verify the second equation in \eqref{eq:aux}.

Differentiating $e_t=\xi_t/|\xi_t|$, we therefore get
\begin{align*}
 \frac{\d e_t}{\d t} =\frac{\dot\xi_t}{|\xi_t|}
       -\frac{\xi_t}{|\xi_t|^2}\frac{\d}{\d t}|\xi_t| =G_K(Y_t)e_t+\vartheta_tJe_t-G_K(Y_t)e_t
 =\vartheta_tJe_t.
\end{align*}
Finally, since $J$ is constant and $J^2=-I_2$,
\[
 \frac{\d}{\d t}(Je_t)=J\dot e_t
 =\vartheta_tJ^2e_t=-\vartheta_te_t.
\]
In particular, $e_t$ has finite variation. Noting that
\[
\d V_t=\left[-\frac{1}{2}V_t+K\zeta(U_t)H(|V_t|^2)J V_t\right] \dd t+\dd\beta_t^v
\]
and $a_t=V_t\cdot e_t$, $c_t=V_t\cdot Je_t$, It\^o's formula and \eqref{eq:vartheta} imply that
\begin{align*}
 & \d a_t = e_t\cdot\d V_t+V_t\cdot\d e_t\\
 ={} & \bigl(-a_t/2-K\zeta(U_t)H(|V_t|^2)c_t+\vartheta_tc_t\bigr)\dd t
+e_t\cdot\dd\beta_t^v\\
={} & \bigl(-a_t/2+ \underbrace{2K\zeta(U_t)H'(|V_t|^2)a_t^2 c_t}_{=\, KF_1(Y_t)}\bigr)\dd t
+e_t\cdot\dd\beta_t^v,
\end{align*}
and
\begin{align*}
& \d c_t = Je_t\cdot\d V_t+V_t\cdot\d(Je_t)\\
={} & \bigl(-c_t/2+K\zeta(U_t)H(|V_t|^2)a_t-\vartheta_ta_t\bigr)\dd t
+Je_t\cdot\dd\beta_t^v\\
={} & \bigl(-c_t/2\underbrace{-2K\zeta(U_t)H'(|V_t|^2)a_t^3}_{=\, KF_2(Y_t)}\bigr)\dd t
+Je_t\cdot\dd\beta_t^v.
\end{align*}
Thus, $Y_t$ satisfies equation
\[
 \d Y_t=(-Y_t/2+KF(Y_t))\dd t+\dd\widehat{\beta}_t,
\]
with
\[
 \widehat\beta_t:=
 \left(\int_0^t e_s\cdot\dd\beta_s^v,
       \int_0^t Je_s\cdot\dd\beta_s^v,
       \beta_t^u\right).
\]
Since $\widehat{\beta}$ is a continuous local martingale under $Q_+$ with $\langle\widehat\beta\rangle_t=tI_d$, L\'evy's characterization shows that it is the Brownian motion in \eqref{eq:aux}.

Finally, a direct calculation shows that $F$ satisfies \eqref{eq:aux-id}.
\end{proof}

To describe the full Lyapunov spectrum, view $\Omega_+$ as $C_0([0,\infty);\mR^d)$ by concatenating the paths $\omega_j$, $j\geq1$. The Brownian shift and the associated skew product are
\[
 (\theta_t\omega)(r):=\omega(t+r)-\omega(t),\quad
 \Theta_t(z,\omega):=\bigl(\Phi_t(z,\omega),\theta_t\omega\bigr),
 \quad t,r\geq0.
\]
By the flow property, we have $\Theta_{t+s}=\Theta_t\circ\Theta_s$. Since $\gamma_d$ is stationary and future Brownian increments are independent of the past, $\Theta_t$ preserves $\gamma_d\times Q_+$; it is ergodic by Lemma~\ref{lem:harris}. Differentiating the flow property with respect to $z$, we obtain
\[
 \nabla\Phi_{t+s}(z,\omega)
 =\nabla\Phi_t\bigl(\Phi_s(z,\omega),\theta_s\omega\bigr)
       \nabla\Phi_s(z,\omega).
\]
Thus, $\nabla\Phi_t$ is a linear cocycle over $\Theta_t$. Since $\nabla A$ is bounded, \eqref{eq:dflow} yields
\[
 \sup_{0\leq t\leq1}
 \bigl(\|\nabla\Phi_t(z,\omega)\|\vee
       \|(\nabla\Phi_t(z,\omega))^{-1}\|\bigr)
       \leq e^{\|\nabla A\|_\infty}.
\]
By the continuous-time Oseledets theorem \cite[Theorem~3.4.1]{A98random}, this cocycle therefore has a Lyapunov spectrum under $\gamma_d\times Q_+$.

Combining the Oseledets theorem with Lemma~\ref{lem:frame} and the ergodic theorem for $Y_t$ (see Lemma \ref{lem:harris}), we express the top Lyapunov exponent as the stationary average of $G_K$ and obtain the full spectrum below.

\begin{proposition}\label{prop:aux}
For each $K>0$, the auxiliary diffusion $Y_t$ has a unique invariant probability $\mu_K(\dd y)=\rho_K(y)\dd y$, with smooth positive density. Under $\mu_K$, the pair $(a^2+c^2,u)$ has distribution
\begin{equation}\label{eq:fixed}
 \frac12e^{-\frac{r}{2}}\1_{(0,\infty)}(r)\dd r \times\gamma_{d-2}(\dd u),
\end{equation}
with the second factor omitted for $d=2$. For every initial point $z\in\mR^d$ and every $0\neq\xi\in\mR^2$,
\begin{equation}\label{eq:lyap}
 \lim_{t\to\infty}\frac1t\log|\xi_t|
       =\lambda_K:=\int G_K\dd\mu_K,
       \quad Q_+\mbox{-a.s.}.
\end{equation}
The Lyapunov exponents of $\nabla\Phi_t$, listed in nonincreasing order are
\begin{equation}\label{eq:spectrum}
 \lambda_K,\quad
 \underbrace{-\tfrac12,\ldots,-\tfrac12}_{d-2\ \mathrm{times}},
 \quad -1-\lambda_K,
 \quad (\gamma_d\times Q_+)\mbox{-a.s.},
\end{equation}
where the middle group is omitted for $d=2$.
\end{proposition}
\begin{proof}
The field $F$ is smooth, compactly supported, and tangent to spheres. By Lemma~\ref{lem:harris}, the invariant probability exists and the diffusion is positive Harris recurrent. To compute this distribution, put $R_t:=a_t^2+c_t^2$. It\^o's formula and the fact that $y\cdot F(y)=0$ imply that
\begin{align*}
 & \d R_t
  = 2a_t\dd a_t+2c_t\dd c_t
       +\dd\langle a\rangle_t+\dd\langle c\rangle_t\\
 ={} & \bigl(2-R_t+2K[a_tF_1(Y_t)+c_tF_2(Y_t)]\bigr)\dd t
       +2a_t\dd\widehat\beta_t^1+2c_t\dd\widehat\beta_t^2\\
 ={} & (2-R_t)\dd t
       +2a_t\dd\widehat\beta_t^1+2c_t\dd\widehat\beta_t^2. 
\end{align*}
Define
\[
 q_t:=\begin{cases}
 (a_t,c_t)/\sqrt{R_t},&R_t>0,\\
 (1,0),&R_t=0,
 \end{cases}
 \quad
 b_t:=\int_0^t q_s\cdot\dd(\widehat\beta_s^1,\widehat\beta_s^2).
\]
L\'evy's characterization shows that $(b,\widehat\beta^u)$ is a $(d-1)$-dimensional Brownian motion. Therefore,
\[
 \d R_t=(2-R_t)\dd t+2\sqrt{R_t}\dd b_t,
 \quad \d U_t=-U_t\dd t/2+\dd\widehat\beta_t^u.
\]
The equations are uncoupled and driven by independent Brownian motions. The equation for $R$ is also satisfied by the squared norm of a two-dimensional Ornstein--Uhlenbeck process. Therefore, the process $(R_t, U_t)$ has a unique invariant probability, which is the product of the law of $|N_2|^2$ and $N_{d-2}$, where $N_2\sim\gamma_2$ and $N_{d-2}\sim\gamma_{d-2}$ are independent. This proves \eqref{eq:fixed}.

The function $G_K$ is bounded for fixed $K$. By the ergodic theorem for the auxiliary diffusion $Y_t$ (Lemma \ref{lem:harris}), applied from its initial state, we have $t^{-1}\int_0^tG_K(Y_s)\dd s\to\int G_K\dd\mu_K$ $Q_+$-a.s.. This and \eqref{eq:aux} prove \eqref{eq:lyap}.

We now compute all Lyapunov exponents. Recall that $E=\mR^2\times\{0\}$, and write $\cH_t(z):=\nabla\Phi_t(z)|_E$ as a $2\times2$ matrix. Applying \eqref{eq:lyap} to the two standard basis vectors, we find $t^{-1}\log\|\cH_t(z)\|\to\lambda_K$ for each $z$, $Q_+$-a.s., and hence $(\gamma_d\times Q_+)$-a.s.. Since $\tr \nabla_v A^v=-1$, by Liouville's formula, we have $\det\cH_t(z)=e^{-t}$. In particular, $\|\cH_t(z)\|\geq e^{-\frac{t}{2}}$, so $\lambda_K\geq-1/2$. The product of its two singular values is therefore $e^{-t}$. Their exponential growth rates are $\lambda_K$ and $-1-\lambda_K$.

For $d>2$, the structure of $\nabla A$ given by \eqref{eq:DA} yields
\[
 \nabla\Phi_t(z)=
 \begin{pmatrix}
 \cH_t(z)& \cK_t(z)\\
 0&e^{-\frac{t}{2}}I_{d-2}
 \end{pmatrix}.
\]
The induced cocycle on the quotient $\mR^d/E$ is thus $e^{-\frac{t}{2}}I_{d-2}$, with all $d-2$ exponents equal to $-1/2$. The full spectrum is the union, with multiplicities, of the spectra on $E$ and $\mR^d/E$. Together with the exponents on $E$ computed above, this proves \eqref{eq:spectrum}.
\end{proof}

\subsection{Estimate for the top Lyapunov exponent}

Write
\[
 h_K:=\frac{\rho_K}{(2\pi)^{-\frac{d}{2}}e^{-\frac{|y|^2}{2}}},\quad
 I_K:=\int\frac{|\nabla h_K|^2}{h_K}\dd\gamma_d.
\]

The following lemma is crucial for proving the existence of a positive Lyapunov exponent for $\Phi_t(z)$.

\begin{lemma}\label{lem:fisher}
Suppose $F\in C_c^\infty(\mR^d;\mR^d)$, $F(y)\cdot y=0$, and the probability measure $h\,\dd\gamma_d$, with smooth positive density $h$ and finite second moment, is stationary for $L_0+KF\cdot\nabla$, where $L_0:=\frac12\Delta-\frac12y\cdot\nabla$. Then $h\in H^1(\gamma_d)$ and
\begin{equation}\label{eq:fisher}
 \frac12\int\frac{|\nabla h|^2}{h}\dd\gamma_d
       =-K\int\div F\,h\dd\gamma_d<\infty.
\end{equation}
\end{lemma}
\begin{proof}
For $\varphi\in C_c^\infty(\mR^d)$, stationarity of $h\gamma_d$ means
\begin{equation}\label{eq:station}
     0=\int hL_0\varphi\dd\gamma_d+K\int hF\cdot\nabla\varphi\dd\gamma_d. 
\end{equation}
Integrating by parts, we obtain
\[
 \int hL_0\varphi\dd\gamma_d
 =-\frac12\int\nabla h\cdot\nabla\varphi\dd\gamma_d.
\]
Hence,
\begin{equation}\label{eq:stat}
 \frac12\int\nabla h\cdot\nabla\varphi\dd\gamma_d
 =K\int Fh\cdot\nabla\varphi\dd\gamma_d,
 \quad \varphi\in C_c^\infty.
\end{equation}
We next prove $h\in H^1(\gamma_d)$ to justify using a logarithmic test function. Set
\[
 T_M(r):=\min\{r,M\},\quad
 J_M(r):=\int_0^rT_M(s)\dd s\leq Mr.
\]
Choose smooth radial cutoffs $\chi_R\uparrow1$, equal to one on $B_R$ and zero outside $B_{2R}$, with $|L_0\chi_R|\leq C\1_{\{|y|\geq R\}}$ for $R\geq1$, where $C$ depends only on $d$ and the fixed cutoff profile. Once $\mathrm{supp} F\subset B_R$, test \eqref{eq:stat} with $\chi_RT_M(h)$, using smooth approximations of the truncation. We have  $\div(F\gamma_d)=(\div F)\gamma_d$ and $F\cdot \nabla\chi_R=0$; since $h\nabla T_M(h)= \frac12 \nabla(T_M(h)^2)$, integration by parts leads to
\[
 \frac12\int\chi_R|\nabla T_M(h)|^2\dd\gamma_d
 =-\frac K2\int(\div F)T_M(h)^2\dd\gamma_d
       +\int J_M(h)L_0\chi_R\dd\gamma_d.
\]
The last term comes from $T_M(h)\nabla h=\nabla J_M(h)$ and tends to zero as $R\to\infty$, we therefore obtain
\[
 \int|\nabla T_M(h)|^2\dd\gamma_d
 \leq C_K:=K \int|\div F|h^2\dd\gamma_d<\infty,
\]
where $C_K$ depends only on $d,K,F$, since $h\gamma_d$ is the unique invariant probability by Lemma~\ref{lem:harris}; it is finite because $\div F$ is compactly supported and $h$ is smooth. Applying the Gaussian Poincar\'e inequality and using $\int T_M(h)\dd\gamma_d\leq1$, we obtain $\int T_M(h)^2\dd\gamma_d\leq1+C_K$. Letting $M\to\infty$ proves $h\in H^1(\gamma_d)$.

Equation \eqref{eq:stat} therefore extends to all test functions in $H^1(\gamma_d)$. For $\eps>0$, the function $\log(h+\eps)$ belongs to this space by the Lipschitz chain rule. Thus,
\begin{align*}
 \frac12\int\frac{|\nabla h|^2}{h+\eps}\dd\gamma_d=K\int\frac{h}{h+\eps}F\cdot\nabla h\dd\gamma_d \to K\int F\cdot\nabla h\dd\gamma_d
 =-K\int(\div F)h\dd\gamma_d.
\end{align*}
The limit on the right follows by dominated convergence on the compact support of $F$; the final equality again uses $F\cdot y=0$. Monotone convergence on the left as $\eps\downarrow0$ proves \eqref{eq:fisher}, including its finiteness.
\end{proof}

Combining Lemma~\ref{lem:fisher}, \eqref{eq:aux-id}, and \eqref{eq:lyap}, we obtain the exact identity
\begin{equation}\label{eq:lyap-fi}
\lambda_K=-\frac12+\frac14I_K.
\end{equation}
The integration-by-parts argument above is similar to the proof of the Fisher information formula for Lyapunov exponents in \cite[Section~3.1]{BBP22regularity}. Bedrossian, Blumenthal, and Punshon-Smith combined their formula with quantitative hypoelliptic estimates to prove positivity of the top Lyapunov exponent for the weakly damped, stochastically forced Lorenz--96 model; see Corollary D therein.

The term $-1/2$ in \eqref{eq:lyap-fi} is the contraction rate contributed by the Ornstein--Uhlenbeck drift. Here $h_K$ is the density, relative to $\gamma_d$, of the invariant law of $Y_t=(a_t,c_t,U_t)$. The coordinates $(a_t,c_t)$ express $V_t$ in the moving frame $(e_t,Je_t)$, so the Gaussian stationary law of $Z_t$ does not imply that $h_K$ is constant.

The reason why $I_K$ grows can be seen where $H'=\zeta=1$. Writing $a=r\cos\theta$ and $c=r\sin\theta$, the drift $KF$ has angular velocity $-2Kr^2\cos^2\theta$. It vanishes when $a=0$ ($\theta=\pi/2$ or $3\pi/2$), that is, when the tangent vector $\xi_t$ is perpendicular to the planar position $V_t$. A nonzero density invariant under this angular transport would have to be proportional to $1/\cos^2\theta$, which is not integrable. The proof below combines this observation with the fixed positive mass of an annular region to show that $I_K\to\infty$.

\begin{proposition}\label{prop:shear}
For the profile in \eqref{eq:shear}, it holds that
\[
 I_K\to\infty,
 \quad \lambda_K\to\infty
 \quad\mbox{as }K\to\infty.
\]
In particular, there exists $K_0<\infty$ such that $\lambda_K>0$ for every $K\geq K_0$. For each such $K$, on the two-sided stationary extension, there is a splitting into measurable invariant subspaces
\[
 \mR^d=E^u(\underline\omega,z)
       \oplus E^{(-\frac{1}{2})}(\underline\omega,z)
       \oplus E^{ss}(\underline\omega,z),
 \quad \mu_{B_K}^*\mbox{-a.s.},
\]
with respective Lyapunov exponents $\lambda_K,-1/2,-1-\lambda_K$ and dimensions $1,d-2,1$. Moreover,
\[
 E^u\oplus E^{ss}=\mR^2\times\{0\},
 \quad E^s:=E^{(-\frac{1}{2})}\oplus E^{ss},
 \quad \dim E^s=d-1.
\]
Here $E^u$ and $E^s$ are the unstable and stable subspaces, respectively.
\end{proposition}
\begin{proof}
Suppose, to the contrary, that $K_n\to\infty$ and $\sup_n I_{K_n}<\infty$. Write $a=r\cos\theta$, $c=r\sin\theta$, and choose $0<r_-<r_+<1$. On the region
\[
 \cA:=\{(a,c,u):r_-<\sqrt{a^2+c^2}<r_+,\ |u|<1\},
\]
we have $H'(r^2)=\zeta(u)=1$ and $F(y)=2a^2(c,-a,0)$. By \eqref{eq:fixed}, $m:=\mu_K(\cA)>0$ depends only on $d,r_-,r_+$. For $d=2$, omit the $u$ coordinate and its restriction.

Put $g_n:=\sqrt{h_{K_n}}$. By the assumed Fisher-information bound, we have
\[
 \int\bigl(g_n^2+|\nabla g_n|^2\bigr)\dd\gamma_d
 =1+\tfrac14 I_{K_n}\leq C.
\]
Here $C$ depends only on the assumed bound $\sup_n I_{K_n}$. On a fixed ball $B\supset\overline{\cA}$, the Gaussian density is bounded above and below by positive constants. By Rellich's theorem, after passing to a subsequence, we have $g_n\to g$ in $L^2(B)$, and hence $h_{K_n}\to h:=g^2$ in $L^1(B)$. In particular,
\begin{equation}\label{eq:mass}
 \int_{\cA}h\dd\gamma_d
 =\lim_{n\to\infty}\mu_{K_n}(\cA)=m\in (0,1].
\end{equation}

For $\varphi\in C_c^\infty(\cA)$, noting that $L_0\varphi$ is bounded, by \eqref{eq:station}, we obtain
\begin{equation}\label{eq:FDphi}
      \int_{\cA}F\cdot\nabla\varphi\,h\dd\gamma_d=-\lim_{n\to\infty}\frac1{K_n}\int L_0\varphi\dd\mu_{K_n}=0.
\end{equation}
On $\cA$, using the explicit formula for $F$, we obtain
\[
 F\cdot\nabla
 =2a^2(c\partial_a-a\partial_c)
 =-2r^2\cos^2\theta\,\partial_\theta.
\]
In polar coordinates, with $h$ and $\varphi$ expressed in $(r,\theta,u)$, equation \eqref{eq:FDphi} becomes
\[
 \int r^3e^{-\frac{r^2+|u|^2}{2}}\cos^2\theta\,h(r,\theta,u)
       \partial_\theta\varphi(r,\theta,u)\dd r\dd\theta\dd u=0.
\]
Here $\varphi$ is any smooth function that is $2\pi$-periodic in $\theta$ and compactly supported in $(r_-,r_+)\times B_1$ in the variables $(r,u)$. Since the weight $r^3e^{-\frac{r^2+|u|^2}{2}}$ is positive and independent of $\theta$, we obtain $\partial_\theta\bigl(\cos^2\theta\,h\bigr)=0$. Thus, for fixed $(r,u)$, the function $\cos^2\theta\,h(r,\theta,u)$ is constant on the circle. Therefore, $h(r,\theta,u)=C(r,u)/\cos^2\theta$. By \eqref{eq:mass},
\[
|\{(r,u)\in (r_-,r_+) \times B_1: C(r,u)>0\}|>0. 
\]
Noting that $1/\cos^2\theta$ is not integrable near $\theta=\pi/2$ or $3\pi/2$, Fubini's theorem then yields that
\[
\int_{\cA}h
=\int_{(r_-,r_+)\times B_1}rC(r,u)\dd r\dd u
 \int_0^{2\pi}\frac{\dd\theta}{\cos^2\theta}=\infty,
\]
which contradicts \eqref{eq:mass}. This proves $I_K\to\infty$. The assertions about $\lambda_K$ follow from \eqref{eq:lyap-fi}.

For $K\geq K_0$, applying the multiplicative ergodic theorem and Proposition~\ref{prop:aux}, we obtain the stated invariant splitting. The invariant subspace $E=\mR^2\times\{0\}$ has Lyapunov exponents $\lambda_K$ and $-1-\lambda_K$, so it equals $E^u\oplus E^{ss}$. The subspaces with negative exponents form the stable subspace $E^s$.
\end{proof}

\section{Random SRB measures and failure of strong solvability}\label{sec:srb}

For the profile $B_K$, Section~\ref{sec:sample} constructs the sample measure $\mu_{\underline\omega^-}$ by conditioning on the past Brownian noise. Proposition~\ref{prop:strong} shows that it is the conditional law of $X_1$ given $W$. Our remaining task is therefore to show that this sample measure is not a Dirac mass when $\lambda_K>0$.

A Dirac sample measure would mean that, after the entire past noise is fixed, evolution from the remote past concentrates the smooth stationary distribution at a single present state. A positive top Lyapunov exponent describes exponential growth of infinitesimal perturbations under the same noise, suggesting that the limiting distribution retains spatial spread along unstable directions.

Under the hypotheses verified below, the theorem of Ledrappier--Young \cite[Theorem~B]{LY88entropy}, with the direct geometric proof of Blumenthal--Young \cite{BY19equivalence}, establishes the following stronger conclusion: the sample measures are random SRB measures, whose conditional measures on local unstable manifolds are absolutely continuous with respect to the induced volume. We state this result for the Brownian paths in Section~\ref{sec:sample}.

\begin{proposition}[Ledrappier--Young, Blumenthal--Young]\label{prop:BY}
Let $M$ be a compact Riemannian manifold and let $\eta\mapsto f_\eta$, $\eta\in\Omega_0$, be a Borel family of $C^2$ diffeomorphisms of $M$. Let $(\omega_j)_{j\in\mZ}$ be independent Brownian paths on $[0,1]$, each with Wiener law $Q$. Assume that
\begin{equation}\label{eq:BY-c2}
 \int\log^+\|f_\eta\|_{C^2}\,Q(\dd\eta)
 +\int\log^+\|f_\eta^{-1}\|_{C^2}\,Q(\dd\eta)<\infty.
\end{equation}
Let $\mu$ be an ergodic stationary probability, with a continuous density relative to Riemannian volume. If the top Lyapunov exponent is positive $(\mu\times Q_+)$-a.s., then, for $Q^{\mZ}$-almost every $\underline\omega$, the weak limit $\mu_{\underline\omega^-}=\lim_{N\to\infty}(f_{\omega_0}\circ\cdots\circ f_{\omega_{-N+1}})_\#\mu$ exists and is a random SRB measure: its conditional measures on local unstable manifolds are absolutely continuous with respect to the induced Riemannian volume. In particular, $\mu_{\underline\omega^-}$ is not a Dirac mass.
\end{proposition}

The last assertion follows because a Dirac mass cannot have absolutely continuous conditional measures on positive-dimensional unstable manifolds.

To apply this theorem, we must first place our diffusion on a compact manifold. Our drift equals $-z/2$ outside a compact set, which allows the diffusion to extend smoothly across the point at infinity. The following lemma verifies the hypotheses of Proposition~\ref{prop:BY} and shows that the Lyapunov exponents are unchanged.

Let $\mS^d\subset\mR^{d+1}$ be the unit sphere, equipped with the Riemannian metric induced by the Euclidean metric on $\mR^{d+1}$, and let $N=(0,\ldots,0,1)$ be its north pole. The inverse stereographic projection is
\begin{equation}\label{eq:stereo}
 \iota:\mR^d\to \mS^d\setminus\{N\},\qquad
 \iota(z):=\left(\frac{2z}{1+|z|^2},
                \frac{|z|^2-1}{1+|z|^2}\right).
\end{equation}
Recall that $B_K\cdot z=\div B_K=0$ by \eqref{eq:prof}, so $\mu_{B_K}=\gamma_d$, due to Lemma~\ref{lem:harris}.

\begin{lemma}\label{lem:sphere}
For $Q$-a.s. $\eta$, the conjugates $\iota\circ f_\eta\circ\iota^{-1}$ extend to a measurable family of smooth diffeomorphisms $\widehat f_\eta$ of $\mS^d$ fixing $N$ satisfying \eqref{eq:BY-c2}. The measure $\widehat\gamma:=\iota_\#\gamma_d$ is an ergodic stationary probability with a smooth density relative to the Riemannian volume measure. The Lyapunov exponents with respect to $\widehat\gamma\times Q_+$ are the same as in \eqref{eq:spectrum}. In particular, its top exponent $\lambda_K>0$, if $K\geq K_0$.
\end{lemma}
\begin{proof}
Let $S:=\iota(0)$ be the south pole and let $\kappa:\mS^d\setminus\{S\}\to\mR^d$ be the stereographic chart obtained by projection from $S$. The two coordinate systems are related by
\[
 \kappa(N)=0,\quad \kappa(\iota(z))=I(z),\qquad
 I:\mR^d\setminus\{0\}\to\mR^d\setminus\{0\},
 \quad I(z):=\frac{z}{|z|^2}.
\]
Thus, $I$ is the coordinate change $z\mapsto y$, with $I^{-1}=I$; the point at infinity in the $z$-coordinates becomes $y=0$.

Choose $R_0>0$ such that $\mathrm{supp}\,B_K\subset B_{R_0}$. The Stratonovich vector fields of the original equation \eqref{eq:log} are $V_0(z)=A(z)=-z/2+B_K(z)$ and $V_i(z)=e_i$, $1\leq i\leq d$. In the $y$-coordinates they become $\widehat V_i(y)=\nabla I(I(y))V_i(I(y))$. Since
\[
 \nabla I(I(y))=|y|^2I_d-2yy^\top
\]
and $B_K(I(y))=0$ for $0<|y|<R_0^{-1}$, we obtain that
\[
 \widehat V_0(y)=\tfrac12y,\quad
 \widehat V_i(y)=|y|^2e_i-2y_i y,\quad 1\leq i\leq d, ~ 0<|y|<R_0^{-1}.
\]
These fields are smooth at $y=0$ and vanish there. Away from $N$ they are smooth by the original coordinates, so they define smooth vector fields on all of $\mS^d$. Their stochastic flow of smooth diffeomorphisms $\widehat\Phi_t$ fixes $N$ and satisfies $\widehat\Phi_t\circ\iota=\iota\circ\Phi_t$ by uniqueness. The time-one maps $\widehat f_\eta$ are therefore the required extensions. By smooth stochastic flow theory on compact manifolds \cite{K90stochastic}, these maps are measurable and satisfy the moment bound
\[
 \int\bigl(\|\widehat f_\eta\|_{C^2}^2
          +\|\widehat f_\eta^{-1}\|_{C^2}^2\bigr)Q(\dd\eta)\leq C_K<\infty,
\]
where $C_K$ depends only on $d,K$ and the fixed profiles $\chi,\zeta$. This implies \eqref{eq:BY-c2}.

The measure $\widehat\gamma=\iota_\#\gamma_d$ is stationary and ergodic, since $\widehat f_\eta\circ\iota=\iota\circ f_\eta$ and $\gamma_d$ is stationary and ergodic for the time-one maps by \eqref{eq:harris}.

The Riemannian volume in the $z$-coordinates is $2^d(1+|z|^2)^{-d}\dd z$, so
\[
 \frac{\dd\widehat\gamma}{\dd\mathrm{vol}}(\iota(z))
 =C_d(1+|z|^2)^d e^{-\frac{|z|^2}{2}}.
\]
Substituting $z=I(y)$, we express the density as $C_d(1+|y|^2)^d|y|^{-2d}e^{-\frac{1}{2|y|^2}}$, which extends smoothly by zero at $y=0$. Thus, $\widehat\gamma$ has a smooth density on $\mS^d$.

Finally, we compare the Lyapunov exponents of the time-one maps by tracking the growth of corresponding tangent vectors. We write $\d$ for the differential of a map between manifolds. Applying the chain rule to $\widehat\Phi_n\circ\iota=\iota\circ\Phi_n$, we obtain
\begin{equation}\label{eq:stereo-deriv}
 \d\widehat\Phi_n(\iota(z))\d\iota(z)v
 =\d\iota(\Phi_n(z))\nabla\Phi_n(z)v,
\end{equation}
where $\d\iota(z):\mR^d\to T_{\iota(z)}\mS^d$ is a linear isomorphism. By \eqref{eq:stereo}, we have
\begin{equation}\label{eq:stereo-norm}
 |\d\iota(x)w|=\frac{2|w|}{1+|x|^2},
 \quad x,w\in\mR^d.
\end{equation}
For $v\ne0$, set $u:=\d\iota(z)v$. Applying \eqref{eq:stereo-norm} to the right-hand side of \eqref{eq:stereo-deriv}, we get
\[
 \frac1n\log|\d\widehat\Phi_n(\iota(z))u|
 =\frac1n\log|\nabla\Phi_n(z)v|
   +\frac{\log2-\log(1+|\Phi_n(z)|^2)}{n}.
\]
It remains to show that the last term in the above equation goes to zero. Under $\gamma_d\times Q_+$, stationarity implies that $\Phi_n(\cdot)$ has law $\gamma_d$, so $\bE_{\gamma_d\times Q_+}|\Phi_n(\cdot)|^2=d$. Hence,
\[
 (\gamma_d\times Q_+)\bigl\{\log(1+|\Phi_n(z)|^2)> \sqrt{n} \bigr\}=(\gamma_d\times Q_+)\bigl\{1+|\Phi_n(z)|^2>e^{\sqrt{n}}\bigr\}
 \leq C e^{-\sqrt{n}}.
\]
Here $C$ depends only on $d$. The Borel--Cantelli lemma yields
\[
 \frac1n\log(1+|\Phi_n(z)|^2)\to 0,
 \quad (\gamma_d\times Q_+)\mbox{-a.s.}
\]
So we complete our proof.
\end{proof}

\begin{corollary}\label{cor:ndir}
If $\lambda_K>0$, then $\mu_{\underline\omega^-}$ is non-Dirac for $Q_-$-almost every realization $\underline\omega^-$ of the past noise.
\end{corollary}
\begin{proof}
By Lemma~\ref{lem:sphere}, Proposition~\ref{prop:BY} applies to the diffeomorphisms $\widehat f_\eta$ of $\mS^d$, with the same Brownian paths as in Section~\ref{sec:sample}. Their sample measures $\widehat\mu_{\underline\omega^-}$ are therefore almost surely non-Dirac.

To relate these measures to $\mu_{\underline\omega^-}$, put $\widehat F_N:=\widehat f_{\omega_0}\circ\cdots\circ \widehat f_{\omega_{-N+1}}$. For every $N$,
\[
 (\widehat F_N)_\#\widehat\gamma
       =\iota_\#\pi_N(\underline\omega^-).
\]
Taking weak limits and using \eqref{eq:sample} and continuity of $\iota$, we obtain
\begin{equation*}
 \widehat\mu_{\underline\omega^-}
       =\iota_\#\mu_{\underline\omega^-}
       \quad Q^{\mZ}\mbox{-a.s.}
\end{equation*}
Thus, the sample measure on $\mS^d$ is exactly the pushforward of the sample measure from Section~\ref{sec:sample}. Since $\iota$ is injective, $\mu_{\underline\omega^-}$ must be non-Dirac as well. Since it depends only on the past, the stated $Q_-$-almost sure conclusion follows.
\end{proof}

\begin{proof}[Proof of Theorem~\ref{thm:counter}]
Choose $K\geq K_0$ so that $\lambda_K>0$, as given by Proposition~\ref{prop:shear}, and set $b:=b_K$ as in \eqref{eq:shear_b}. In Section~\ref{sec:profile}, we established the $L^\infty_tL^d_x$ norm bound and the divergence-free property. By Proposition~\ref{prop:entry}, weak existence and uniqueness in law hold on $C([0,1];\mR^d)$, with the Gaussian one-time marginals in \eqref{eq:marg}.

For any weak solution, $\Xi(W)$ has law $Q_-$. Combining Corollary~\ref{cor:ndir} and \eqref{eq:cond}, we therefore obtain
\[
 \Law(X_1\mid\sF_1^W)=\mu_{\Xi(W)},
       \quad \mu_{\Xi(W)}\mbox{ is non-Dirac almost surely}.
\]
Thus, $X_1$ cannot be a measurable function of the Brownian path. Proposition~\ref{prop:strong} excludes strong existence on $[0,1]$. Since weak existence and pathwise uniqueness would imply strong existence by the Yamada--Watanabe theorem, pathwise uniqueness fails on $[0,1]$ as well.
\end{proof}

\appendix
\renewcommand{\theequation}{\thesection.\arabic{equation}}
\section{Some parabolic estimates}\label{app:pde}

We use the notation of Section~\ref{sec:pde}. In particular, $a$ is bounded and uniformly elliptic, and $c,\nabla a,\widehat b$ belong to $\bF_{C\delta}(Cg)$. All estimates below are uniform under admissible approximation and the change of variables $r\mapsto t-r$. For the weight exponent $\kappa=2q$ used in the applications below and in Part~\ref{part:well}, the required smallness of $\delta$ is included in the choice of $\overline\delta$ in Theorem~\ref{thm:flow}.

\subsection{Energy and $L^\infty$ estimates}\label{app:energy}

\begin{proof}[Proof of Lemma~\ref{lem:pde}]
We first assume that $F$ is smooth in space. Test the equation with $q|u|^{q-2}u\rho_z$ and its $k$th derivative with $q|\nabla u|^{q-2}\partial_ku\rho_z$, and sum over $k$. By integration by parts, ellipticity and Young's inequality, we obtain
\begin{align}
 &\frac{\d}{\d t}\int(|u|^q+|\nabla u|^q)\rho_z
 +c_0\int\bigl(|u|^{q-2}|\nabla u|^2 +|\nabla u|^{q-2}|\nabla^2u|^2\bigr)\rho_z
 \nonumber\\
 \leq{} & C\int(1+|c|^2+|\nabla a|^2) (|u|^q+|\nabla u|^q)\rho_z+C\int|F|^2\bigl(|u|^{q-2}+|\nabla u|^{q-2}\bigr)\rho_z.
 \label{eq:coerc}
\end{align}
Here $c_0>0$ depends only on $d,\lambda$, and $C$ depends only on $d,\lambda,\kappa$. The weight terms are controlled by $|\nabla\rho_z|\leq C\rho_z$. These tests are justified by spatial cutoffs and regularization of the powers at zero.

For $v=u$ or $v=\nabla u$, we have $|\nabla|v|^{\frac{q}{2}}|^2\leq(q^2/4)|v|^{q-2}|\nabla v|^2$. Choose $\overline\delta=\overline\delta(d,\lambda,\kappa)>0$ small enough that, for $\delta<\overline\delta$, Lemma~\ref{lem:wform} absorbs the coefficient terms in \eqref{eq:coerc} into one half of the coercive term, leaving $C(1+g)\int(|u|^q+|\nabla u|^q)\rho_z$. If $q>2$, use
\begin{equation*}
 |F|^2|v|^{q-2}\leq h|F|^2|v|^q+C_h|F|^2,
 \quad h>0.
\end{equation*}
Applying Lemma~\ref{lem:wform} to $|v|^{\frac{q}{2}}$ and choosing $h$ sufficiently small in terms of $d,\lambda,\kappa,\beta$ absorbs the gradient term from the source into one half of the remaining coercivity. Its remaining contribution is bounded by $C(1+g_F)\int|v|^q\rho_z+C\int|F|^2\rho_z$, where $C$ depends only on $d,\lambda,\kappa,\beta$. For $q=2$, the source term is already bounded by $C\int|F|^2\rho_z$. Applying Gronwall's inequality, we obtain \eqref{eq:energy}.

When $\varphi=0$, the local Sobolev embedding, $q>d$, and \eqref{eq:energy} imply
\begin{align*}
 & \|u\|_{L^\infty([0,S]\times\mR^d)}^q
  \leq C\sup_z\sup_{t\leq S}
                \int(|u(t)|^q+|\nabla u(t)|^q)\rho_z\\
 \leq{} & C\sup_z\int_0^S\int|F|^2\rho_z
 \leq C\sup_z\int_0^S\int_{B_2(z)}|F|^2.
\end{align*}
The last inequality follows by a covering with unit cubes and $\kappa>d$. This proves \eqref{eq:pot}. By truncating and mollifying in space, we extend the estimates to $F\in\bF_\beta(g_F)$.
\end{proof}

\subsection{Transition semigroups and the backward equation}\label{app:semi}

\begin{proof}[Proof of Proposition~\ref{prop:semi}]
We choose $\overline\delta=\overline\delta(d,\lambda)>0$ small enough to satisfy the energy estimates in all three steps below.

\emph{(i) Convergence of the transition semigroups.}
Fix $\kappa=2q$. The functions $u_n$ solve $\partial_su_n+ \nabla\cdot(a^n \nabla u^n) + c^n \cdot \nabla u^n =0$, $u_n(t)=\varphi$. By Lemma~\ref{lem:pde} and the maximum principle,
\begin{equation}\label{eq:semi-u}
 \sup_n\sup_{s\leq t}\|u_n(s)\|_{W^{1,q}_{\rho_z}}<\infty,
 \quad \|u_n\|_\infty\leq\|\varphi\|_\infty.
\end{equation}
The first bound depends only on $d,\lambda,T$, $\|g\|_{L^1(0,T)}$ and $\|\varphi\|_{W^{1,q}}$. Fix $R>0$. Since $q>d$, we deduce from \eqref{eq:semi-u} that $\sup_n \sup_{s\leq t} \|u_n(s)\|_{C^\gamma(B_R)}<\infty$ with $\gamma=1/2-d/(2q)$.

It remains to establish time equicontinuity. We first estimate the time derivative in $H^{-1}(B_R)$. Since $q=d+1$, we have $L^{\frac{2q}{q+2}}(B_R)\hookrightarrow H^{-1}(B_R)$ by the Sobolev embedding. Hence 
\begin{align*}
 & \|\partial_su_n(s)\|_{H^{-1}(B_R)}
  \leq C_R\|a^n(s)\nabla u_n(s)\|_{L^2(B_R)}
       +C_R\|c^n(s)\|_{L^2(B_R)}\|\nabla u_n(s)\|_{L^q(B_R)}\\
 \leq{} & C_R\bigl(1+\sqrt{g(s)}\bigr)\quad\text{for a.e. }s.
\end{align*}
The last inequality uses \eqref{eq:semi-u}, boundedness of $a^n$ and the form bound for $c^n$, and $C_R$ depends only on $d,\lambda,R,T$, $\|g\|_{L^1(0,T)}$ and $\|\varphi\|_{W^{1,q}}$. Since $g\in L^1(0,T)$, we therefore have
\[
 \sup_n\|\partial_su_n\|_{L^2(0,t;H^{-1}(B_R))}\leq C_R.
\]
By the Cauchy--Schwarz inequality, for $0\leq r\leq s\leq t$,
\[
 \|u_n(s)-u_n(r)\|_{H^{-1}(B_R)}
 \leq |s-r|^{\frac{1}{2}}\|\partial_su_n\|_{L^2(0,t;H^{-1}(B_R))}
 \leq C_R|s-r|^{\frac{1}{2}}.
\]
Noting that $W^{1,q}(B_R)\subset C^{2\gamma}(\overline B_R)\hookrightarrow H^{-1}(B_R)$, and using \eqref{eq:semi-u}, we have 
\begin{align*}
 & \|u_n(s)-u_n(r)\|_{C^{\gamma}(\overline B_R)}
  \leq \eps\|u_n(s)-u_n(r)\|_{W^{1,q}(B_R)}
       +C_{R,\eps}\|u_n(s)-u_n(r)\|_{H^{-1}(B_R)}\\
 \leq{} & C_R\eps+C_{R,\eps}|s-r|^{\frac{1}{2}},
\end{align*}
Choosing first $\eps$ and then $|s-r|$ small proves time equicontinuity, uniformly in $n$ and $x\in\overline B_R$.

Together with the spatial equicontinuity established above, this establishes joint equicontinuity on $[0,t]\times\overline B_R$. By the uniform bound in \eqref{eq:semi-u}, the Arzel\`a--Ascoli theorem and a diagonal argument over $R$, we obtain relative compactness for locally uniform convergence on $[0,t]\times\mR^d$.

Every limit $u$ satisfies
\begin{equation}\label{eq:back}
 \partial_su+\nabla\cdot(a\nabla u)+c\cdot\nabla u=0,
 \quad u(t)=\varphi,
\end{equation}
in distributions, with the bounds \eqref{eq:semi-u}. This follows from the strong local $L^2$ convergence of $c^n$, weak local $L^2$ convergence of $\nabla u_n$, and uniform boundedness of $a^n$ together with its local convergence to $a$ in measure. For the difference $v$ of two limits, after reversing time, we test with $v\rho$ and apply Lemma~\ref{lem:wform} to obtain
\[
 \frac{\d}{\d r}\|v(r)\|_{2,\rho}^2
 +c_1\|\nabla v(r)\|_{2,\rho}^2
 \leq C(1+g(t-r))\|v(r)\|_{2,\rho}^2,
 \quad v(0)=0.
\]
Here $c_1>0$ and $C$ depend only on $d,\lambda$. Time regularization and spatial cutoffs justify this test. Gronwall's inequality implies $v=0$. Hence the whole sequence converges, and the limit is independent of the approximation.

By approximation on balls and the moment bound \eqref{eq:disp-u}, the convergence extends to bounded continuous terminal functions. By the same argument, we also have $P^n_{r,t}\varphi_n\to P_{r,t}\varphi$ locally uniformly whenever $\varphi_n\in C_b(\mR^d)$ are uniformly bounded and converge locally uniformly to $\varphi$. The limiting operators define probability kernels by tightness of the transition probabilities and continuity of $P_{r,t}\varphi$. For $s<t_1<\cdots<t_k$ and $\varphi_j\in C_b(\mR^d)$, the Markov property expresses $\bE\prod_{j=1}^k\varphi_j(X^n_{s,t_j}(x))$ by successive application of $P^n_{t_i,t_{i+1}}$. Passing to the limit from $t_k$ to $s$ shows that all subsequential limits of the laws on path space have the same finite-dimensional distributions. Hence the law on $C([s,T];\mR^d)$ is unique.

\emph{(ii) Smoothing and the local $L^2$ Krylov estimate.}
Let $S:=t-r$, $v(\tau):=P^n_{t-\tau,t}\varphi$ for $0\leq\tau\leq S$, and $h(\tau):=g(t-\tau)$. For $p\geq2$, test the equation with $p|v|^{p-2}v$. Writing $w:=|v|^{\frac{p}{2}}$, we obtain
\[
 \frac{\d}{\d\tau}\|v\|_p^p
 +\frac{4(p-1)}p\int a^n\nabla w\cdot\nabla w
 =2\int c^n\cdot\nabla w\,w.
\]
Since $4(p-1)/p\geq2$, we apply Young's inequality and the form bound to obtain
\begin{equation*}
 \frac{\d}{\d\tau}\|v\|_p^p
 +c_1\|\nabla|v|^{\frac{p}{2}}\|_2^2
 \leq A h(\tau)\|v\|_p^p,
\end{equation*}
where $c_1,A>0$ and the smallness threshold for $\delta$ depend only on $d,\lambda$. Put
\[
 V(\tau):=\exp\left(-\frac A2\int_0^\tau h(\theta)\dd\theta\right)v(\tau).
\]
Then, for every $p\geq2$,
\begin{equation}\label{eq:damp}
 \frac{\d}{\d\tau}\|V\|_p^p
 +c_1\|\nabla|V|^{\frac{p}{2}}\|_2^2\leq0.
\end{equation}
Apply the Nash inequality
\[
 \|w\|_2^{2+\frac{4}{d}}\leq C_d\|\nabla w\|_2^2\|w\|_1^{\frac{4}{d}}
\]
to $w=|V|^p$. Using \eqref{eq:damp} at exponents $p$ and $2p$ and integrating the resulting differential inequality, we obtain
\[
 \|V(\tau)\|_{2p}
 \leq C^{\frac{1}{2p}}(\tau-\sigma)^{-\frac{d}{4p}}\|V(\sigma)\|_p,
 \quad 0\leq\sigma<\tau\leq S.
\]
Here $C$ depends only on $d,\lambda$. Iterate with $p_k:=2^{k+1}$ and $\tau_k:=S(1-2^{-k})$. Since $\sum_{k\geq0}p_k^{-1}=1$ and $\sum_{k\geq0}k/p_k<\infty$, the $L^p$ contraction in \eqref{eq:damp} yields
\[
 \|V(S)\|_\infty\leq CS^{-\frac{d}{4}}\|\varphi\|_2.
\]
Returning to $v$ proves the estimate for $P^n$, with $C$ depending only on $d,\lambda$ and $\|g\|_{L^1(0,T)}$. By the convergence in part~(i), this estimate passes to $P$.

Let $Z$ be the canonical process under $\mathsf Q$. For $t>s$,
\[
 \bigl|\bE_{\mathsf Q}\varphi(Z_t)\bigr|
 =|P_{s,t}\varphi(x)|
 \leq C(t-s)^{-\frac{d}{4}}\|\varphi\|_2.
\]
By the Riesz representation theorem, the marginal law has a density satisfying \eqref{eq:dens}. For continuous compactly supported $f$, Fubini's theorem and Cauchy--Schwarz imply
\[
 \bE_{\mathsf Q}\int_a^T|f(r,Z_r)|\dd r
 \leq C\left(\int_a^T(r-s)^{-\frac{d}{2}}\dd r\right)^{\frac{1}{2}}
          \|f\|_{L^2((a,T)\times\mR^d)},
 \quad a>s.
\]
By marginal absolute continuity, $L^2$-approximation and Fatou's lemma, we extend the same inequality to every Borel $f\in L^2$. The time integral is finite, so every weak solution with path law $\mathsf Q$ belongs to $\sK_2(s,x)$.

\emph{(iii) Regularity for the backward equation.}
Set $v(\tau):=u_n(t-\tau)$. Testing its equation with $v$ and $-\Delta v$, we obtain
\begin{align*}
 \frac{\d}{\d\tau}\|v\|_2^2+c_0\|\nabla v\|_2^2
 &\leq Cg(t-\tau)\|v\|_2^2,\\
 \frac{\d}{\d\tau}\|\nabla v\|_2^2+c_0\|\nabla^2v\|_2^2
 &\leq Cg(t-\tau)\|\nabla v\|_2^2.
\end{align*}
Here $c_0>0$ and $C$ depend only on $d,\lambda$. For the second estimate, the error terms are bounded by
\[
 C\int(|\nabla a^n|+|c^n|)|\nabla v|\,|\nabla^2v|
 \leq (\eps+C_\eps\delta)\|\nabla^2v\|_2^2
       +C_\eps g(t-\tau)\|\nabla v\|_2^2,
\]
by Young's inequality and the form bounds applied to $\nabla v$. Choose $\eps=\eps(d,\lambda)>0$ and decrease $\overline\delta$ if necessary so that these terms can be absorbed whenever $\delta<\overline\delta$. The tests are justified by spatial cutoffs. Gronwall's inequality yields the $L^\infty H^1$ and $L^2H^2$ bounds. Moreover,
\[
 \int_0^t\|\widehat b^n\cdot\nabla u_n\|_2^2\dd r
 \leq C\delta\int_0^t\|\nabla^2u_n\|_2^2\dd r
       +C\int_0^t g(r)\|\nabla u_n\|_2^2\dd r.
\]
Since $a^n$ is bounded, the equation in nondivergence form also bounds $\partial_ru_n$ in $L^2$. This proves \eqref{eq:l2back} uniformly in $n$.

By weak compactness and part~(i), we obtain the same bounds for $u$. The form bounds for $\nabla a$ and $\widehat b$ allow us to write \eqref{eq:back} as $\partial_ru+a:\nabla^2u+\widehat b\cdot\nabla u=0$ in space-time $L^2$. The bounded continuous representative and terminal value follow from part~(i), while $L^2(0,t;H^2)\cap W^{1,2}(0,t;L^2)\hookrightarrow C([0,t];H^1)$ justifies the supremum in \eqref{eq:l2back}.
\end{proof}

\section*{Acknowledgments}
The authors thank Damir Kinzebulatov for informing them of his recent work \cite{BK26strong}, which motivated them to revisit and complete a manuscript that had been left unfinished for some time. The second author also thanks Damir for many helpful discussions on analysis.

\medskip

\paragraph{\bf AI Statement:} ChatGPT (Sol and Astra within Codex) assisted with language editing, organization, bibliographic checks, and checking calculations, based on a draft prepared entirely by the authors containing a complete proof. The authors take full responsibility for the manuscript in its present form.

\small
\bibliographystyle{alphaabbrv}
\bibliography{mybib.bib}

\end{document}